\documentclass[a4paper,leqno]{article}

\usepackage{amsmath,amssymb,amsthm,ifthen,verbatim}

\newcounter{mylisti} \newcounter{mylistii}
\newcounter{nest}
\newcommand{\defaultlabel}{}
\newlength{\toowide}

\newenvironment{mylist}[1]{%
  \addtocounter{nest}{1}
  \ifthenelse{\value{nest}=1}{%
    \renewcommand{\defaultlabel}{(\roman{mylisti})\hfill}}{%
    \renewcommand{\defaultlabel}{(\alph{mylistii})\hfill}}
  \settowidth{\toowide}{(iii)}
  \begin{list}{\defaultlabel}{%
      \ifthenelse{\value{nest}=1}{\usecounter{mylisti}}{%
        \usecounter{mylistii}}
      
      \addtolength{\topsep}{1ex}
      \addtolength{\itemsep}{0.5ex}
      \settowidth{\labelwidth}{#1}
      \ifthenelse{\lengthtest{\labelwidth > \toowide}}{%
        \setlength{\labelwidth}{\toowide}}{}
      \setlength{\leftmargin}{\labelwidth}
      \addtolength{\leftmargin}{\labelsep}}}{\addtocounter{nest}{-1}
\end{list}}

\newenvironment{myeqnarray}{%
  \begin{list}{}{%
      \setlength{\partopsep}{.6ex}
      \setlength{\labelsep}{1em}
      \settowidth{\labelwidth}{{\rmfamily \upshape (99)}}
      \setlength{\leftmargin}{\labelwidth}
      \addtolength{\leftmargin}{\labelsep}
      
      \renewcommand{\\}[1][0pt]{\end{math} \vspace{##1}
  \refstepcounter{equation} \item \begin{math} \ds}}%
    \refstepcounter{equation} \item \begin{math} \ds}{%
\end{math} \end{list}}

\theoremstyle{plain}
\newtheorem{thm}{Theorem}
\newtheorem{lemma}[thm]{Lemma}
\newtheorem{prop}[thm]{Proposition}
\newtheorem{cor}[thm]{Corollary}
\newtheorem{defn}[thm]{Definition}
\newtheorem{problem}[thm]{Problem}

\theoremstyle{definition}
\newtheorem{ex}[thm]{Example}

\newtheorem{rem}{Remark}
\newtheorem{rems}{Remarks}
\newtheorem{claim}{Claim}

\newcommand{\bn}{{\mathbb N}}

\newcommand{\br}{{\mathbb R}}
\newcommand{\ba}{\ensuremath{\boldsymbol{a}}}
\newcommand{\bb}{\ensuremath{\boldsymbol{b}}}
\newcommand{\bc}{\ensuremath{\boldsymbol{c}}}
\newcommand{\bi}{\ensuremath{\mathbb{I}}}

\newcommand{\vegtelen}[1]{{#1}^{(\omega)}}
\newcommand{\veges}[1]{{#1}^{(<\omega)}}
\newcommand{\supp}[1]{\ensuremath{\mathrm{supp}(#1)}}

\newcommand{\osc}[2]{\ensuremath{\mathrm{osc}(#1,#2)}}

\newcommand{\kin}{\!\in\!}
\newcommand{\knotin}{\!\notin\!}
\newcommand{\keq}{\!=\!}
\newcommand{\kneq}{\!\neq\!}
\newcommand{\kge}{\!>\!}
\newcommand{\kle}{\!<\!}
\newcommand{\kgeq}{\!\geq\!}
\newcommand{\kleq}{\!\leq\!}
\newcommand{\kplus}{\!+\!}
\newcommand{\kminus}{\!-\!}
\newcommand{\ksubset}{\!\subset\!}
\newcommand{\ksupset}{\!\supset\!}
\newcommand{\kprec}{\!\prec\!}

\newcommand{\ksetminus}{\!\setminus\!}

\newcommand{\norm}[1]{\lVert #1\rVert}
\newcommand{\tnorm}[1]{\lvert\mspace{-1mu}\lvert\mspace{-1mu}\lvert
  #1\rvert\mspace{-1mu}\rvert\mspace{-1mu}\rvert}

\newcommand{\Btnorm}[1]{\Big\lvert\mspace{-1mu}\Big\lvert\mspace{-1mu}%
  \Big\lvert #1\Big\rvert\mspace{-1mu}\Big\rvert\mspace{-1mu}%
  \Big\rvert}
\newcommand{\bnorm}[1]{\big\lVert #1\big\rVert}
\newcommand{\Bnorm}[1]{\Big\lVert #1\Big\rVert}
\newcommand{\abs}[1]{\lvert #1\rvert}
\newcommand{\babs}[1]{\big\lvert #1\big\rvert}
\newcommand{\Babs}[1]{\Big\lvert #1\Big\rvert}

\newcommand{\mc}{\mathrm{c}}

\newcommand{\ds}{\displaystyle}
\newcommand{\phtm}[1]{\text{\makebox[0pt]{\phantom{$#1$}}}}
\newcommand{\eg}{\text{e.g.}\ }
\newcommand{\ie}{\text{i.e.}\ }
\newcommand{\cf}{\text{c.f.}\ }

\title{Partial Unconditionality}
\author{\noindent S.~J.~Dilworth \and
  E.~Odell\footnotemark[2]\setcounter{footnote}{1} \and
  Th.~Schlumprecht\thanks{Research of the second and third named
    authors was supported by NSF} \and A.~Zs\' ak}
\date{\today\\[16pt]\textit{Dedicated to Haskell Rosenthal on
    the occasion of his $65^{\text{th}}$ birthday}}

\begin{document}

\maketitle

\begin{abstract}
  J.~Elton proved that for $\delta\kin (0,1]$ there exists
  $K(\delta)\kle\infty$ such that every normalized weakly null
  sequence in a Banach space admits a subsequence $(x_i)$ with the
  following property: if $a_i\kin [-1,1]$ for all $i\kin \bn$ and
  $E\ksubset \{i\kin\bn:\,\abs{a_i}\kgeq\delta\}$, then
  $\norm{\sum_{i\in E}a_ix_i}\kleq K(\delta)\norm{\sum_{i}a_ix_i}$. It
  is unknown if $\sup_{\delta >0}K(\delta)\kle\infty$. This problem
  turns out to be closely related to the question whether every
  infinite-dimensional Banach space contains a quasi-greedy basic
  sequence. The notion of a quasi-greedy basic sequence was introduced
  recently by S.~V.~Konyagin and V.~N.~Temlyakov. We present an
  extension of Elton's result which includes Schreier
  unconditionality. The proof involves a basic framework which we show
  can be also employed to prove other partial unconditionality results
  including that of convex unconditionality due to Argyros,
  Mercourakis and Tsarpalias. Various constants of partial
  unconditionality are defined and we investigate the relationships
  between them. We also explore the combinatorial problem underlying
  the $\sup_{\delta >0}K(\delta)\kle\infty$ problem and show that
  $\sup_{\delta >0}K(\delta)\kgeq 5/4$.
\end{abstract}

\tableofcontents

\section{Introduction}

Given a weakly null, normalized sequence in a Banach space, can we
pass to a subsequence that is a basic sequence and is in some sense
close to being unconditional? There are various ways in which one
can make this vague question precise, and in many situations one has a
positive answer. There are important cases, however, for which the
corresponding question is still open. In this paper we will study such
questions and provide some partial answers. We will also revisit known
results and discuss the relationship (\eg duality) between the various
notions of partial unconditionality.

As usual, we denote by $\mc _{00}$ the space of scalar sequences that
are eventually zero. Given a basic sequence $(x_i)$ in a Banach space
and $\delta\kin(0,1]$, we say $(x_i)$ is
\emph{$\delta$-near-unconditional with constant $C$} if its basis
constant is at most $C$ and
\begin{equation}
  \label{eq:near-unc}
  \Big\lVert \sum _{i\in E} a_ix_i \Big\rVert \leq C \Big\lVert
  \sum_{i=1}^{\infty} a_ix_i \Big\rVert
\end{equation}
for all $(a_i)\kin\mc_{00}$ with $\abs{a_i}\kleq 1$ for all
$i\kin\bn$, and for all $E\ksubset
\{i\kin\bn:\,\abs{a_i}\kgeq\delta\}$. Roughly speaking, this says that
we are allowed to project vectors onto sets of co-ordinates with
``large'' coefficients. A basic sequence is called
\emph{$\delta$-near-unconditional} if for some $C$ it is
$\delta$-near-unconditional with constant $C$; it is called
\emph{near-unconditional} if it is $\delta$-near-unconditional for all
$\delta\kin(0,1]$. The following result is due to J.~Elton.
\begin{thm}[Elton~\cite{El}]
  \label{thm:near-unc}
  For each $\delta\kin(0,1]$, every normalized, weakly null sequence
  has a $\delta$-near-unconditional subsequence. In particular, every
  normalized, weakly null sequence has a near-unconditional
  subsequence.
\end{thm}
For each $\delta\kin (0,1]$ let $K(\delta)$ be the infimum of the set
of real numbers~$K$ such that every normalized, weakly null sequence
has a $\delta$-near-unconditional subsequence with
constant~$K$. An upper bound of order $\log\big(1/\delta\big)$ for
$K(\delta)$ follows from the proof of Theorem~\ref{thm:near-unc}
presented in~\cite{O1}. This was first pointed out by Dilworth, Kalton
and Kutzarova~\cite{DKK}. It is unknown whether there is in fact a
uniform upper bound.
\begin{problem}
  \label{problem:near-unc}
  Let $K$ be the function defined above. Is
  $\sup_{\delta>0}K(\delta)\kle\infty?$
\end{problem}
Additional motivation for this problem comes from approximation
theory. A positive answer to Problem~\ref{problem:near-unc} would
imply the existence of a quasi-greedy basic sequence in every
infinite-dimensional Banach space. A basic sequence $(x_i)$ in a
Banach space is called \emph{quasi-greedy} if there exists a constant
$C$ such that for all $\delta\kge0$ and for all $(a_i)\in\mc_{00}$,
\eqref{eq:near-unc} above holds with
$E\keq\{i\kin\bn:\,\abs{a_i}\kgeq\delta\}$. In other words,
we can project with a \emph{uniform} constant onto sets consisting of
\emph{all} co-ordinates with ``large'' coefficients. This concept was
introduced by Konyagin and Temlyakov~\cite{KT}. One of the main
results in this paper, Theorem~\ref{thm:bdd-osc-unc}, gives a positive
answer to Problem~\ref{problem:near-unc} under some additional
assumptions on the sets of co-ordinates onto which we can project.

We will now place the above notions in a wider context. We will
explain the term `partial unconditionality' and discuss further
examples. Let $(x_i)$ be a sequence of non-zero vectors in a Banach
space. Then $(x_i)$ is a basic sequence with constant $C$ if and only
if~\eqref{eq:near-unc} holds for all $(a_i)\kin\mc_{00}$ and whenever
$E\keq\{1,\ldots,n\}$ for some $n\kin\bn$. Moreover, $(x_i)$ is an
\emph{unconditional basic sequence} if and only if~\eqref{eq:near-unc}
holds for all $(a_i)\kin\mc_{00}$ and for all finite subsets $E$ of
$\bn$. Thus for a basic sequence we can uniformly project onto initial
segments of $\bn$, whereas for an unconditional sequence we can
uniformly project onto \emph{all} finite (or indeed infinite) subsets
of $\bn$. By partial unconditionality we mean a property of a sequence
of non-zero vectors in a Banach space that lies between these two
extremes. We next describe one way in which this idea can be
formalized.

Let $\mathcal F$ be a collection of finite subsets of $\bn$. Given a
sequence $(x_i)$ of non-zero vectors in a Banach space, we say that
$(x_i)$ is \emph{$\mathcal F$-unconditional with constant $C$}
if~\eqref{eq:near-unc} holds for all $(a_i)\kin\mc_{00}$ and for all
finite sets $E$ such that either $E\kin\mathcal F$ or $E$ is an
initial segment of $\bn$. Our opening question can now be made
precise: Does every normalized, weakly null sequence have an $\mathcal
F$-unconditional subsequence?

If $\mathcal F\keq\emptyset$, then $(x_i)$ is $\mathcal
F$-unconditional with constant $C$ if and only if it is a basic
sequence with constant $C$. It is well known that for any
$\epsilon\kge0$ every normalized, weakly null sequence has a
subsequence that is a basic sequence with constant
$1\kplus\epsilon$. On the other hand if $\mathcal F$ is the set of all
finite subsets of $\bn$, then $(x_i)$ is $\mathcal F$-unconditional
with constant $C$ if and only if it is an unconditional sequence with
constant $C$. In this case our question has a negative answer: in 1974
Maurey and Rosenthal constructed a Banach space with a
normalized, weakly null basis which has no unconditional
subsequence. Note that by Rosenthal's $\ell_1$-theorem~\cite{R}, if a
space contains no normalized, weakly null sequence, then it contains
$\ell_1$ and, in particular, an unconditional basic sequence. Thus,
given a collection $\mathcal F$ of finite subsets of $\bn$, a more
general question would be to ask if every infinite-dimensional Banach
space contains an $\mathcal F$-unconditional sequence. For
unconditional sequences it was not until 1993 that the more general
question was also answered in the negative by Gowers and
Maurey \cite{GM}. They constructed a Banach space that contains no
unconditional basic sequence.

Because of the Maurey-Rosenthal and Gowers-Maurey counterexamples it
is an interesting problem to search for non-trivial examples of
partial unconditionality that lead to positive answers to the
questions we raised above. As it happens such examples occur naturally
in various contexts. We give two examples which are relevant in the
study of spreading models and asymptotic structures in Banach space
theory. A finite subset $E$ of $\bn$ is a \emph{Schreier set} if
$|E|\kleq\min E$. The collection of all Schreier sets is denoted by
$\mathcal S_1$. A sequence of non-zero vectors in a Banach space is
called \emph{Schreier-unconditional} if it is $\mathcal
S_1$-unconditional. The following result was announced in~\cite{MR}, a
proof is given in~\cite{O2}.
\begin{thm}
  \label{thm:schreier-unc}
  For each $\epsilon\kge0$, every normalized weakly null sequence in a
  Banach space has a Schreier-unconditional subsequence with constant
  $2\kplus\epsilon$.
\end{thm}
One could generalize Schreier-unconditionality by considering
higher-order Schreier families that were introduced by Alspach and
Odell~\cite{AO} and by Alspach and
Argyros~\cite{AA}. For example $\mathcal S_2$ can be defined as the
collection of disjoint unions $\bigcup_{i=1}^nF_i$ of Schreier sets
$F_1,\ldots,F_n$ with $\{\min F_1,\dots,\min F_n\}\kin\mathcal
S_1$. Unfortunately, the questions corresponding to $\mathcal S_2$
already have negative answers: the basis in the example of Maurey and
Rosenthal has no $\mathcal S_2$-unconditional subsequence, and the
space of Gowers and Maurey contains no $\mathcal S_2$-unconditional
basic sequence. However, it is worth mentioning two positive results
here. Let $\alpha$ be a countable ordinal and let $\mathcal
S_{\alpha}$ denote the Schreier family of order $\alpha$. It is shown
in~\cite{AG} that if the normalized weakly null sequence $(x_i)$ is an
$\ell_1^{\alpha}$-spreading model, then $(x_i)$ admits an $\mathcal
S_{\alpha}$-unconditional subsequence. Moreover, in~\cite{GOW} it is
shown that an $\mathcal S_{\alpha}$-unconditional normalized weakly
null sequence in $C(\mathcal S_{\alpha})$ admits an unconditional
subsequence.

The next example is about projecting onto ``$\ell_1$-subsets''. Before
giving it we need a definition. Let $X$ and $Y$ be Banach spaces, and
let $(x_i)$ and $(y_i)$ be sequences in~$X$ and in~$Y$, respectively
(either both infinite, or both finite of the same length). For
$C\kge0$ we say that $(x_i)$ and $(y_i)$ are \emph{$C$-equivalent},
written $(x_i)\stackrel{C}{\sim}(y_i)$ if there exist constants
$A\kge0$ and $B\kge0$ with $B/A\kleq C$ such that
\[
A\Bnorm{\sum_ia_ix_i}\leq\Bnorm{\sum_ia_iy_i}\leq
B\Bnorm{\sum_ia_ix_i}
\]
for all $(a_i)\kin\mc_{00}$. If only the second inequality holds, then
we say $(x_i)$ \emph{$B$-dominates} $(y_i)$, and write
$(y_i)\stackrel{B}{\lesssim}(x_i)$. Let $(e_i)$ be the unit vector
basis of $\ell_1$. Given a real number $\delta\kge 0$ and a sequence
$(x_i)$ in a Banach space, set $\mathcal
F\big(\delta,(x_i)\big)\keq\big\{ E\kin\veges{\bn}:\,(x_i)_{i\in E}
\stackrel{1/\delta}{\sim} (e_i)_{i=1}^{|E|}\big\}$. In
Section~\ref{section:convex-unc} we will present a result due to
Argyros, Mercourakis, Tsarpalias~\cite{AMT} of which the following is
an immediate consequence.
\begin{thm}
  \label{thm:ell_1-unc}
  For each $\delta\kin (0,1]$ there exists a constant $C$ such that
  every normalized, weakly null sequence has a subsequence $(x_i)$
  that is $\mathcal F\big(\delta,(x_i)\big)$-unconditional with
  constant $C$. Moreover, $C\kleq 16\log_2\big(1/\delta\big)$ for
  $\delta\kle 1/4$.
\end{thm}
As we shall later see, finding the best constant $C$ in the above
result is closely related to Problem~\ref{problem:near-unc}. Indeed,
if Problem~\ref{problem:near-unc} has a positive answer, then
the above theorem is valid with a constant $C$ not depending on
$\delta$. Another problem of interest (although we shall not address
it in this paper) is to determine which symmetric
bases could replace the unit vector basis of~$\ell_1$ in the
definition of $\mathcal F\big(\delta,(x_i)\big)$. We note that
projecting onto ``$\mc_0$-subsets'' can always be done: every basic
sequence dominates the unit vector basis of $\mc_0$. In fact, by
Theorem~\ref{thm:schreier-unc;n-form} below, for every $\epsilon\kge0$
every normalized, weakly null sequence has a basic subsequence that
$(1\kplus\epsilon)$-dominates the unit vector basis of $\mc_0$.

We now describe a different scheme for defining partial
unconditionality from the one above. We will denote by $\veges{\bn}$
the set of all finite subsets of $\bn$. Let $\mathcal F$ be a subset
of $\mc_{00}\times\veges{\bn}$. We say that the sequence $(x_i)$ is
\emph{$\mathcal F$-unconditional with constant $C$} if
\begin{equation}
  \label{eq:part-unc-refined}
  \Bnorm{\sum _{i\in E} a_ix_i}\leq C \Bnorm{\sum_{i=1}^{\infty}
    a_ix_i}
\end{equation}
holds whenever $\ba\keq(a_i)\kin\mc_{00}$, and either
$(\ba,E)\kin\mathcal F$ or $\ba$ is arbitrary and $E$ is an initial
segment of $\bn$. Observe that such a sequence is a basic sequence
with constant~$C$, \ie we can uniformly project onto
initial segments with constant~$C$. However, in general, for a given
finite set $E\ksubset \bn$ we can only project certain vectors onto
$E$ with uniform constant; namely the vectors $\sum_{i}a_ix_i$ for
which the pair $\big((a_i),E\big)$ belongs to $\mathcal F$. So this
kind of partial unconditionality is of a non-linear nature. Both
$\delta$-near-unconditionality and the quasi-greedy property are
examples of this. If we let $\mathcal F$ to be the set of all pairs
$(\ba,E)$ such that $\ba\keq(a_i)\kin\mc_{00}$ and
$E\keq\{i\kin\bn:\,\lvert a_i\rvert\kgeq\delta\}$ for some
$\delta\kge0$, then $(x_i)$ is $\mathcal F$-unconditional if and only
if it is quasi-greedy. If for a fixed $\delta\kin(0,1)$ we let
$\mathcal F_{\delta}$ be the set of pairs $(\ba,E)$ such that
$\ba\keq(a_i)\kin\mc_{00}$, $\abs{a_i}\kleq 1$ for all $i\kin\bn$, and
$E\ksubset\{i\kin\bn:\,\abs{a_i}\kgeq\delta\}$, then $(x_i)$ is
$\mathcal F_{\delta}$-unconditional if and only if it is
$\delta$-near-unconditional.
\begin{problem}
  \label{problem:quasi-greedy}
  Does every normalized, weakly null sequence have a quasi-greedy
  subsequence, or more generally, does every infinite-dimensional
  Banach space contain a quasi-greedy basic sequence?
\end{problem}
Dilworth, Kalton and Kutzarova~\cite[Theorem~5.4]{DKK} proved that if
a normalized, weakly null sequence $(x_i)$ has a spreading model not
equivalent to the unit vector basis of $\mc_0$, then for any
$\epsilon\kge0$ there is a quasi-greedy subsequence of $(x_i)$ with
constant $3\kplus\epsilon$. This is not too surprising: if we are in
some sense far from $\mc_0$, then we expect a uniform bound on the
number of large coefficients in a norm-$1$ vector, from which the
result follows by Schreier-unconditionality. This argument also shows
(using a version of Schreier-unconditionality,
Theorem~\ref{thm:schreier-unc;n-form} below) that if $(x_i)$ is a
normalized, weakly null sequence with spreading model not equivalent
to the unit vector basis of $\mc_0$, then for any $\epsilon\kge 0$ and
for any $\delta\kin (0,1)$ there is a $\delta$-near-unconditional
subsequence of $(x_i)$ with constant $1\kplus\epsilon$.

Thus Problems~\ref{problem:near-unc} and~\ref{problem:quasi-greedy}
have positive answers if we are ``far'' from $\mc_0$. However, they
are still open in general. What we do know is that one cannot hope to
find for any $\epsilon\kge 0$ subsequences of normalized, weakly null
sequences that are $\delta$-near-unconditional or quasi-greedy with
constant $1\kplus\epsilon$. We are going to prove this in
Section~\ref{section:resolutions} (Example~\ref{ex:elton>1}). We will
also show in Section~\ref{section:schreier} that a positive answer to
Problem~\ref{problem:near-unc} implies a positive answer to
Problem~\ref{problem:quasi-greedy}.

One could be forgiven for thinking that a positive answer to
Problem~\ref{problem:near-unc} would easily
imply that every normalized, weakly null sequence has an unconditional
subsequence. It is certainly true that in a
$\delta$-near-unconditional sequence we can project onto any
\emph{subset} of the co-ordinates with `large' coefficients (unlike in
a quasi-greedy sequence). However, there are two restrictions. First,
there is a normalization: $\lvert a_i\rvert\kleq1$ for all $i\kin\bn$
whenever $(\ba,E)\kin\mathcal F_{\delta}$ (where $\mathcal F_{\delta}$
is defined just before the statement of
Problem~\ref{problem:quasi-greedy}). Without this condition, for any
pair $(\ba,E)$, there would exist a positive real number~$r$ such that
$(r\ba,E)\kin\mathcal F_{\delta}$, and hence a
$\delta$-near-unconditional sequence would indeed be
unconditional. Second, even if there is a constant $K$ such that
$K(\delta)\kle K$ for all $\delta\kge0$, the subsequence that is
$\delta$-near-unconditional with constant $K$, and that we can find in
a given normalized, weakly null sequence may very well depend on
$\delta$. In other words there is no obvious reason why a positive
answer to Problem~\ref{problem:near-unc} would find, in every
normalized weakly null sequence, a subsequence that is
$\delta$-near-unconditional with constant $K$ for
\emph{all}~$\delta\kge0$ (which again would be unconditional). Note
that the standard diagonal argument would give a subsequence that is
$\delta$-near-unconditional with constant $N(\delta)\kplus 2K$ for
\emph{all}~$\delta\kge0$, where $N$ is an integer-valued function with
$\lim _{\delta\to0}N(\delta)\keq\infty$.

This paper will be organized as follows. In the next section we
introduce the concept of a bounded-oscillation-unconditional basic
sequence, which is a new type of partial unconditionality. We then
prove our main result (Theorem~\ref{thm:bdd-osc-unc}) that states that
every normalized, weakly null sequence has a
bounded-oscillation-unconditional subsequence. The combinatorial
machinery that we set up in order to prove our main result will be
subsequently  employed to prove other partial unconditionality
results. We will use it in Section~\ref{section:schreier} to give a
new proof of Schreier unconditionality. Here we will also deduce
Elton's theorem from our main result,
Theorem~\ref{thm:bdd-osc-unc}. We will then prove that a positive
answer to Problem~\ref{problem:near-unc} implies a positive answer to
Problem~\ref{problem:quasi-greedy}.

In Section~\ref{section:variants} we introduce various constants
similar to the constant $K(\delta)$ defined above. These will allow us
to quantify the relationships between various notions of partial
unconditionality. We will also show that for solving
Problem~\ref{problem:near-unc} one can restrict attention to the
Banach spaces of continuous functions on countable, compact, Hausdorff
spaces. In Section~\ref{section:c_0-problem} we raise the question
whether there is a uniform constant~$C$ such that every sequence
equivalent to the unit vector basis of $\mc_0$ has an unconditional
subsequence with constant~$C$. This turns out to be closely related to
Problem~\ref{problem:near-unc}. The proof will again use our
combinatorial machinery.

In the following two sections we revisit convex unconditionality of
Argyros, Mercourakis, Tsarpalias~\cite{AMT}, and unconditionality of
certain sequences in spaces of continuous functions. Using our
approach we give new proofs of known results and establish a duality
between them and near-unconditionality.

In the final section we will have a closer look at our combinatorial
machinery. We give a necessary and sufficient condition for a positive
answer to Problem~\ref{problem:near-unc} (\cf
Proposition~\ref{prop:matching}). To decide if this condition can be
satisfied in general one is lead to consider certain combinatorial
data attached to subsequences of a normalized, weakly null
sequence. We will study this data on its own right as a purely
combinatorial object. Our results will be used at the end to give an
example that among other thing shows that $\sup_{\delta>0}K(\delta)$
is strictly greater than~$1$.

\section{Main results}

\label{main}
Given a sequence $\ba\keq(a_i)$ of real
numbers, we define its \emph{support} to be the set
$\supp{\ba}\keq\{i\kin\bn:\,a_i\kneq 0\}$. If this set is finite we
call~$\ba$ \emph{finitely supported}. Recall that $\mc_{00}$ denotes
the space of finitely supported sequences of real numbers. Given
$\ba\keq(a_i)\kin\mc_{00}$ and a subset $E$ of $\bn$ we define the
\emph{oscillation} $\text{osc}(\ba,E)$ \emph{of} $\ba$ \emph{over} $E$
as
\[
\osc{\ba}{E}=\sup\left\{\frac{|a_i|}{|a_j|}:\,i, j\in E,\ 
a_j\neq 0\right\}.
\]
For subsets $E$ and $F$ of $\bn$ we write $E\kle F$ if $m\kle n$ for
all $m\kin E$ and for all $n\kin F$. We say that a sequence
$E_1,\dots,E_n$ of subsets of $\bn$ is \emph{successive} if
$E_1\kle\dots\kle E_n$. A decomposition $E\keq\bigcup_{j=1}^nE_j$ of a
finite set $E$ will be called a \emph{Schreier decomposition} if
$E_1\kle\dots\kle E_n$ is a successive sequence of non-empty sets such
that $n\kleq\min E_1$, \ie the set $\{\min E_1\dots\min E_n\}$ belongs
to $\mathcal S_1$.

We now come to the main definition. Let $C,D,d\kin[1,\infty)$. We
say that a basic sequence $(x_i)$ in a Banach space $X$ is
\emph{$(D,d)$-bounded-oscillation-unconditional with constant $C$}
if for every $\ba\keq(a_i)\in\mc_{00}$, and for every
finite set $E\ksubset \bn$ with $\osc{\ba}{E}\kleq D$, we have
\[
\Big\| \sum _{i\in E} a_i x_i \Big\| \leq C \Big\| \sum _{i=1}
^{\infty} a_i x_i \Big\|
\]
provided $E$ has a Schreier decomposition $E\keq\bigcup_{j=1}^nE_j$ such
that $\osc{\ba}{E_j}\kleq d$ for each $j=1,\dots,n$. Note that without
this proviso the sequence $(x_i)$ would be a $1/D$-near-unconditional
sequence.

Our main theorem is the following.
\begin{thm}
  \label{thm:bdd-osc-unc}
  For all $d\kin[1,\infty)$, there is a constant $C\kleq 8d$ such that
  for all $D\kin [1,\infty)$ and for any $\epsilon\kge 0$ every
  normalized, weakly null sequence has a subsequence that is a
  $(D,d)$-bounded-oscillation-unconditional basic sequence with
  constant $C\kplus\epsilon$.
\end{thm}
Note that if $\ba\keq(a_i)\kin\mc_{00},\ E\kin\veges{\bn}$ and
$\osc{\ba}{E}\kleq D$, then we can write $E$ as the disjoint union of
$n\kleq \big\lfloor\log_2\big(D\big)\big\rfloor\kplus 1$ sets
$E_1,\ldots,E_n$ such that $\osc{\ba}{E_j}\kleq 2$ for each
$j\keq1,\ldots,n$. So without the assumption that the sets in a
Schreier decomposition are successive the above result would be a
positive answer to Problem~\ref{problem:near-unc}.

A key ingredient in the proof of Theorem~\ref{thm:bdd-osc-unc} is a
purely combinatorial result which we call the Matching Lemma
(Theorem~\ref{thm:matching}). In its proof and in much of this paper
we will be making heavy use of infinite Ramsey theory. For this reason
we now recall some notation and results from the subject. For a subset
$M$ of $\bn$ we denote by $\veges{M}$ the set of all finite subsets of
$M$ and by $\vegtelen{M}$ the set of all infinite subsets of $M$. The
power-set $2^{\bn}$ of $\bn$ is equipped with the product topology,
and all subspaces will carry the subspace topology. A collection
$\mathcal U\ksubset\vegtelen{\bn}$ is said to be \emph{Ramsey} if for
all $L\kin\vegtelen{\bn}$ there exists $M\kin\vegtelen{L}$ such that
either $\vegtelen{M}\ksubset\mathcal U$ or
$\vegtelen{M}\ksubset\mathcal U^{\complement}$, where $\mathcal
U^{\complement}\keq\vegtelen{\bn}\backslash \mathcal U$ denotes the
complement of $\mathcal U$. One example of an infinite Ramsey theorem,
due to Galvin and Prikry~\cite{GP}, states that every Borel subset
$\mathcal U$ of $\vegtelen{\bn}$ is Ramsey. More generally, whenever
$\vegtelen{\bn}$ is partitioned into finitely many Borel sets, every
infinite subset $L$ of $\bn$ has an infinite subset $M$ such that
$\vegtelen{M}$ is contained in one of the Borel sets of the
partition. The strongest result of this type was proved by
Ellentuck~\cite{Ell}; his result concerns topological
characterizations of Ramsey sets. In all our applications (and indeed
in most applications to Banach space theory) it will suffice to know
that open sets (and hence closed sets) are Ramsey. This was first
proved by Nash-Williams~\cite{NW}. Following tradition we will often
talk about colourings instead of partitions. This and other pieces of
terminology will be introduced as we go along. For a very good
introduction to infinite Ramsey theory see~\cite{B}. An extensive
account is presented in~\cite{GRS}.

We need one final piece of notation before stating
Lemma~\ref{thm:matching}. For subsets $A,B$ of $\bn$ we write $A\prec
B$ if $A$ is an initial segment of $B$.
\begin{thm}[Matching Lemma]
  \label{thm:matching}
  Let $n\kin\bn$. Assume that for every infinite subset $M$ of $\bn$
  we are given a successive sequence
  \[
  F_1^M<\dots < F_n^M
  \]
  of non-empty, finite subsets of $M$. Further assume that for each
  $j=1,\dots,n$ the function $F_j\colon\vegtelen{\bn}\to\veges{\bn},\
  M\mapsto F^M_j$, is continuous. Then for all $N\kin\vegtelen{\bn}$
  there exist $L,M\kin\vegtelen{N}$ such that
 \begin{mylist}{(ii)}
 \item
   for each $j=1,\dots,n$ either $F_j^L\prec F_j^M$ or
   $F_j^M\prec F_j^L$, and
 \item
   $L\cap M=\ds \bigcup _{j=1}^nF_j^L\cap F_j^M$.
 \end{mylist}
\end{thm}

\begin{proof}
  We begin by setting up some notation. Let $F_L\keq\bigcup _{j=1}^n
  F_j^L$ for each $L\kin\vegtelen{\bn}$ . We are going to define a
  \emph{finite colouring} $c$ of pairs $(L,l)$, where $L$ is an infinite
  subset of $\bn$ and $l\kin L$. In other words we are going to define
  a function $c$ on the set of all such pairs taking values in some
  finite set whose elements will be referred to as \emph{colours}. So
  fix $L\kin\vegtelen{\bn}$ and $l\kin L$. If $l\kin F_i^L$ for some
  $i\keq 1,\dots,n$, then we set $c(L,l)\keq i$. If $l\knotin F_L$,
  and the minimum of $\{l'\kin L:\,l'\kge l\}\cap F_L$ belongs to
  $F_i^L$, then we set $c(L,l)\keq i+$. Finally, if $l\kge\max
  F_L$, then we set $c(L,l)\keq +$. Clearly there exists $l_0\kin
  L$ with $c(L,l_0)\keq +$, and for such an $l_0$ we have
  $c(L,l)\keq +$ for all $l\kin L$ with $l\kgeq l_0$.

  We now prove a preliminary result.
  \begin{claim}
    For all pairs $(F,X)$, where $F\kin\veges{\bn}$ and
    $X\kin\vegtelen{\bn}$, there exist $Y\kin\vegtelen{X}$ and a
    colour $\lambda$ such that $F\kle Y$ and $c(F\cup V,\min
    V)\keq\lambda$ for all $V\kin\vegtelen{Y}$.
  \end{claim}
  To see this define a finite colouring $d$ of $\vegtelen{\bn}$ by
  setting $d(V)\keq c(F\cup V,\min V)$ for every
  $V\kin\vegtelen{\bn}$. It follows from the continuity of the maps
  $F_j$ that if $\lambda$ is a colour other than $+$, the
  corresponding \emph{colour-class}, \ie the collection
  $\{V\kin\vegtelen{\bn}:\,d(V)\keq\lambda\}$ is an open subset of
  $\vegtelen{\bn}$. It follows that the colour-class of $+$ is
  closed. Since open sets and closed sets are Ramsey, it follows that
  there is an infinite subset $Y$ of $X$ all whose infinite subsets
  have the same colour. Replacing $Y$ by a smaller set if necessary we
  may clearly assume that $F\kle Y$.

  We now turn to the proof of Theorem~\ref{thm:matching}. Fix
  $N\kin\vegtelen{\bn}$.  We shall build infinite subsets $L$ and $M$
  of $N$ from recursively constructed sequences $l_1\kleq
  l_2\kleq\dots,\ m_1\kleq m_2\kleq\dots$ of positive integers in
  $N$. Along the way we shall also construct a sequence $P_0\ksupset
  P_1\ksupset P_2\ksupset\dots$ of infinite subsets of $N$, and
  sequences $(\lambda_k)_{k=0}^{\infty}$ and $(\mu_k)_{k=0}^{\infty}$
  of colours. To start the construction apply the Claim with $F\keq
  \emptyset$ and $X\keq N$. This yields an infinite subset $Y$ of $X$
  and a colour $\lambda$ such that $c(V,\min V)\keq\lambda$ for all
  $V\kin\vegtelen{Y}$. Let us set $P_0\keq Y$ and
  $\lambda_0\keq\mu_0\keq\lambda$.

  For the recursive step suppose that $k\kgeq0$ and that $l_r,\ m_r$
  for $1\kleq r\kleq k$ and $P_r,\ \lambda_r,\ \mu_r$ for $0\kleq
  r\kleq k$ have been chosen. We also assume that setting
  $A_k\keq\{l_r:\,1\kleq r\kleq k\}$ and $B_k=\{m_r:\,1\kleq r\kleq
  k\}$ the following hold.
  \begin{myeqnarray}
    \label{eq:matching:develop}
    A_k<P_k\quad\text{and}\quad B_k<P_k,\\
    \label{eq:matching:predict}
    c(A_k\cup Q,\min Q)\keq\lambda_k,\ c(B_k\cup Q,\min
    Q)\keq\mu_k\quad\text{for all}\quad Q\kin\vegtelen{P_k}.
  \end{myeqnarray}
  Note that when $k\keq0$ these assumptions are satisfied by the
  choice of $P_0$. To choose $l_{k+1}$ and $m_{k+1}$ we consider four
  cases.
  \begin{mylist}{\textit{Case 4}}
  \item[\textit{Case 1.}]
    If $\lambda _k\keq\mu_k\keq i$ for some $i\kin\{1,\dots,n\}$, then
    we choose $l_{k+1}\keq m_{k+1}$ to be an arbitrary element of
    $P_k$.
  \item[\textit{Case 2.}]
    If one of
    \begin{mylist}{(a)}
    \item
      neither $\lambda_k$ nor $\mu_k$ belongs to $\{1,\dots,n\}$,
    \item
      $\{\lambda_k,\mu_k\}\keq\{i,j+\}$ for some $1\kleq i\kle j$, or
    \item
      at least one of $\lambda_k$ and $\mu_k$ is $+$
    \end{mylist}
    holds, then we choose $l_{k+1}$ and $m_{k+1}$ to be distinct
    elements of $P_k$.
  \item[\textit{Case 3.}]
    If $\lambda_k\keq i$ and either $\mu_k\keq j$ for some $1\kleq
    j\kle i$ or $\mu_k\keq j+$ for some $1\kleq j\kleq i$, then we set
    $l_{k+1}\keq l_k$ and choose $m_{k+1}$ to be an arbitrary element
    of $P_k$.
  \item[\textit{Case 4.}]
    If $\mu_k\keq i$ and either $\lambda_k\keq j$ for some $1\kleq
    j\kle i$ or $\lambda_k\keq j+$ for some $1\kleq j\kleq i$, then we
    set $m_{k+1}\keq m_k$ and choose $l_{k+1}$ to be an arbitrary
    element of $P_k$.
  \end{mylist}
  Note that when $k\keq0$ only Cases~1 and~2 can arise, since
  $\lambda_0\keq\mu_0$. When $k\kgeq 1$ we have $l_k\kleq l_{k+1}$ and
  $m_k\kleq m_{k+1}$ in all the cases, as required. Let us at this
  point set $l_0=m_0=0$ in order to avoid having to consider the first
  step of the construction separately from the recursive
  steps. Observe that for any $k\kgeq 0$, if $l_{k+1}\kge l_k$, then
  $l_{k+1}\kin P_k$. Similarly, if $m_{k+1}\kge m_k$, then
  $m_{k+1}\kin P_k$.

  To complete the recursive step we need to choose $P_{k+1},\
  \lambda_{k+1}$ and $\mu_{k+1}$. First set $A_{k+1}\keq
  A_k\cup\{l_{k+1}\}$ and $B_{k+1}\keq B_k\cup\{m_{k+1}\}$. Then apply
  the Claim with $F\keq A_{k+1}$ and $X\keq P_k$ to obtain an infinite
  subset $\tilde{P}$ of $P_k$ and a colour $\lambda_{k+1}$ such that
  $A_{k+1}\kle\tilde{P}$ and $c(A_{k+1}\cup Q,\min
  Q)\keq\lambda_{k+1}$ for all $Q\kin\vegtelen{\tilde{P}}$. Now apply
  the Claim again with $F\keq B_{k+1}$
  and $X\keq\tilde{P}$ to obtain an infinite subset $P_{k+1}$ of
  $\tilde{P}$ and a colour $\mu_{k+1}$ such that $B_{k+1}\kle P_{k+1}$
  and $c(B_{k+1}\cup Q,\min Q)\keq\mu_{k+1}$ for all
  $Q\kin\vegtelen{P_{k+1}}$. With these choices it is clear that the
  assumptions for the next recursive step
  (\ie~\eqref{eq:matching:develop} and~\eqref{eq:matching:predict}
  with~$k$ replaced by~$k\kplus 1$) are satisfied. Observe that if
  $l_{k+1}\keq l_k$, then $\lambda_{k+1}\keq\lambda_k$, and if
  $m_{k+1}\keq m_k$, then $\mu_{k+1}\keq\mu_k$.

  Having completed the recursive construction let us put
  $L\keq\{l_r:\,r\kin\bn\}$ and $M\keq\{m_r:\,r\kin\bn\}$.  Notice
  that for any $k\kgeq 0$, if $l_{k+1}\kge l_k$, then
  $L_k\keq\{l_r:\,r\kge k\}$ is a subset of $P_k$. Indeed, for $r\kge
  k$ we have $l_r\keq l_{s+1}\kge l_s$ for some $s$ with $k\kleq s\kle
  r$, and hence $l_r\kin P_s\ksubset P_k$. So if in addition $L$ is
  infinite, then $c(L,l_{k+1})\keq c(A_k\cup L_k,\min
  L_k)\keq\lambda_k$. Similarly, for any $k\kgeq0$ if $m_{k+1}\kge
  m_k$, then $M_k\keq\{m_r:\,r\kge k\}$ is a subset of $P_k$, and if
  in addition $M$ is infinite, then we have $c(M,m_{k+1})\keq
  c(B_k\cup M_k,\min M_k)\keq\mu_k$.

  We will now verify that $L$ and $M$ are indeed infinite sets. We
  argue by contradiction. Assume, for example, that $L$ is
  finite. Then for some $k_0\kin\bn$ we have $l_{k+1}\keq l_k$ for all
  $k\kgeq k_0$. It follows that for every $k\kgeq k_0$ we applied Case~3
  in the $k^{\text{th}}$ step of the recursion. Hence for every $k\kgeq
  k_0$ we have $m_{k+1}\kge m_k$ and $c(M,m_{k+1})\keq\mu_k\neq+$. This
  contradiction shows that $L$ is infinite. Similar reasoning gives
  that $M$ must also be infinite.

  Next let us fix $i\kin\{1,\dots,n\}$. We need to show that
  either $F_i^L\prec F_i^M$ or $F_i^M\prec F_i^L$. We argue by
  contradiction. Suppose that there exist  $l\kin F_i^L$ and $m\kin
  F_i^M$ such that
  \begin{equation}
    \label{eq:matching:init-seg}
    l\neq m\quad \text{and}\quad \{l'\in F_i^L:\,l'<l\}=\{m'\in
    F_i^M:\,m'<m\}.
  \end{equation}
  For some $k\kgeq0$ and $k'\kgeq0$ we have $l\keq l_{k+1}\kge
  l_k$ and $m\keq m_{k'+1}\kge m_{k'}$. Then $\lambda_k\keq c(L,l)\keq
  i$ and $\mu_{k'}\keq c(M,m)\keq i$. From now on assume that $k\kleq
  k'$ (the case $k\kgeq k'$ is similar). There exists $k''$ with
  $k\kleq k''\kleq k'$ such that $m_k\keq m_{k''}\kle m_{k''+1}$, and
  so $\mu_{k}\keq\mu_{k''}\keq c(M,m_{k''+1})$. From $m_{k''+1}\kleq
  m$ and $c(M,m)\keq i$ we deduce that the colour $\mu_k$ is either
  $j$ or $j+$ for some $j$ with $1\kleq j\kleq i$. Hence in the
  $k^{\text{th}}$ recursive step we applied either Case~1 or
  Case~3. Case~1 leads to $l\keq l_{k+1}\keq m_{k+1}\kin F_i^M$ and
  $l\kleq m$ which contradicts~\eqref{eq:matching:init-seg}, whereas
  Case~3 gives $l_{k+1}\keq l_k$ contradicting the choice of~$k$.

  We are left to show that $L\cap M\ksubset\bigcup_{i=1}^nF_i^L\cap
  F_i^M$ (the reverse inclusion being obvious). Let $l$ belong to
  $L\cap M$. There exist $k\kgeq0$ and $k'\kgeq0$ such that $l\keq
  l_{k+1}\keq m_{k'+1}$ and $l_{k+1}\kge l_k,\ m_{k'+1}\kge
  m_{k'}$. Then $l\kin P_k\backslash P_{k+1}$ and $l\kin
  P_{k'}\backslash P_{k'+1}$, from which we get $k\keq k'$. So we have
  $l\keq l_{k+1}\keq m_{k+1}$ and $l_{k+1}\kge l_k,\ m_{k+1}\kge
  m_k$. It follows immediately that in the $k^{\text{th}}$ step of the
  recursion we must have been in Case~1. Hence for some $i\keq
  1,\dots,n$ we have $c(L,l)\keq\lambda_k\keq i$ and
  $c(M,l)\keq\mu_k\keq i$, \ie $l\kin F_i^L\cap F_i^M$, as
  required. This completes the proof of Theorem~\ref{thm:matching}.
\end{proof}
Some minor modifications of the proof and a simple diagonalization
procedure yields a corollary that we shall refer to as the Schreier
version of the Matching Lemma. The diagonalization process will be
used later on, so we state it separately as an abstract principle. A
family $\mathcal A$ of finite subsets of $\bn$ is \emph{thin} if no
element of $\mathcal A$ is the proper initial segment of another
element of $\mathcal A$. The following result was proved by
Nash-Williams~\cite{NW}: if a thin family $\mathcal A$ is finitely
coloured, then for all $L\kin\vegtelen{\bn}$ there exists
$M\kin\vegtelen{L}$ such that $\veges{M}\cap\mathcal A$ is
monochromatic. To see this, simply give an infinite set $L$ the colour
of its unique initial segment in $\mathcal A$ (introducing a new
colour for infinite sets with no initial segment in $\mathcal
A$). Clearly, each colour-class is either open or closed, so the
result follows. An easy diagonalization argument then gives the
following result. (A much stronger statement is given by Pudl\' ak and
R\" odl~\cite{PR}.)
\begin{prop}
  \label{prop:diagonal-thin-ramsey}
  Let $\mathcal A\ksubset\veges{\bn}$ be a thin family. For each
  $k\kin\bn$ let $S_k$ be a finite set, and let $c\colon\mathcal
  A\to\bigcup_{k=1}^{\infty}S_k$ be a colouring of $\mathcal A$ so
  that for all $F\kin\mathcal A$ we have $c(F)\kin S_k$, where
  $k\keq\min F$. Then for all $L\kin\vegtelen{\bn}$ there exists
  $M\kin\vegtelen{L}$ such that if $A,B\kin\veges{M}\cap\mathcal A$
  and $\min A\keq\min B$, then $c(A)\keq c(B)$.

  \qed
\end{prop}
We are now ready to state and prove the promised corollary to (the
proof of) Theorem~\ref{thm:matching}.
\begin{cor}[Schreier version of the Matching Lemma]
  \label{cor:schreier-matching-lemma}
  Assume that for each $M\kin\vegtelen{\bn}$ we have a positive
  integer $n_M$ and non-empty finite subsets $A_M,\ F^M_1\kle\dots\kle
  F^M_{n_M}$ of $M$ such that
  \[
  \bigcup_{j=1}^{n_M}F^M_j\subset A_M\qquad\textrm{and}\qquad
  n_M\leq\min F^M_1=\min A_M.
  \]
  Further assume that the function $M\mapsto
  A_M\colon\vegtelen{\bn}\to\veges{\bn}$ is continuous, that the
  family $\mathcal A\keq\{A_M:\,M\kin\vegtelen{\bn}\}$ is thin, and
  that for all $L,M\kin\vegtelen{\bn}$ if $A_L\keq A_M$, then $n_L\keq
  n_M$ and $F^L_j\keq F^M_j$ for each $j\keq1,\dots,n_L$. Then for all
  $N\kin\vegtelen{\bn}$ there exists $L, M\kin\vegtelen{N}$ with
  $n_L\keq n_M$ such that 
  \begin{mylist}{(ii)}
  \item
    for each $j\keq1,\dots,n_L$ either $F_j^L\prec F_j^M$, or
    $F_j^M\prec F_j^L$, and
  \item
    $L\cap M =\displaystyle \bigcup _{j=1}^{n_L} F_j^L\cap F_j^M$.
  \end{mylist}
\end{cor}
\begin{proof}
  We first define a colouring of $\mathcal A$ by giving each $A_M$,
  $M\kin\vegtelen{\bn}$,  the colour $n_M$. This is well-defined by
  the  assumptions. By Proposition~\ref{prop:diagonal-thin-ramsey}
  there exists $N_1\kin\vegtelen{N}$ such that for all
  $L,M\kin\vegtelen{N_1}$ if $\min F^L_1\keq\min F^M_1$, then $n_L\keq
  n_M$.

  We now follow the proof of Theorem~\ref{thm:matching}. We define the
  colouring~$c$ on pairs $(L,l)$ as before. Although this time~$c$ is
  a possibly infinite colouring, the colouring~$d$ used in the proof
  of the Claim is finite, so the Claim remains valid. We then carry
  out the recursive construction that produces the sets $L$ and
  $M$. The only changes we need is to work inside $N_1$ (rather than
  $N$), and to replace in Cases~1--4 each occurence of $\{1,\dots,n\}$
  by $\bn$. The verification that $L$ and $M$ are infinite is the same
  as before.

  At this point we need to insert the
  observation that $\min F^L_1\keq\min F^M_1$. To see this choose
  $k\kgeq0$ and $k'\kgeq0$ such that $\min F^L_1\keq l_{k+1}\kge
  l_{k}$ and $\min F^M_1\keq m_{k'+1}\kge m_{k'}$ so that
  $\lambda_k\keq \mu_{k'}\keq1$. Assume that $k\kleq k'$ (the case
  $k'\kleq k$ is similar). Then $m_k\kleq m_{k'}\kle m_{k'+1}$, and
  so $\mu_k$ is either~$1$ or~$1+$. It follows that in the
  $k^{\text{th}}$ step of the recursion we were either in Case~1, in
  which case we have $m_{k+1}\keq l_{k+1}$ (and $k\keq k'$), as
  required, or we were in Case~3, in which case we obtain $l_k\keq
  l_{k+1}$, which contradicts the choice of~$k$.
 
  We now have $n_L\keq n_M$ by our initial application of
  Proposition~\ref{prop:diagonal-thin-ramsey}. To finish the proof we
  verify properties~(i) and~(ii) exactly as in the proof
  Theorem~\ref{thm:matching} (letting~$n$ in the proof stand
  for~$n_L$).
\end{proof}
Applications of the Matching Lemma and of its Schreier version will
require two further lemmas. To motivate the first one of these we now
give a preview of the type of argument that will follow. Consider the
general problem of starting with a normalized, weakly null sequence
$(x_i)$ and seeking a subsequence with a certain desired
property. Arguing by contradiction, we assume that for all
$M\kin\vegtelen{\bn}$ we have a witness $w_M$ to the lack of the
desired property in the subsequence $(x_i)_{i\in M}$. The witness
$w_M$ will then give rise in a very natural way to finitely many
subsets $F^M_1\kle F^M_2\kle\dots$ of $M$. Lemma~\ref{lem:selection}
below will allow us to choose $w_M$ from the set of all possible
witnesses for $M$ in a ``continuous'' way so that among other things
the assumptions of the Matching Lemma or its corollary are
satisfied. In typical examples a witness $w_M$ has as a constituent
part some functional $x^{\ast}_M$. A priori we will not be able to assume
that the support of $x^{\ast}_M$, \ie the set
$\supp{x^{\ast}_M}\keq\{i\kin\bn:\,x^{\ast}_M(x_i)\kneq0\}$ is
contained in $M$, precisely because we lack unconditionality. In
Lemma~\ref{lem:support-condition} we show that we can stabilize, \ie
we can pass to some infinite set with respect to which the property
$\supp{x^{\ast}_M}\ksubset M$ can be assumed (provided the choice of
$x^{\ast}_M$ had already been made in a ``continuous'' manner).
\begin{lemma}
  \label{lem:selection}
  Let $\Omega\keq\bigcup _{r=1}^{\infty}\Omega _r$ be an arbitrary set
  written as the union of a countably infinite collection of its
  subsets. Let
  \[
  \Phi \colon \vegtelen{\bn} \to 2^{\Omega}\backslash \{ \emptyset \}
  \]
  be a function into the set of non-empty subsets of $\Omega$. Assume
  that for all $r\kin\bn$ and for all $L,M\kin\vegtelen{\bn}$ we have
  \[
  L\cap\{1,\dots,r\}=M\cap\{1,\dots,r\}\quad\implies\quad\Phi
  (L)\cap\Omega_r=\Phi (M)\cap\Omega _r.
  \]
  Then there is a function
  \[
  \phi \colon \vegtelen{\bn} \to \Omega
  \]
  such that
  \begin{mylist}{(ii)}
    \item
      $\phi(M)\kin\Phi(M)$ for all $M\kin\vegtelen{\bn}$, and
    \item
      $\phi$ is continuous if $\Omega$ is given the discrete topology.
  \end{mylist}
\end{lemma}

\begin{proof}
  Fix a well-ordering of $\Omega$. For $M\in \vegtelen{\bn}$ let
  \[
  r(M)=\min\{r\in\bn:\, \Phi(M)\cap\Omega _r\neq\emptyset\}.
  \]
  Define $\phi(M)$ to be the least element of
  $\Phi(M)\cap\Omega_{r(M)}$ in our chosen well-ordering. We claim
  that $\phi\colon\vegtelen{\bn}\to\Omega$ has the required
  properties.

  Clearly $\phi(M)\kin\Phi(M)$ for all infinite subsets $M$ of
  $\bn$. To verify continuity fix $M\kin\vegtelen{\bn}$ and set $r\keq
  r(M)$. Let $[r]\keq\{1,\dots,r\}$. If $L\kin\vegtelen{\bn}$
  satisfies $L\cap[r]\keq M\cap[r]$, then for each $1\kleq r'\kleq r$
  we have $\Phi(L)\cap\Omega_{r'}\keq\Phi(M)\cap\Omega_{r'}$, which is
  the empty set for $r'\kle r$ and is not empty for $r'\keq r$. It
  follows that $r(L)\keq r(M)$, which in turn implies that
  $\phi(L)\keq\phi(M)$. This shows that $\phi$ maps the neighbourhood
  $\{L\kin\vegtelen{\bn}:\,L\cap[r]\keq M\cap[r]\}$ of $M$ onto
  $\phi(M)$.
\end{proof}

\begin{lemma}
  \label{lem:support-condition}
  Let $\mc_0$ be equipped with the topology of pointwise convergence
  on $\bn$. Let $f\colon\vegtelen{\bn}\to\mc_0,\ M\mapsto f_M$, be a
  continuous function such that every sequence in the image of $f$ has
  a cluster point in $\mc_0$. Then for every $\epsilon\kge0$ and for
  every $M\kin\vegtelen{\bn}$ there exists $N\kin\vegtelen{M}$ such
  that for all $P\kin\vegtelen{N}$ we have
  \[
  \sum _{i\in N\backslash P} |f_P(i)|\leq \epsilon,
  \]
  \ie the support $\supp{f_P}\keq\{i\kin\bn:\,f_P(i)\kneq0\}$ of $f_P$
  relative to the set $N$ is contained in $P$ up to a small
  perturbation.
\end{lemma}

\begin{proof}
  For $L\kin\vegtelen{\bn}$ let us write $L'$ as a temporary notation
  for $L\backslash\{\min L\}$. For $F\kin\veges{\bn}$ and
  $\delta\kge0$ let $\mathcal U_{F,\delta}$ be the collection of all
  infinite subsets $L$ of $\bn$ for which we have
  \[
  |f_{F\cup L'}(\min L)|<\delta.
  \]
  As a preliminary step we first prove the following claim. Given
  $F\kin \veges{\bn}$ and $L\kin\vegtelen{\bn}$, there exists
  $\tilde{L}\kin\vegtelen{L}$ such that
  $\vegtelen{\tilde{L}}\ksubset\mathcal U_{F,\delta}$. Indeed, the
  continuity of $f$ implies that $\mathcal U_{F,\delta}$ is an open
  set, and hence it is Ramsey. Thus there exists
  $\tilde{L}\kin\vegtelen{L}$ such that either
  $\vegtelen{\tilde{L}}\ksubset\mathcal U_{F,\delta}$ or
  $\vegtelen{\tilde{L}}\ksubset \mathcal
  U_{F,\delta}^{\complement}$. So to prove the claim we need to
  exclude the second alternative. We argue by contradiction. Assume
  that $\vegtelen{\tilde{L}}\ksubset\mathcal
  U_{F,\delta}^{\complement}$. Let $l_1\kle l_2\kle\dots$ be an
  enumeration of $\tilde{L}$, and for $n\kin\bn$ let
  $L_n\keq\{l_i:\,i\kge n\}$. Then $L_n\cup\{l_i\}\kin\mathcal
  U_{F,\delta}^{\complement}$, and hence
  \[
  |f_{F\cup L_n}(l_i)|\geq \delta\quad \text{whenever}\ 1\leq i\leq n.
  \]
  Let $x\kin\mc_0$ be a cluster point of the sequence
  $(f_{F\cup L_n})_{n=1}^{\infty}$. From the above we have
  $|x(l_i)|\kgeq\delta$ for all $i\kin\bn$ contradicting that $x$ is
  an element of $\mc_0$. This completes the proof of the claim.

  To prove Lemma~\ref{lem:support-condition} let us fix
  $\epsilon\kge0$ and $M\kin\vegtelen{\bn}$. Choose real numbers
  $\epsilon_i\kge0,\ i\keq1,2,\dots$, such that
  $\sum_{i=1}^{\infty}\epsilon_i\kle\epsilon$. We shall now
  recursively construct a sequence $n_1\kle n_2\kle\dots$ of positive
  integers, and a sequence $L_0\ksupset L_1\ksupset L_2\ksupset\dots$
  of infinite subsets of $\bn$ as follows. To start with, set $L_0\keq
  M$. Assume that for some $k\kgeq1$ we have chosen $n_i$ for $1\kleq
  i\kle k$ and $L_i$ for $0\kleq i\kle k$. Let $F_1,\dots,F_K$ be an
  enumeration of the power-set of $\{n_1,\dots,n_{k-1}\}$. Then choose
  a chain $L_{k-1}\keq\tilde{L}_0\ksupset\tilde{L}_1\ksupset\dots
  \ksupset\tilde{L}_K$ of infinite sets such that for each
  $j\keq1,\dots,K$ we have $\vegtelen{\tilde{L}_j}\ksubset\mathcal
  U_{F_j,\epsilon_k}$. This can be done by our preliminary claim. Now
  set $n_k\keq\min\tilde{L}_K$ and
  $L_k\keq\tilde{L}_K\backslash\{n_k\}$. Note that $n_k\kin L_{k-1},\
  L_k\ksubset L_{k-1},\ n_k\kle L_k$ and
  \begin{equation}
    \label{eq:small_on_next_coord}
    |f_{F\cup Q}(n_k)|<\epsilon _k \qquad \text{for all}\
     F\subset\{n_1,\dots,n_{k-1}\},\ \text{and}\ Q\in\vegtelen{L_k}.
  \end{equation}
  Having completed the recursive construction, set
  $N\keq\{n_1,n_2,\dots\}$. It is clear that
  $N\kin\vegtelen{M}$. Given any $P\kin\vegtelen{N}$, if $k\kin\bn$
  with $n_k\knotin P$, then $P\keq F\cup Q$, where
  \[
  F=P\cap\{n_1,\dots ,n_{k-1}\},\ \text{and}\ Q=P\backslash F\in
  \vegtelen{L_k}.
  \]
  Hence from~\eqref{eq:small_on_next_coord} we have
  $|f_P(n_k)|\kle\epsilon_k$. It follows that
  \[
  \sum_{n\in N\backslash P} |f_P(n)|<\sum_{k=1}^{\infty}
  \epsilon_k<\epsilon,
  \]
  as required.
\end{proof}

We are now ready to present a proof for
Theorem~\ref{thm:bdd-osc-unc}. It will be convenient to use the
following definition of an $\epsilon$-net $F$ for a subset $S$ of
$\br^d$, where $\epsilon\kge0$ and $d\kin\bn$: for every
$(\alpha_j)_{j=1}^d\kin S$ there exists $(\beta_j)_{j=1}^d\kin F$ such
that $\beta_j\kleq\alpha_j\kleq\beta_j\kplus\epsilon$ for each
$j\keq1,\dots,d$.

\begin{proof}[Proof of Theorem~\ref{thm:bdd-osc-unc}]
  Fix $C,D,d\kin[1,\infty)$. Assume that $(x_i)$ is a normalized,
  weakly null sequence no subsequence of which is
  $(D,d)$-bounded-oscillation-unconditional basic sequence with
  constant $C$. We shall deduce that $C\kleq 8d$. Fix
  $\epsilon\kin(0,1)$ and then choose an increasing function
  $\gamma\colon\bn\to\bn$ such that
  $\ds\lim_{k\to\infty}\gamma(k)\keq\infty$ and
  \begin{equation}
    \label{eq:bdd-osc-unc:growth-of-gamma}
    \gamma(k)+Dk\leq (1+\epsilon) \gamma(k-1)\qquad \text{for all}\
    k\geq 2.
  \end{equation}
  For example, we can take $\gamma(k)\keq k^2$ for $k\kgeq k_0$ and
  $\gamma(k)\keq k_0^2$ for $k\kle k_0$, where $k_0$ is sufficiently
  large.

  After passing to a subsequence we may assume that
  $(x_i)$ is a basic sequence with constant~$1\kplus\epsilon$. Then in
  particular for all $\ba\keq(a_i)\kin\mc_{00}$ we have
  \begin{equation}
    \label{eq:bdd-osc-unc:schreier-unc}
    \lVert \ba\rVert _{\ell _{\infty}}\leq 2(1+\epsilon)\Big
    \lVert \sum _{i=1}^{\infty} a_ix_i\Big \rVert \leq 4\Big
    \lVert \sum _{i=1}^{\infty} a_ix_i\Big \rVert.
  \end{equation}
  We now show that for every infinite subset $M$ of $\bn$ there
  exists a triple $(\ba,x^{\ast},F)$, which we shall call a \emph{witness
  for $M$}, with the following properties.
  \begin{myeqnarray}
    \ba=(a_i)\kin \mc _{00},\quad x^{\ast}\kin
    B_{X^{\ast}}\quad\textrm{and}\quad F\kin\veges{\bn};\\
    \label{eq:bdd-osc-unc:witness-A}
    F\ksubset A\ksubset M\quad \text{and}\quad \min F\keq\min A,\quad
    \text{where}\ A=\supp{\ba};\\
    \label{eq:bdd-osc-unc:witness-d}
    F\  \text{has a Schreier decomposition}\ F=\bigcup_{j=1}^nF_j\
    \text{such that}\newline
    \osc{\ba}{F_j}\leq d\quad \text{for each}\quad j=1,\dots,n;\\
    \label{eq:bdd-osc-unc:witness-D}
    1\leq a_i\leq D\quad \text{and}\quad x^{\ast}(x_i)>0\quad \text{for
    all}\ i\in F;\\
    \label{eq:bdd-osc-unc:witness-g(k)}
    \frac{C}{2(1+\epsilon)(2+\epsilon)}\lVert x\rVert < \sum _{i\in F}
    a_ix^{\ast}(x_i)\leq \gamma(k)+D,\quad \text{where}\newline k=\min
    F,\quad \text{and}\quad x=\sum _{i\in M} a_ix_i.
  \end{myeqnarray}
  To see this let us fix $M\kin\vegtelen{\bn}$. Since $(x_i)_{i\in M}$
  is not $(D,d)$-bounded-oscillation-unconditional with constant
  $C$, there exist $\bb\keq(b_i)\kin\mc_{00}$ with $\supp{\bb}\ksubset
  M$, and a finite
  subset $E$ of $M$ with a Schreier decomposition $E\keq\bigcup
  _{j=1}^nE_j$ such that $\osc{\bb}{E}\kleq D,\ \osc{\bb}{E_j}\kleq d$
  for each $j\keq 1,\dots,n$, and
  \[
  \Big\lVert \sum _{i\in E}b_ix_i \Big\rVert >C \lVert y\lVert,
  \]
  where $y\keq\sum_{i\in M} b_ix_i$. We may then choose $x^{\ast}\kin
  B_{X^{\ast}}$ such that
  \[
  \sum _{i\in E} b_ix^{\ast}(x_i)>C \lVert y\rVert.
  \]
  Replacing $\bb$ and $x^{\ast}$ by $-\bb$ and $-x^{\ast}$ if
  necessary, we may assume that if we let $E'\keq\{i\kin E:\,b_i\kge
  0,\ x^{\ast}(x_i)\kge 0\}$, then we still have
  \[
  \sum _{i\in E'} b_ix^{\ast}(x_i)>\frac{C}{2} \lVert y\rVert.
  \]
  By homogeneity, we may also assume that $\min\{b_i:\,i\kin
  E'\}\keq 1$, and hence $1\kleq b_i\kleq D$ for all $i\kin
  E'$. Finally, let $k'\kgeq \min E'$ be minimal so that
  \[
  \sum _{%
    \begin{subarray}{l}
      i\in E'\\[1pt]
      i>k'
    \end{subarray}}
  b_i x^{\ast}(x_i)\leq\gamma(k').
  \]
  Then we have
  \[
  \sum _{i\in E'} b_ix^{\ast}(x_i)\leq (1+\epsilon) \sum _{%
    \begin{subarray}{l}
      i\in E'\\[1pt]
      i\geq k'
    \end{subarray}}  b_i x^{\ast}(x_i).
  \]
  Indeed, this is clear when $k'\keq\min E'$, whereas if $k'\kge\min
  E'$, then by the triangle-inequality,
  by~\eqref{eq:bdd-osc-unc:growth-of-gamma} and by the choice of~$k'$
  we have
  \begin{eqnarray*}
    \sum _{i\in E'} b_ix^{\ast}(x_i) & \leq & Dk'+\gamma(k')\leq
    (1+\epsilon)\gamma(k'-1)\\
    & \leq & (1+\epsilon) \sum _{%
      \begin{subarray}{l}
        i\in E'\\[1pt]
        i\geq k'
      \end{subarray}}  b_i x^{\ast}(x_i),\\[-12pt]
  \end{eqnarray*}
  as claimed. Set $F\keq\{i\kin E':\,i\kgeq k'\}$. For each $i\kin\bn$
  set $a_i\keq b_i$ when $i\kgeq \min F$ and $a_i\keq 0$ when $i\kle
  \min F$, and let $\ba \keq(a_i)_{i\in M}$. It is now routine to
  verify that $(\ba,x^{\ast},F)$ is a witness for~$M$ as defined above.

  The next step is to select witnesses in a continuous manner using
  Lemma~\ref{lem:selection}. Let $\Omega$ be the set of all witnesses
  of all infinite subsets of $\bn$, and for each $M\kin
  \vegtelen{\bn}$ let $\Phi (M)$ be the (non-empty) set of all
  witnesses for $M$. For each $r\kin\bn$ let $\Omega_r$ be the set
  of elements $(\ba,x^{\ast},F)$ of $\Omega$ that satisfy
  $\max\supp{\ba}\kleq r$. It is easy to verify
  that the conditions of Lemma~\ref{lem:selection} are satisfied. It
  follows that there exists a function $\phi \colon \vegtelen{\bn}\to
  \Omega$ such that $\phi (M)\kin \Phi (M)$, \ie $\phi (M)$ is a
  witness for $M$ for all $M\kin \vegtelen{\bn}$, and $\phi$ is
  continuous if $\Omega$ is given the discrete topology. For each
  $M\kin \vegtelen{\bn}$ let $\phi(M)\keq(\ba_M,x^{\ast}_M,F_M)$, and let
  $n_M$ be the positive integer such that $F_M$ has a Schreier
  decomposition $F_M\keq\bigcup_{j=1}^{n_M}F^M_j$ with
  $\osc{\ba_M}{F^M_j}\kleq d$ for each $j\keq 1,\dots,n_M$. We will
  also use the notation
  \[
  \ba _M=(a^M_i),\quad x_M=\sum _{i\in M}a^M_ix_i,\quad
  \text{and}\quad A_M=\supp{\ba _M}.
  \]
  By the proof of Lemma~\ref{lem:selection} we may assume that for
  each $M\kin \vegtelen{\bn}$ there is an $r\kin \bn$ such that $\Phi
  (M)\cap \Omega _s\keq\emptyset$ if $1\kleq s\kle r$ and $\phi (M)$
  is the least element of $\Phi(M)\cap \Omega_r$ with respect to some
  fixed well-ordering of $\Omega$. It follows that for $L, M\kin
  \vegtelen{\bn}$ if $A_L$ is an initial segment of $A_M$, then we
  must have $\phi (L)\keq\phi (M)$. In particular
  \[
  \mathcal A=\{ A\in \veges{\bn}:\, A=A_M\ \text{for some}\ M\in
  \vegtelen{\bn} \}
  \]
  is a thin family, and we are in the situation of
  Corollary~\ref{cor:schreier-matching-lemma}.

  We shall now select infinite subsets $N_1\ksupset N_2\ksupset N_3$
  of  $\bn$ stabilizing various parameters. To select $N_1$ we use
  Lemma~\ref{lem:support-condition}. Let $f\colon\vegtelen{\bn}\to\mc
  _0$ be the function mapping $M\kin\vegtelen{\bn}$ to
  $\big(x^{\ast}_M(x_i)\big)\kin\mc_0$. Note that
  this is the only place, where we use the weakly null property of the
  sequence $(x_i)$. It follows easily from the continuity of $\phi$
  and from the w$^{\ast}$-compactness of $B_{X^{\ast}}$ that $f$ is
  continuous with respect to the topology of pointwise convergence on
  $\mc_0$, and that the image of $f$ has compact closure. Hence, by
  Lemma~\ref{lem:support-condition}, there exists an infinite subset
  $N_1$ of $\bn$ such that for all $P\kin\vegtelen{N_1}$ we have
  \begin{equation}
    \label{eq:bdd-osc-unc:support-condition}
    \sum _{i\in N_1\backslash P} \lvert x^{\ast}_P(x_i)\rvert <\epsilon.
  \end{equation}
  We next choose an infinite subset $N_2$ of $N_1$ using infinite
  Ramsey theory. We colour $\mathcal A$ by giving $A_M$,
  $M\kin\vegtelen{\bn}$, colour $(r_j)_{j=1}^{n_M}\kin\bn^{n_M}$ if
  \[
  (1+\epsilon)^{r_j-1}\leq\min\{a^M_i:\,i\in F^M_j\}<
  (1+\epsilon)^{r_j}\qquad\text{for each}\ j=1,\dots,n_M.
  \]
  This colouring is well-defined, \ie the colour of $A\kin\mathcal A$
  does not depend on the choice of infinite set $M$ with $A\keq A_M$.  
  Note that for each $k\kin\bn$ the family $\{A\kin\mathcal A:\,\min
  A\keq k\}$ is finitely coloured. An application of
  Proposition~\ref{prop:diagonal-thin-ramsey} now gives
  $N_2\kin\vegtelen{N_1}$ such that for all $L,M\kin\vegtelen{N_2}$ if
  $\min F_L\keq\min F_M$, then $n_L\keq n_M$ and
  \begin{equation}
    \label{eq:bdd-osc-unc:fix-minima}
    \frac{a^M_i}{a^L_i}\geq \frac{1}{d(1+\epsilon)}\qquad \text{for
      all}\ i\in F^L_j\cap F^M_j,\ j=1,\dots,n_L.
  \end{equation}
  For our final stabilization we choose for each $k\kin\bn$ an
  $\epsilon/k$-net $S_k$ of $[0,\gamma(k)\kplus D]^k$ (in the sense
  defined just before the start of this proof) together with an
  ordering of its elements. Given $M\kin\vegtelen{\bn}$, let
  $k\keq\min F_M$ and let $(w_j)_{j=1}^k$ be the least element of
  $S_k$ satisfying
  \[
  w_j\leq \sum _{i\in F^M_j} a^M_ix^{\ast}_M(x_i) \leq w_j+\epsilon/k
  \qquad  \text{for each}\ j=1,\dots,n_M.
  \]
  We shall refer to $(w_j)_{j=1}^k$ as the \emph{weight-colour} of
  $M$. This colouring of $\vegtelen{\bn}$ induces a colouring of the
  family $\mathcal A$ satisfying the assumptions of
  Proposition~\ref{prop:diagonal-thin-ramsey}. Hence there is an
  infinite subset $N_3$ of $N_2$ so that for all
  $L,M\kin\vegtelen{N_3}$ if $\min F_L\keq\min F_M$ then $L$ and $M$
  have the same weight-colour.

  To finish the proof we apply the Schreier version of the Matching
  Lemma (Corollary~\ref{cor:schreier-matching-lemma}). As observed
  earlier, the assumptions of the corollary are satisfied. So we can
  find $L,M\kin\vegtelen{N_3}$ with $n_L\keq n_M$ such that
  \begin{mylist}{(ii)}
  \item
    for each $j\keq 1,\dots,n_L$ either $F^L_j\prec F^M_j$, or
    $F^M_j\prec F^L_j$, and
  \item
    $L\cap M =\displaystyle \bigcup _{j=1}^{n_L} F^L_j\cap F^M_j$.
  \end{mylist}
  Note that in particular $\min F_L\keq\min F_M$, and hence $L$ and
  $M$ have the same weight-colour, say $(w_j)_{j=1}^k$, where
  $k\keq\min F_L$. Set
  \[
  J=\big\{ j\in \{1,\dots,n_L\}:\, F^L_j\prec F^M_j\big\}.
  \]
  Interchanging $L$ and $M$ and replacing $J$ by $\{1,\dots,n_L\}
  \ksetminus J$ if necessary, we may assume that
  \begin{equation}
    \label{eq:bdd-osc-unc:half}
    \sum _{j\in J} w_j\geq \frac{1}{2} \sum _{j=1}^{n_L} w_j.
  \end{equation}
  We now establish a number of inequalities. First we have
  \begin{multline*}
    \lVert x_M\rVert \geq x^{\ast}_L(x_M) = \sum _{i\in M} a^M_i
    x^{\ast}_L(x_i) \geq \sum _{i\in L\cap M} a^M_i x^{\ast}_L(x_i)
    -\sum _{i\in M\setminus L} \lvert a^M_i\rvert\, \lvert
    x^{\ast}_L(x_i)\rvert \\
    \geq \sum _{i\in L\cap M} a^M_i x^{\ast}_L(x_i) -
    4\epsilon\lVert x_M \rVert,
  \end{multline*}
  where the last inequality comes
  from~\eqref{eq:bdd-osc-unc:schreier-unc} and
  from~\eqref{eq:bdd-osc-unc:support-condition} applied with $P\keq
  L$. We now obtain the following sequence of inequalities (the steps
  are justified below).

  \begin{eqnarray*}
    \Big( 1+4\epsilon\Big) \lVert x_M \rVert & \geq &
    \sum _{i\in L\cap M} a^M_i x^{\ast}_L(x_i)\\
    & \geq & \frac{1}{d(1+\epsilon)} \sum
    _{\phtm{i\in F^L_j\cap F^M_j j=1}j=1}^{n_L} \sum _{\phtm{i\in
        F^L_j\cap F^M_j j=1}i\in F^L_j\cap F^M_j} a^L_i x^{\ast}_L(x_i)\\
    & \geq & \frac{1}{d(1+\epsilon)} \sum _{\phtm{i\in F^L_j j\in
        J}j\in J} \sum _{\phtm{i\in F^L_j j\in J}i\in F^L_j} a^L_i
    x^{\ast}_L(x_i)\\
    & \geq & \frac{1}{d(1+\epsilon)} \sum _{j\in J} w_j \geq
    \frac{1}{2d(1+\epsilon)} \sum _{j=1}^{n_L} w_j\\
    & \geq & \frac{1}{2d(1+\epsilon)} \Big( \sum _{\phtm{i\in F^M_j
        j=1}j=1}^{n_L} \sum _{\phtm{i\in F^M_j j=1} i\in F^M_j} a^M_i
    x^{\ast}_M(x_i) -\epsilon\Big)\\
    & \geq & \frac{1}{2d(1+\epsilon)}
    \frac{C}{2(1+\epsilon)(2+\epsilon)}\lVert x_M\rVert -
    \frac{\epsilon}{2d(1+\epsilon)}\,4 \lVert
    x_M\rVert.
  \end{eqnarray*}

  \noindent
  The second line uses~\eqref{eq:bdd-osc-unc:fix-minima} and the third
  line uses the definition of~$J$. For the next two lines we use the
  fact that $L$ and $M$ both have weight-colour $(w_j)_{j=1}^{n_L}$,
  and we also use~\eqref{eq:bdd-osc-unc:half}. For the last inequality
  we apply~\eqref{eq:bdd-osc-unc:witness-g(k)} from the definition of
  a witness, and inequality~\eqref{eq:bdd-osc-unc:schreier-unc}
  (from~\eqref{eq:bdd-osc-unc:witness-D} we have $\lVert
  \ba\rVert_{\ell _{\infty}}\kgeq 1$). We have thus shown that
  \[
  C\leq 4d(1+\epsilon)^2(2+\epsilon) \Big( 1+4\epsilon
  +\frac{2\epsilon}{d(1+\epsilon)}\Big).
  \]
  Since $\epsilon$ was arbitrary it follows that $C\kleq 8d$, as
  claimed.
\end{proof}

\section{Schreier- and near-unconditionality}
\label{section:schreier}

In this section we give new proofs of two results quoted in the
Introduction. We begin with Schreier-unconditionality. It is not
difficult to apply Theorem~\ref{thm:bdd-osc-unc} with $d\keq 1$ and a
diagonal process to show that for any $\epsilon\kge 0$, every
normalized, weakly null sequence has a Schreier-unconditional
subsequence with constant $8\kplus\epsilon$. The better constant
claimed in Theorem~\ref{thm:schreier-unc} follows by a straightforward
diagonal argument from the statement below. For $M\ksubset\bn$ and
$n\kin\bn$ we denote by $M^{(\leq n)}$ the collection of subsets of
$M$ of size at most~$n$. So a sequence is $\bn^{(\leq
  n)}$-unconditional if we can uniformly project onto sets of size at
most~$n$.
\begin{thm}
  \label{thm:schreier-unc;n-form}
  Fix $n\kin\bn$ and $\epsilon\kge 0$. Every normalized weakly null
  sequence has a $\bn^{(\leq n)}$-unconditional subsequence with
  constant $1\kplus\epsilon$.
\end{thm}
\begin{proof}
  Let $C\kin[1,\infty)$ and assume that $(x_i)$ is a normalized weakly
  null sequence no subsequence of which is $\bn^{(\leq
  n)}$-unconditional with constant~$C$. We need to show that $C\kleq
  1$.

  Let $M\kin\vegtelen{\bn}$. By our assumption there exists a triple
  $(\ba,F,x^{\ast})$, called a \emph{witness for~$M$}, such that
  \begin{myeqnarray}
    \ba\keq(a_i)\kin\mc_{00},\quad F\kin M^{(\leq n)}\quad
    \text{and}\ x^{\ast}\kin B_{X^{\ast}};\\
    \sum _{i\in M} a_ix_i\in S_X;\\
    \sum _{i\in F} a_ix^{\ast}(x_i)>C. \label{eq:schreier-unc;large-on-F}
  \end{myeqnarray}
  Let $\Omega$ be the set of all witnesses of all infinite subsets
  of $\bn$ equipped with the discrete topology. By
  Lemma~\ref{lem:selection} we obtain a continuous function
  $\phi\colon\vegtelen{\bn}\to\Omega$ such that $\phi (M)$ is a
  witness for $M$ for all $M\kin\vegtelen{\bn}$. For each
  $M\kin\vegtelen{\bn}$ we write
  \[
  \phi (M)=(\ba_M,F_M,x^{\ast}_M),
  \]
  where $\ba_M\keq(a^M_i)$, and we let $x_M\keq\sum_{i\in M}
  a^M_ix_i$.

  We will now select infinite subsets $N_1\ksupset N_2\ksupset N_3$ of
  $\bn$. We first choose $N_1$ so that $(x_i)_{i\in N_1}$ is a basic
  sequence with basis constant at most~$2$, say. Then in particular
  for any $M\kin\vegtelen{N_1}$ we have
  \begin{equation}
    \label{eq:schreier-unc;basic-const}
    \sup_{i\in M}\abs{a^M_i}\leq 4\Big\lVert\sum_{i\in M} a^M_{i}
    x_{i}\Big\rVert =4.
  \end{equation}
  We next fix an arbitrary positive real number $\delta$. We then
  select $N_2\kin\vegtelen{N_1}$ so that for all
  $L,M\kin\vegtelen{N_2}$ we have $\abs{F_L}\keq\abs{F_M}$ and
  \begin{equation}
    \label{eq:schreier-unc;stab-coeff}
    \big\lVert (a^L_i)_{i\in F_L} -(a^M_i)_{i\in F_M} \big\rVert
    _{\ell _1} < \delta.
  \end{equation}
  This is done by a straightforward use of infinite Ramsey
  theory. Finally, using Lemma~\ref{lem:support-condition} we obtain
  $N_3\kin\vegtelen{N_2}$ such that for all $P\kin\vegtelen{N_3}$ we
  have
  \begin{equation}
    \label{eq:schreier-unc;small-on-rest}
    \sum _{i\in N_3\backslash P}\lvert x^{\ast}_P(x_i)\rvert <\delta.
  \end{equation}
  After these stabilizations we apply Theorem~\ref{thm:matching} with
  $n\keq 1$ to obtain infinite subsets $L,M$ of $N_3$ such that
  either $F_L\kprec F_M$ or $F_L\kprec F_M$, and $L\cap M\keq F_L\cap
  F_M$. The choice of $N_2$ implies that in fact we have $F_L\keq
  F_M\keq L\cap M$.  We now estimate $x^{\ast}_L(x_M)$ to obtain the
  required inequality. First, we write $x^{\ast}_L(x_M)$ as
  \begin{multline}
    \label{eq:schreier-unc;L-on-M}
    \sum _{i\in M} a^M_i x^{\ast}_L(x_i) = \sum _{i\in F_L} a^L_i
    x^{\ast}_L(x_i) + \sum _{i\in F_M} \big( a^M_i-a^L_i \big)
    x^{\ast}_L(x_i)\\ + \sum _{i\in M\backslash F_M} a^M_i
    x^{\ast}_L(x_i).
  \end{multline}
  We then estimate the three terms on the right-hand side
  of~\eqref{eq:schreier-unc;L-on-M} as follows.  Applying
  property~\eqref{eq:schreier-unc;large-on-F} of a witness to~$L$
  gives~$C$ as a lower bound on the first term. 
  Applying~\eqref{eq:schreier-unc;stab-coeff} to the second term,
  and~\eqref{eq:schreier-unc;basic-const},
  \eqref{eq:schreier-unc;small-on-rest} to the third term give upper
  bounds leading to
  \[
  1=\lVert x_M \rVert \geq \big\lvert x^{\ast}_L(x_M) \big\rvert \geq
  C - \delta -4\delta.
  \]
  Since $\delta$ was arbitrary, it follows that $C\kleq 1$, as
  claimed.
\end{proof}
\begin{rem}
  If $(x_i)$ is a normalized basic sequence with basis constant $C$,
  then for all $(a_i)\kin\mc_{00}$ we have
  \[
  \abs{a_n}=\Bnorm{\sum_{i=1}^na_ix_i-\sum_{i=1}^{n-1}a_ix_i}\leq
  2C\Bnorm{\sum_{i=1}^{\infty}a_ix_i}\qquad\text{for all}\ n\kin\bn.
  \]
  We shall often use this to assume after passing to a subsequence
  $(x_i)$ of a given normalized, weakly null sequence that
  $\abs{a_n}\kleq 4\bnorm{\sum_{i=1}^{\infty}a_ix_i}$, say, for all
  $(a_i)\kin\mc_{00}$ and for all $n\kin\bn$. The constant~$4$ is
  often adequate, however, sometimes we will need to be able to
  replace~$4$ by $1\kplus\epsilon$ for any given $\epsilon\kge 0$. We
  can do this by applying Theorem~\ref{thm:schreier-unc;n-form} with
  $n\keq 1$.
\end{rem}
We now turn to Elton's theorem on near-unconditional sequences. As
mentioned in the Introduction, it is this result that raises
Problem~\ref{problem:near-unc} --- the main focus in this paper.
\begin{proof}[Proof of Theorem~\ref{thm:near-unc}]
  Let $\delta\kin (0,1]$ and let $(x_i)$ be a normalized, weakly null
  sequence in a Banach space. An application of our main result,
  Theorem~\ref{thm:bdd-osc-unc}, with $d\keq D\keq 1/\delta$ gives,
  for each $\epsilon\kge 0$, a $\delta$-near-unconditional subsequence
  of $(x_i)$ with constant~$8/\delta\kplus\epsilon$. As we mentioned
  in the Introduction, a better constant of order
  $\log\big(1/\delta\big)$ can be obtained as follows. Set $d\keq
  D\keq 2$ and pass to a $(D,d)$-bounded-oscillation-unconditional
  subsequence $(y_i)\ksubset (x_i)$ with constant~$17$, say. We show
  that $(y_i)$ is $\delta$-near-unconditional with constant $17k$,
  where $k\keq\big\lfloor\log_2\big(1/\delta\big)\big\rfloor\kplus
  1$. Indeed, let $(a_i)\kin\mc_{00}$ with $\abs{a_i}\kleq 1$ for all
  $i\kin\bn$, and let $E\ksubset
  \{i\kin\bn:\,\abs{a_i}\kgeq\delta\}$. Set
  \[
  E_j=\big\{i\in E:\,2^{-j}<\abs{a_i}\leq 2^{-(j-1)}\big\},\qquad
  \text{for each}\ j=1,\dots,k.
  \]
  Since $\osc{(a_i)}{E_j}\kleq 2$ we have
  \[
  \Bnorm{\sum_{i\in E_j} a_iy_i}\leq 17\Bnorm{\sum_{i=1}^{\infty}
    a_iy_i}\qquad \text{for each}\ j=1,\dots,k.
  \]
  Hence, by the triangle-inequality we get
  \[
  \Bnorm{\sum_{i\in E} a_iy_i}\leq 17k\Bnorm{\sum_{i=1}^{\infty}
    a_iy_i},
  \]
  as claimed. Note that $17k<18\log_2\big(1/\delta\big)$ if $\delta$
  is sufficiently small.
\end{proof}
Let us mention that recently Lopez-Abad and Todorcevic~\cite{LAT} also
gave new proofs of Theorems~\ref{thm:near-unc}
and~\ref{thm:schreier-unc} based on results on pre-compact families of
finite subsets of $\bn$.

We conclude this section by proving that a positive answer to
Problem~\ref{problem:near-unc} implies a positive answer to
Problem~\ref{problem:quasi-greedy}.
\begin{prop}
  If $\sup_{\delta>0}K(\delta)\kle \infty$, then there exists a
  constant $C$ such that every normalized, weakly null sequence has a
  quasi-greedy subsequence with constant $C$.
\end{prop}
\begin{proof}
  Let $C\kge 2\sup_{\delta>0}K(\delta)\kplus1$. Fix $\epsilon\kin
  (0,1)$, and for each $n\kin\bn$ set
  $\delta_n\keq\epsilon/n$. Given a normalized, weakly null sequence
  $(x_i)$, we apply a diagonal procedure to extract a subsequence
  $(y_i)$ such that for each $n\kin\bn$ the tail $(y_i)_{i=n}^{\infty}$ is
  $\delta_n$-near-unconditional with constant
  $K(\delta_n)\kplus\epsilon$. Passing to a further subsequence, if
  necessary, we may assume that $(y_i)$ is a basic sequence with
  constant $1\kplus\epsilon$, and moreover $\abs{a_n}\kleq
  (1\kplus\epsilon)\bnorm{\sum_{i=1}^{\infty}a_iy_i}$ for all
  $(a_i)\kin\mc_{00}$ and for all $n\kin\bn$. For the latter property
  we used Theorem~\ref{thm:schreier-unc;n-form} with $n\keq 1$. We
  will now show that $(y_i)$ is quasi-greedy with constant $C$
  provided $\epsilon$ is sufficiently small.

  Given $(a_i)\kin\mc_{00}$ and $\delta\kin (0,1]$, we need to show
  that
  \begin{equation}
    \label{eq:quasi-need}
    \Bnorm{\sum_{i\in E}a_iy_i}\leq C\norm{x},
  \end{equation}
  where $E\keq\{i\kin\bn:\,\abs{a_i}\kgeq\delta\}$ and
  $x\keq\sum_{i=1}^{\infty}a_iy_i$. We may clearly assume that
  $\sup_i\abs{a_i}\keq 1$, which implies that $\norm{x}\kgeq
  (1\kplus\epsilon)^{-1}\kge 1/2$. Now choose the smallest $n\kin\bn$
  such that $\delta_n\kleq\delta$. Note that $(n-1)\delta\kleq\epsilon
  \kle 2\epsilon\norm{x}$. Hence
  \begin{eqnarray*}
    \Bnorm{\sum_{i\in E,\ i<n}a_iy_i} &\leq& \Bnorm{\sum_{i=1}^{n-1}
      a_iy_i}+\Bnorm{\sum_{i\notin E,\ i<n}a_iy_i}\\[6pt]
    &\leq& (1+\epsilon)\norm{x} + (n-1)\delta \leq
    (1+3\epsilon)\norm{x}.
  \end{eqnarray*}
  On the other hand, since $(y_i)_{i=n}^{\infty}$ is
  $\delta_n$-near-unconditional with constant
  $K(\delta_n)\kplus\epsilon$, we have
  \[
  \Bnorm{\sum_{i\in E,\ i\geq n}a_iy_i}\leq \big(K(\delta_n) \kplus
  \epsilon \big)\Bnorm{\sum_{i=n}^{\infty}a_iy_i}\leq \big(K(\delta_n)
  \kplus \epsilon \big)(2\kplus\epsilon)\norm{x}.
  \]
  Now~\eqref{eq:quasi-need} follows for suitable $\epsilon$ by the
  triangle inequality.
\end{proof}

\section{Variants of near-unconditionality}
\label{section:variants}

In the following sections we will be considering various problems that
turn out to be related to the Elton problem. In order to make this
relationship precise we will now introduce some variants of the
constant $K(\delta)$, and explain the relationships between them.
To begin with we recall the definition of $K(\delta)$ in a slightly
different way. Given $\delta\kin (0,1]$ and a normalized, weakly null
sequence $(x_i)$, let $K((x_i),\delta)$ be the least real number~$C$
such that $(x_i)$ is $\delta$-near-unconditional with constant
$C$,\ie for all $(a_i)\kin\mc_{00}$ and for all
$E\ksubset\{i\kin\bn:\,\abs{a_i}\kgeq\delta\}$ if 
$\sup_i\abs{a_i}\kleq 1$, then $\bnorm{\sum_{i\in E}a_ix_i}\kleq
C\bnorm{\sum_{i=1}^{\infty}a_ix_i}$. Observe that
\[
K(\delta)=\sup_{(x_i)}\underset{(y_i)\subset(x_i)}{\inf\phantom{p}}
K((y_i),\delta),
\]
where the supremum is taken over all normalized, weakly null sequences
$(x_i)$ and the infimum over all subsequences $(y_i)$ of $(x_i)$.
Recall that the normalization $\sup_i\abs{a_i}\kleq 1$ in the
definition is essential (see remarks in the Introduction). We will now
introduce three other constants $K', L$ and $L'$. For $K'$ we will use
the normalization $\bnorm{\sum_ia_ix_i}\kleq1$, whereas in the
definition of $L, L'$ we restrict to vectors all whose non-zero
coefficients are ``large''. Below we repeated the definition of $K$
for the convenience of the reader.
\enlargethispage{28pt}
\begin{defn}
  \label{defn:constants}
  Let $\delta\kin(0,1]$ and let $(x_i)$ be a normalized, weakly null
  sequence in a Banach space. Each supremum below is over all
  normalized, weakly null sequences $(y_i)$ and the
  infimum is taken over all subsequences $(z_i)$ of $(y_i)$.
  \begin{multline*}
    K((x_i),\delta)=\inf\bigg\{C:\,\Bnorm{\sum_{i\in E}a_ix_i}\kleq
    C\Bnorm{\sum_{i=1}^{\infty}a_ix_i}\ \text{whenever}\
    (a_i)\kin\mc_{00},\\
    \shoveright{%
      E\ksubset\{i\kin\bn:\,\abs{a_i}\kgeq\delta\},\
      \sup_i\abs{a_i}\kleq1\bigg\}}\\[12pt]
    \shoveleft{%
      K(\delta)=\sup_{(y_i)}\underset{(z_i)\subset(y_i)}{\inf\phantom{p}}
      K((z_i),\delta)}\\[12pt]
    \shoveleft{%
      K'((x_i),\delta)=\inf\bigg\{C:\,\Bnorm{\sum_{i\in E}a_ix_i}\kleq
      C\Bnorm{\sum_{i=1}^\infty a_ix_i}\ \text{whenever}\
      (a_i)\kin\mc_{00},}\\
    \shoveright{%
      E\ksubset\{i\kin\bn:\,\abs{a_i}\kgeq\delta\},\
      \Bnorm{\sum_{i=1}^\infty a_ix_i}\kleq1\bigg\}}\\[12pt]
    \shoveleft{%
      K'(\delta)=\sup_{(y_i)}\underset{(z_i)\subset(y_i)}{\inf\phantom{p}}
      K'((z_i),\delta)}\\[12pt]
    \shoveleft{%
    L((x_i),\delta)=\inf\bigg\{C:\,\Bnorm{\sum_{i\in E}a_ix_i}\kleq
      C\Bnorm{\sum_{i=1}^{\infty}a_ix_i}\ \text{whenever}\
      \ba\keq(a_i)\kin\mc_{00},}\\
    \shoveright{%
      \abs{a_i}\kgeq\delta\ \forall i\kin\supp{\ba},\
      E\kin\veges{\bn},\ \sup_i\abs{a_i}\kleq1\bigg\}}\\[12pt]
    \shoveleft{%
      L(\delta)=\sup_{(y_i)}\underset{(z_i)\subset(y_i)}{\inf\phantom{p}}
      L((z_i),\delta)}\\[12pt]
  \end{multline*}
  \begin{multline*} 
    L'((x_i),\delta)=\inf\bigg\{C:\,\Bnorm{\sum_{i\in E}a_ix_i}\kleq
    C\Bnorm{\sum_{i=1}^\infty a_ix_i}\ \text{whenever}\
    \ba\keq(a_i)\kin\mc_{00},\\
    \shoveright{%
      \abs{a_i}\kgeq\delta\ \forall i\kin\supp{\ba},\
      E\kin\veges{\bn},\ \Bnorm{\sum_{i=1}^\infty
        a_ix_i}\kleq1\bigg\}}\\[12pt]
    \shoveleft{%
      L'(\delta)=\sup_{(y_i)}\underset{(z_i)\subset(y_i)}{\inf\phantom{p}}
      L'((z_i),\delta)}\\[-12pt]
  \end{multline*}
\end{defn}
The following result establishes some relationships between the
constants we just introduced. It shows in particular that for solving
Problem~\ref{problem:near-unc} we are free to choose the
normalization. In many situations it is more convenient to work with
the constants $K'$ and $L'$ instead of $K$ and $L$.
\begin{prop}
  \label{prop:near-unc-constants}
  Let $K, K', L$ and $L'$ be the functions defined above.
  \begin{mylist}{(ii)}
  \item
    If $0\kle\delta_1\kle\delta_2\kleq 1$, then
    $K'(\delta_2)\kleq K(\delta_1)$ and $L'(\delta_2)\kleq
    L(\delta_1)$.
  \item
    If $0\kle\delta\kleq 1$, then $L(\delta)\kleq K(\delta)$ and
    $L'(\delta)\kleq K'(\delta)$.
  \end{mylist}
  In particular we have
  $\sup_{\delta>0}K(\delta)\keq\sup_{\delta>0}K'(\delta)
  \kgeq\sup_{\delta>0}L(\delta)\keq\sup_{\delta>0}L'(\delta)$.
\end{prop}
\begin{proof}
  (ii) is clear from definition. To see~(i) let $(x_i)$ be a
  normalized, weakly null sequence. By
  Theorem~\ref{thm:schreier-unc;n-form} we may assume, after passing
  to a subsequence if necessary, that
  \[
  \sup_i\abs{a_i}\kleq \frac{\delta_2}{\delta_1}
  \Bnorm{\sum_{i=1}^{\infty}a_ix_i}\qquad \text{for all}\
  (a_i)\kin\mc_{00}.
  \]
  Now given $(a_i)\kin\mc_{00}$ with $\bnorm{\sum_ia_ix_i}\kleq 1$ and
  $E\ksubset\{i\kin\bn:\,\abs{a_i}\kgeq\delta_2\}$, if $b_i\keq
  \frac{\delta_1}{\delta_2}a_i$ for all $i\kin\bn$, then
  $\sup_i\abs{b_i}\kleq 1$ and
  $E\ksubset\{i\kin\bn:\,\abs{b_i}\kgeq\delta_1\}$. It follows that
  for any subsequence $(y_i)$ of $(x_i)$ we have
  $K'((y_i),\delta_2)\kleq K((y_i),\delta_1)$ and
  $L'((y_i),\delta_2)\kleq L((y_i),\delta_1)$.

  It remains to show that $\sup_{\delta}K(\delta)\kleq
  \sup_{\delta}K'(\delta)$ and $\sup_{\delta}L(\delta)\kleq
  \sup_{\delta}L'(\delta)$. We show the second inequality (the proof
  of the first one is similar). Assume that $L'\keq
  \sup_{\delta}L'(\delta)\kle\infty$, and let $\delta\kin (0,1]$. We
  will show that $L(\delta)\kleq L'$. Let $(x_i)$ be a normalized,
  weakly null sequence. Fix $\epsilon\kin (0,1]$ and positive real
  numbers $M_n$ such that $n\kle\epsilon (M_n\kminus 1)$ for all
  $n\kin\bn$. After passing to a subsequence, if necessary, we may
  assume that $(x_i)$ is a basic sequence with constant
  $1\kplus\epsilon$. Then using a standard diagonal argument we pass
  to a subsequence $(y_i)$ of $(x_i)$ such that for each $n\kin\bn$
  we have $L'\big((y_i)_{i\geq n},\delta/M_n\big)\kleq
  L'\kplus\epsilon$. We claim that $L((y_i),\delta)\kleq (L'\kplus
  2\epsilon)(1\kplus 3\epsilon)$.

  Given $\ba\keq(a_i)\kin\mc_{00}$ with $\delta\kleq\abs{a_i}\kleq 1$
  for all $i\kin\supp{\ba}$, and $E\kin\veges{\bn}$, set
  $x\keq\sum_ia_iy_i$ and choose $n\kin\bn$ minimal so that
  \[
  \Bnorm{\sum_{i\geq n}a_iy_i}\leq M_n.
  \]
  Note that $\bnorm{\sum_{i\geq n}a_iy_i}\kgeq (n\kminus
  1)/\epsilon$. Now by the choice of $(y_i)$ we have
  \[
  \Bnorm{\sum_{%
      \begin{subarray}{l}
        i\in E\\[1pt]
        i\geq n
    \end{subarray}}a_iy_i}\leq
  (L'+\epsilon)\Bnorm{\sum_{i\geq n}a_iy_i},
  \]
  and hence
  \begin{eqnarray*}
    \Bnorm{\sum_{i\in E}a_iy_i} &\leq&
    (n-1)+(L'+\epsilon)\Bnorm{\sum_{i\geq n}a_iy_i}\\
    &\leq& (L'+2\epsilon)\Bnorm{\sum_{i\geq n}a_iy_i}\\
    &\leq& (L'+2\epsilon)\big(\norm{x}+n\kminus 1\big)\\[5pt]
    &\leq& (L'+2\epsilon)(1+3\epsilon)\norm{x}.
  \end{eqnarray*}
  Indeed, $n\kminus 1\kleq\epsilon\bnorm{\sum_{i\geq n}a_iy_i}\kleq
  \epsilon(2\kplus\epsilon)\norm{x}\kleq 3\epsilon\norm{x}$ since
  $(y_i)$ is a basic sequence with constant $1\kplus\epsilon$. We have
  thus proved our claim from which $L(\delta)\kleq L'$ follows.
\end{proof}
To conclude this section we show that the various constants we
introduced remain the same if we restrict to the class of Banach
spaces $C(S)$, where $S$ is a countable, compact metric space. Recall
that such a space~$S$ is homeomorphic to a countable successor ordinal
in its order topology.
\begin{thm}
  For each $\delta\kin (0,1]$, we have
  \[
  K(\delta)=\sup_{\alpha,(x_i)}
  \underset{(y_i)\subset(x_i)}{\inf\phantom{p}} K((y_i),\delta),
  \]
  where the supremum is taken over all countable, successor ordinals
  $\alpha$ and all normalized, weakly null sequences $(x_i)$ in
  $C(\alpha)$, and the infimum is taken over all subsequences $(y_i)$
  of $(x_i)$. The analogous statements for the functions $K', L$ and
  $L'$ also hold.
\end{thm}
\begin{proof}
  We prove the result only for $K$. The argument for the other
  functions is similar. Fix $\epsilon\kin (0,1]$. By the definition of
  $K(\delta)$ there is a normalized, weakly null sequence $(x_i)$ in
  some Banach space such that $K((y_i),\delta)\kge
  K(\delta)\kminus\epsilon$ for every subsequence $(y_i)$ of
  $(x_i)$. Thus for each $M\kin\vegtelen{\bn}$ we have a triple
  $(\ba,E,x^{\ast})$ witnessing $K((x_i)_{i\in M},\delta)\kge
  K(\delta)\kminus\epsilon$, that is
  \begin{myeqnarray}
    \ba\keq (a_i)\kin\mc_{00},\quad E\ksubset\{i\kin
    M:\,\abs{a_i}\kgeq\delta\},\quad x^{\ast}\kin B_{X^{\ast}},\\
    \sup_{i\in M}\abs{a_i}\keq 1,\\
    \sum_{i\in E}a_ix^{\ast}(x_i)\kge
    (K(\delta)\kminus\epsilon)\norm{x},\quad
    \text{where}\ x\keq\sum_{i\in M}a_ix_i.
  \end{myeqnarray}
  We now use Lemma~\ref{lem:selection} to obtain a continuous
  selection $M\mapsto (\ba_M,E_M,x^{\ast}_M)$ of witnesses in the usual
  way. We set $\ba_M\keq(a^M_i)$ and $x_M\keq\sum_{i\in
  M}a^M_ix_i$.

  Next we pass to infinite subsets $N_1\ksupset N_2\ksupset N_3$ of
  $\bn$. First, there exists $N_1\kin\vegtelen{\bn}$ such that
  $(x_i)_{i\in N_1}$ is a basic sequence. Then we use
  Theorem~\ref{thm:schreier-unc;n-form} with $n\keq 1$ to find
  $N_2\kin\vegtelen{N_1}$ such that $\sup_{i\in N_2}\abs{a_i}\kleq
  (1\kplus\epsilon)\bnorm{\sum_{i\in N_2}a_ix_i}$ for all
  $(a_i)\kin\mc_{00}$. Note that in particular we have
  $\norm{x_M}\kgeq
  (1\kplus\epsilon)^{-1}$ for all $M\kin\vegtelen{N_2}$. Finally, by
  Lemma~\ref{lem:support-condition} there exists $N_3\kin\vegtelen{N_2}$
  such that $\sum_{i\in N_3\setminus P}\abs{x^{\ast}_P(x_i)}\kle\epsilon$
  for all $P\kin\vegtelen{N_3}$. 

  After relabelling, if necessary, we can take $N_3\keq \bn$. Set
  $X\keq[x_i]_{i=1}^{\infty}$, and let $(x^{\ast}_i)$ be the biorthogonal
  functionals to $(x_i)$. Note that $\norm{x^{\ast}_i}\kleq
  1\kplus\epsilon$ for all $i\kin\bn$ by the choice of $N_2$. For each
  $i\kin\bn$ and for $t\kin [-1,1]$ define $\rho_i(t)\keq
  \frac{\epsilon}{2^i}\big\lfloor \frac{2^it}{\epsilon}\big\rfloor
  \vee (-1)$, and note that $\abs{\rho_i(t)-t}\kleq\epsilon/2^i$.
  Now for each $M\kin\vegtelen{\bn}$ we have $\sum_{i\in\bn\setminus
  M}\abs{x^{\ast}_M(x_i)}\kle\epsilon$ by the choice of $N_3$, and
  $x^{\ast}_M\keq \sum_{i=1}^{\infty} x^{\ast}_M(x_i)x^{\ast}_i$ in the
  weak-$*$-sense since $(x_i)$ is a basis for $X$. It follows that
   \[
  \widetilde{x^{\ast}_M}=\sum_{i\in M}\rho_i(x^{\ast}_M(x_i))x^{\ast}_i,
  \]
  converges in the weak-$*$ sense, and moreover
  \[
  \bnorm{\widetilde{x^{\ast}_M}}\kleq 1\kplus
  2\epsilon(1\kplus\epsilon)\kleq 1\kplus 4\epsilon,\qquad\text{for
  all}\ M\kin\vegtelen{\bn}.
  \]
  Now define $S$ to be the closure of $U\keq\big\{
  \widetilde{x^{\ast}_M}:\,M\kin\vegtelen{\bn}\big\}\cup\{x^{\ast}_i
  :\,i\kin\bn\}$ in the weak-$*$ topology. Since $U$ is bounded in
  norm, $S$ is a compact metric space. The continuity of the choice of
  witnesses implies that $U$ is countable, and hence, because of the
  discretization of coefficients using the functions $\rho_i$, $S$ is
  also countable.

  Let $T\colon X\to C(S)$ be the canonical map, \ie
  $T(x)(y^{\ast})\keq y^{\ast}(x)$ for all $x\kin X,\ y^{\ast}\kin S$,
  and note that $\norm{T}\kleq 1\kplus 4\epsilon$. Set $f_i\keq
  T(x_i)$ for all $i\kin\bn$, and $f_M\keq T(x_M)$ for all
  $M\kin\vegtelen{\bn}$. Then $(f_i)$ is a normalized, weakly null
  sequence in $C(S)$. We claim that $K((f_i)_{i\in M},\delta)\kgeq
  \big(K(\delta)\kminus 3\epsilon\big)(1\kplus 4\epsilon)^{-1}$ for
  all $M\kin\vegtelen{\bn}$, which proves the assertion of the
  theorem.

  For $M\kin\vegtelen{\bn}$ we have $f_M\keq\sum_{i\in M}a^M_if_i$
  and $\norm{f_M}\kleq (1\kplus 4\epsilon)\norm{x_M}$. We also have
  $E_M\ksubset\{i\kin M:\,\abs{a^M_i}\kgeq\delta\}$, and
  \begin{eqnarray*}
    \Bnorm{\sum_{i\in E_M}a^M_if_i} &\geq& \sum_{i\in E_M}a^M_if_i
    \big(\widetilde{x^{\ast}_M}\big)=\sum_{i\in
      E_M}a^M_i\rho_i(x^{\ast}_M(x_i))\\
    &\geq& \sum_{i\in E_M}a^M_ix^{\ast}_M(x_i)-\epsilon\\
    &>& (K(\delta)-\epsilon)\norm{x_M}-\epsilon\geq
    (K(\delta)-3\epsilon)\norm{x_M}\\
    &\geq& \big(K(\delta)\kminus 3\epsilon\big)(1\kplus
    4\epsilon)^{-1} \norm{f_M},
  \end{eqnarray*}
  as required.
\end{proof}

\section{The $\mc_{0}$-problem}
\label{section:c_0-problem}

In this short section we consider the following intriguing
question which, to our knowledge, has not been raised elsewhere.
\begin{problem}
  \label{problem:c-0}
  Is there a real number~$C$ such that every sequence equivalent to
  the unit vector basis of $\mc_0$ has an unconditional subsequence
  with constant~$C$?
\end{problem}
Let $Y$ be the space $\mc_0$ or $\ell_p$ for some $p\kin[1,\infty)$,
and let $(e_i)$ be the unit vector basis of $Y$. Let $(x_i)$ be a
sequence in a Banach space equivalent to $(e_i)$. A well known result
of James~\cite{J} says that if $Y\keq\mc_0$ or $Y\keq\ell_1$, then for
any $\epsilon\kge 0$ there is a block basis of $(x_i)$ that is
$(1\kplus\epsilon)$-equivalent to $(e_i)$, and so in particular there
is a block basis of $(x_i)$ that is unconditional with constant
$(1\kplus\epsilon)$. Both these conclusions fail spectacularly if
$Y\keq\ell_p$ for some $p\kin(1,\infty)$: for any constant $C$ there
is an equivalent norm on $Y$ so that it contains no unconditional
basic sequence with constant $C$. This follows from the solution of
the distortion problem by Odell and Schlumprecht~\cite{OS}. For
$\mc_0$ and $\ell_1$ one can go further and consider subsequences
instead of block bases. However, if $Y\keq\mc_0$, then for any~$C$
there are easy examples that show that $(x_i)$ does not need to have a
subsequence $C$-equivalent to $(e_i)$. If $Y\keq\ell_1$, then for any
constant $C$ there are easy examples that show that $(x_i)$ does not
even need to have an unconditional subsequence with constant $C$. The
only remaining question in this context is raised in
Problem~\ref{problem:c-0}, which is still
open. Example~\ref{ex:elton>1} in Section~\ref{section:resolutions}
will show (among other things) that Problem~\ref{problem:c-0} cannot
have a positive answer with $C\kle 5/4$. However, it is possible that
a uniform constant $C$ exists. Indeed, this happens if and only if
$\sup_{\delta>0}L'(\delta)\kle\infty$, where $L'$ is the function
given in Definition~\ref{defn:constants}. Our aim in this section is
to prove this equivalence.

For each $\delta\kin (0,1]$ let us define $C(\delta)$ to be the
infimum of the set of real numbers~$C$ such that every normalized
sequence $1/\delta$-equivalent to the unit vector basis of $\mc_0$ has
an unconditional subsequence with constant $C$. So a positive answer
to Problem~\ref{problem:c-0} is equivalent to the statement that
$\sup_{\delta>0}C(\delta)$ is finite.
\begin{thm}
  \label{thm:c_0-iff-L}
  Let $\delta,\delta_1\kin(0,1]$.
  \begin{mylist}{(ii)}
  \item
    If $\delta_1\kleq\delta$, then $\ds C(\delta)\leq
    L(\delta_1)\cdot\Big(1+\frac{\delta_1}{\delta}\Big)+
    \frac{\delta_1}{\delta}$.
  \item
    If $\ds\delta_1\kle\frac{\delta}{2L'(\delta)}$, then
    $L'(\delta)\kleq C(\delta_1)$.
  \end{mylist}
  In particular we have
  $\sup_{\delta>0}C(\delta)\keq\sup_{\delta>0}L'(\delta)\keq
  \sup_{\delta>0}L(\delta)$.
\end{thm}
\begin{proof}
  To verify~(i) fix $\epsilon\kin (0,1]$, and assume that $(x_i)$ is a
  normalized sequence $1/\delta$-equivalent to the unit vector basis
  of $\mc_0$. So for some constants $A\kge0$ and $B\kge0$ with
  $B/A\kleq 1/\delta$ we have
  \begin{equation}
    \label{eq:c-0-equiv}
    A\sup_i\abs{a_i}\leq \Bnorm{\sum_{i=1}^{\infty}a_ix_i}\leq
    B\sup_i\abs{a_i}
  \end{equation}
  for all $(a_i)\kin\mc_{00}$. After passing to a subsequence, if
  necessary, we may assume that $L((x_i),\delta_1)\kleq
  L(\delta_1)\kplus\epsilon$. We claim that under these circumstances
  $(x_i)$ is unconditional with constant
  $C\keq\big(L(\delta_1)\kplus\epsilon\big)\Big(1\kplus
  \frac{\delta_1}{\delta}\Big)\kplus\frac{\delta_1}{\delta}$, from
  which~(i) follows.

  Given $(a_i)\kin\mc_{00}$ and $A\kin\veges{\bn}$, we need to show
  that $\bnorm{\sum_{i\in A}a_ix_i} \kleq C\norm{x}$, where
  $x\keq\sum_ia_ix_i$. We may clearly assume that $\sup_i\abs{a_i}\keq
  1$. Then it follows from~\eqref{eq:c-0-equiv} that
  \[
  B\delta\leq A=A\sup_i\abs{a_i}\leq\norm{x},
  \]
  and hence for every $F\ksubset\bn$ we have
  \begin{equation}
    \label{eq:c-0-equiv;small-coeff}
    \Bnorm{\sum_{i\in F}a_ix_i}\leq B\sup_{i\in F}\abs{a_i} \leq
    \frac{\norm{x}}{\delta}\,\sup_{i\in F}\abs{a_i}.
  \end{equation}
  Set $E\keq\{i\kin\bn:\,\abs{a_i}\kgeq\delta_1\}$. The definition of
  $L((x_i),\delta_1)$ gives
  \[
  \Bnorm{\sum_{i\in A\cap E}a_ix_i}\leq\big(L(\delta_1)+\epsilon\big)
  \Bnorm{\sum_{i\in E}a_ix_i}.
  \]
  Then from~\eqref{eq:c-0-equiv;small-coeff} we get
  \begin{eqnarray*}
    \Bnorm{\sum_{i\in E}a_ix_i} &\leq& \norm{x}+\Bnorm{\sum_{i\notin
	E}a_ix_i} \leq \Big( 1+\frac{\delta_1}{\delta}\Big)
    \norm{x},\\
    \Bnorm{\sum_{i\in A\setminus E}a_ix_i} &\leq&
    \frac{\delta_1}{\delta}\,\norm{x}.
  \end{eqnarray*}
  Finally, by the triangle-inequality we obtain
  \[
  \Bnorm{\sum_{i\in A}a_ix_i}\leq \bigg[\big(L(\delta_1)\kplus
  \epsilon\big)\Big(1\kplus\frac{\delta_1}{\delta}\Big)\kplus
  \frac{\delta_1}{\delta}\bigg]\cdot\norm{x}
  \]
  as required.

  We now prove~(ii). Fix $\epsilon\kge 0$. By the definition of
  $L'(\delta)$ there is a normalized, weakly null sequence $(x_i)$
  such that $L'((y_i),\delta)\kge L'(\delta)\kminus\epsilon$ for all
  subsequences $(y_i)$ of $(x_i)$. So for each $M\kin\vegtelen{\bn}$
  there is a triple $(\ba,E,x^{\ast})$ that we shall call a \emph{witness
  for~$M$}, where
  \begin{myeqnarray}
    \ba\keq(a_i)\kin\mc_{00},\quad E\ksubset\supp{\ba}\ksubset M,\quad
    x^{\ast}\kin B_{X^{\ast}};\\
    \label{eq:c_0-iff-L;witness2}
    \abs{a_i}\kgeq\delta\ \text{for all}\
    i\kin\supp{\ba},\quad\text{and}\quad \norm{x}\kleq 1,\
    \text{where}\ x\keq\sum_{i\in M}a_ix_i;\\
    \label{eq:c_0-iff-L;witness3}
    \sum_{i\in E}a_ix^{\ast}(x_i)>L'(\delta)-\epsilon.
  \end{myeqnarray}
  Let $\Omega$ be the set of all witnesses of all infinite subsets
  of~$\bn$, and for $M\kin\vegtelen{\bn}$ let $\Phi(M)$ be the
  (nonempty) set of all witnesses for~$M$. For $r\kin\bn$ let
  $\Omega_r$ be the set of all triples $(\ba,E,x^{\ast})\kin\Omega$ such
  that $\max\supp{\ba}\kleq r$. By Lemma~\ref{lem:selection} there is
  a function $\phi\colon\vegtelen{\bn}\to\Omega$ such that $\phi(M)$
  is a witness for~$M$ for all $M\kin\vegtelen{\bn}$, and~$\phi$ is
  continuous when~$\Omega$ is given the discrete topology. We let
  $\phi(M)=(\ba_M,E_M,x^{\ast}_M)$ and let
  \[
  A_M=\supp{\ba_M},\quad \ba_M=(a^M_i),\quad x_M=\sum_{i\in
  M}a^M_ix_i.
  \]
  By the proof of Lemma~\ref{lem:selection} we can choose~$\phi$ so
  that for all $M\kin\vegtelen{\bn}$ there exists $r\kin\bn$ such that
  $\Phi(M)\cap\Omega_s\keq\emptyset$ whenever $1\kleq s\kle r$, and
  $\phi(M)$ is the least element of $\Phi(M)\cap\Omega_r$ in some
  well-ordering of $\Omega$ fixed in advance. It is then easy to
  verify that for all $L,M\kin\vegtelen{\bn}$ if $A_M\kprec L$, then
  $\phi(L)\keq\phi(M)$.

  We now pass to infinite subsets $N_1\ksupset N_2\ksupset N_3\ksupset
  N_4$ of~$\bn$. Let $f\colon \vegtelen{\bn}\to\mc_0$ be the function
  that maps $M\kin\vegtelen{\bn}$ to $f_M\keq
  \big(x^{\ast}_M(x_i)\big)\kin\mc_0$. It follows from the continuity
  of~$\phi$ that~$f$ is continuous with respect to the topology of
  pointwise convergence on $\mc_0$ and that its image has compact
  closure. Hence by Lemma~\ref{lem:support-condition} there exists
  $N_1\kin\vegtelen{\bn}$ such that
  \begin{equation}
    \label{eq:c_0-iff-L;support-condition}
    \sum_{i\in N_1\setminus P}\abs{x^{\ast}_P(x_i)}<\epsilon\qquad
    \text{for all}\ P\in\vegtelen{N_1}.
  \end{equation}
  Now choose arbitrary $N_2\kin\vegtelen{N_1}$ with $N_1\ksetminus
  N_2$ of infinite size. Given $M\kin\vegtelen{N_2}$, we can choose
  $L\kin\vegtelen{N_1}$ such that $A_M\kprec L$ and $L\ksetminus
  A_M\ksubset N_1\ksetminus N_2$. Then $\phi(L)\keq\phi(M)$, and
  applying~\eqref{eq:c_0-iff-L;support-condition} with $P\keq L$ we
  obtain
  \begin{equation}
    \label{eq:c_0-iff-L;support-condition2}
    \sum_{i\in N_2\setminus A_M}\abs{x^{\ast}_M(x_i)}= \sum_{i\in
    N_2\setminus A_M}\abs{x^{\ast}_L(x_i)}<\epsilon
  \end{equation}
  since $N_2\ksetminus A_M\ksubset N_1\ksetminus L$. In other words,
  relative to~$N_2$ and up to a small error, we have
  $\supp{x^{\ast}_M}\ksubset A_M$ for all $M\kin\vegtelen{N_2}$.

  By the definition of~$L'(\delta)$ there exists
  $N_3\kin\vegtelen{N_2}$ such that $L'((x_i)_{i\in N_3},\delta)\kleq
  L'(\delta)\kplus\epsilon$. Finally, we apply
  Theorem~\ref{thm:schreier-unc;n-form} with $n\keq 1$ to obtain
  $N_4\kin\vegtelen{N_3}$ such that for all $M\kin\vegtelen{N_4}$ we
  have $\abs{a^M_i}\kleq 1\kplus\epsilon$ for all $i\kin M$.

  We now relabel so that we can take $N_4\keq\bn$, and define a new
  norm on $\mc_0$ by setting
  \[
  ||| \bb |||=\norm{\bb}_{\mc_0}\vee \sup_{M\in\vegtelen{\bn}}
  \Babs{\sum_{i=1}^{\infty}b_ix^{\ast}_M(x_i)}\qquad \text{for}\
  \bb\keq(b_i)\kin\mc _0.
  \]
  Let $(y_i)$ be the unit vector basis of $\mc_0$ considered with its
  new norm. It
  follows from~\eqref{eq:c_0-iff-L;support-condition2}
  and the choice of $N_3$ that
  \[
  \delta\norm{x^{\ast}_M}_{\ell _1}\leq \sum_{i\in
  A_M}\babs{a^M_i}\babs{x^{\ast}_M(x_i)}+\epsilon\delta \leq
  2\big(L'(\delta)+\epsilon\big)+\epsilon\delta
  \]
  for all $M\kin\vegtelen{\bn}$. Hence $(y_i)$ is $D$-equivalent to
  $(e_i)$, where
  \[
  D=\frac{2\big(L'(\delta)+\epsilon\big)+\epsilon\delta}{\delta}
  <\frac{1}{\delta_1}
  \]
  provided $\epsilon$ is sufficiently small. We claim that $(y_i)$ has
  no unconditional subsequence with constant
  $C\keq\big(L'(\delta)\kminus\epsilon\big)/(1\kplus\epsilon)$. Fix
  $M\kin\vegtelen{\bn}$. We have
  \[
  \Babs{\sum_{i\in M}a^M_ix^{\ast}_L(x_i)}=
  \abs{x^{\ast}_L(x_M)}\leq 1\qquad \text{for all}\
  L\kin\vegtelen{\bn},
  \]
  and hence, by the choice of $N_4$, we have $||| \ba_M |||\kleq
  1\kplus\epsilon$. On the other hand,
  property~\eqref{eq:c_0-iff-L;witness3} of a witness applied to~$M$
  gives
  \[
  \Big|\Big|\Big| \sum_{i\in E_M}a^M_iy_i \Big|\Big|\Big|\geq
  \sum_{i\in E_M} a^M_ix^{\ast}_M(x_i)>L'(\delta)-\epsilon \geq C||| \ba_M
  |||,
  \]
  which shows the claim. Since $\epsilon$ was arbitrary,
  $C(\delta_1)\kgeq L'(\delta)$ follows.

  Parts (i) and (ii) show that $\sup_{\delta>0}L'(\delta)\kleq
  \sup_{\delta>0}C(\delta)\kleq\sup_{\delta>0}L(\delta)$. That we have
  equality throughout follows from
  Proposition~\ref{prop:near-unc-constants}.
\end{proof}

\section{Convex-unconditionality and duality}
\label{section:convex-unc}

The following notion of partial unconditionality was introduced by
Argyros, Mercourakis, and Tsarpalias~\cite{AMT}. Given $\delta\kin
(0,1]$, we say that a basic sequence $(x_i)$ is
\emph{$\delta$-convex-unconditional with constant~$A$} if for all
$(a_i)\kin\mc_{00}$ and for all $E\kin\veges{\bn}$ if
\[
\delta \sum _{i\in E} \lvert a_i\rvert \leq \Big\lVert \sum _{i\in E}
a_ix_i \Big\rVert,
\]
then we have
\[
\Big\lVert \sum _{i\in E} a_ix_i \Big\rVert \leq A \Big\lVert \sum
_{i=1}^{\infty} a_ix_i \Big\rVert.
\]
The definition in~\cite{AMT} is actually slightly different, but it is
equivalent to ours (they express unconditionality in terms of
sign-changes rather than projections). Theorem~\ref{thm:ell_1-unc}
on $\ell_1$-projections follows immediately from the next result.
\begin{thm}[Argyros, Mercourakis, and Tsarpalias~\cite{AMT}]
  \label{thm:convex-unc}
  Given \mbox{$\delta\kin (0,1]$} there is a constant $A$ such that
  every normalized weakly null sequence has a
  $\delta$-convex-unconditional subsequence with constant
  $A$. Moreover, $A\kleq 16\log_2\big(1/\delta\big)$ for
  $\delta\kle1/4$.
\end{thm}
\begin{proof}
  Given $\delta\kin (0,1]$, define
  $l\keq\big\lfloor\log_2\big(1/\delta\big)\big\rfloor\kplus 2$ and
  fix $A\kin[1,\infty)$. Assume that $(x_i)$ is a normalized, weakly
  null sequence, which has no $\delta$-convex-unconditional
  subsequence with constant $A$. We will show that $A\kleq 8l$.
  
  Without loss of generality $(x_i)$ is a basic sequence with
  constant~$2$, say. So for all $(a_i)\kin\mc_{00}$ we have
  \begin{equation}
    \label{eq:conv-unc;basis-const}
    \sup_i\abs{a_i}\leq 4\Bnorm{\sum_{i=1}^{\infty} a_ix_i}.
  \end{equation}
  Let $M\kin\vegtelen{\bn}$. Since $(x_i)_{i\in M}$ is not
  $\delta$-convex-unconditional with constant~$A$, there exist
  $(a_i)\kin\mc_{00}$ and $E\kin\veges{M}$ such that
  \[
  \delta\sum_{i\in E}\abs{a_i} \leq \Bnorm{\sum_{i\in E}a_ix_i},\quad
  \text{and}\quad A\norm{x}<\Bnorm{\sum_{i\in E}a_ix_i},
  \]
  where $x\keq\sum_{i\in M}a_ix_i$. Rescaling, considering
  appropriate subsets of $E$, and replacing $(a_i)$ by $(-a_i)$ if
  necessary, we conclude that for every $M\kin\vegtelen{\bn}$ there
  exists a quadruple $(\ba,F,x^{\ast},k)$, called a \emph{witness for~$M$},
  with the following properties.
  \begin{myeqnarray}
    \label{eq:conv-unc;witness1}
    \ba\keq(a_i)\kin\mc_{00},\quad F\kin\veges{M},\quad
    x^{\ast}\kin B_{X^{\ast}},\quad k\kin\{1,\dots,l\},\\[3pt]
    \label{eq:conv-unc;witness2}
    a_i\kge 0\quad \text{and}\quad 2^{-k}\kle x^{\ast}(x_i)\kleq
    2^{-k+1}\quad \text{for all}\ i\kin F,\\[3pt]
    \label{eq:conv-unc;witness3}
    \frac{\delta}{4l}\leq \sum_{i\in F} a_ix^{\ast}(x_i),\\[3pt]
    \label{eq:conv-unc;witness4}
    A\norm{x}< 2l\sum_{i\in F} a_ix^{\ast}(x_i)+\frac{\delta}{2},\qquad
    \text{where}\ x=\sum _{i\in M} a_ix_i.
  \end{myeqnarray}
  We now use Lemma~\ref{lem:selection} in the usual way to select a
  witness $(\ba_M,F_M,x^{\ast}_M,k_M)$ for each $M\kin\vegtelen{\bn}$ in a
  continuous way, where the set of all witnesses is given the discrete
  topology. We write $\ba_M\keq (a^M_i)$ and $x_M\keq
  \sum_{i\in M}a^M_ix_i$.

  We now carry out stabilizations. Fix $\epsilon\kge 0$, and pass to
  an infnite subset $N$ of $\bn$ such that for all $P\in \vegtelen{N}$
  we have
  \begin{equation}
    \label{eq:conv-unc;supp-condition}
    \sum_{i\in N\setminus P} \babs{x^{\ast}_P(x_i)}\leq \epsilon,
  \end{equation}
  and for all $L,M\kin\vegtelen{N}$ we have $k_L\keq k_M$. The first
  property is achieved by Lemma~\ref{lem:support-condition}, whereas
  the second uses infinite Ramsey theory. Observe that for all $L,
  M\kin\vegtelen{N}$ we have
  \begin{equation}
    \label{eq:const-functional}
    x^{\ast}_L(x_i)\geq \frac{1}{2} x^{\ast}_M(x_i)\qquad
    \text{whenever}\ i\in F_L\cap F_M.
  \end{equation}
  We finally apply Theorem~\ref{thm:matching} with $n\keq 1$ to find
  $L,M\kin\vegtelen{N}$ such that
  \begin{equation}
    \label{eq:conv-unc;matching}
    L\cap M=F_L\cap F_M=F_M.
  \end{equation}
  We now estimate $x^{\ast}_L(x_M)$. On the one hand,
  using~\eqref{eq:conv-unc;matching} followed
  by~\eqref{eq:const-functional}, \eqref{eq:conv-unc;basis-const},
  and~\eqref{eq:conv-unc;supp-condition}, we have
  \begin{multline*}
    x^{\ast}_L(x_M)=\sum_{i\in M}a^M_ix^{\ast}_L(x_i)=\sum_{i\in F_M}a^M_i
  x^{\ast}_L(x_i) + \sum_{i\in M\setminus L} a^M_ix^{\ast}_L(x_i)\\[2pt]
  \geq \frac{1}{2}\sum_{i\in F_M} a^M_i x^{\ast}_M(x_i) - 4\epsilon
  \norm{x_M}.
  \end{multline*}
  On the other hand, property~\eqref{eq:conv-unc;witness4} applied to
  the witness of~$M$ gives
  \[
  x^{\ast}_L(x_M)\leq\norm{x_M} < \frac{2l}{A}\sum_{i\in F_M}
  a^M_ix^{\ast}_M(x_i) +\frac{\delta}{2A}.
  \]
  The last two inequalities together with
  property~\eqref{eq:conv-unc;witness3} of the witness of~$M$ show
  that
  \[
  \left(\frac{1}{2}-(1+4\epsilon)\frac{2l}{A}\right)\frac{\delta}{4l}
  \leq (1+4\epsilon)\frac{\delta}{2A}.
  \]
  Since~$\epsilon$ was arbitrary, it follows that $A\kleq 8l$, as claimed.
\end{proof}
Given a normalized, weakly null sequence $(x_i)$ and $\delta\kin
(0,1]$ let $A((x_i),\delta)$ be the least real number~$A$ such that
$(x_i)$ is $\delta$-convex-unconditional with constant~$A$. Then
define
\[
A(\delta)=\sup _{(x_i)}\underset{(y_i)\subset
  (x_i)}{\inf\phantom{p}}A((y_i),\delta),
\]
where the supremum is taken over all normalized, weakly null sequences
$(x_i)$ and the infimum over all subsequences $(y_i)$ of $(x_i)$.
Theorem~\ref{thm:convex-unc} yields an upper bound of order
$\log\big(1/\delta\big)$ on $A(\delta)$. We are now going to prove
that the question whether $\sup_{\delta}A(\delta)\kle\infty$ is
equivalent to Problem~\ref{problem:near-unc} using the function $K'$
defined on page~\pageref{defn:constants}. As the proof shows the two
problems are in some sense dual to each other.
\begin{thm}
  \label{thm:A-iff-K}
  For $0\kle\delta_1\kle\delta\kleq 1$ we have
  \begin{mylist}{(ii)}
  \item
    $A(\delta)\leq\frac{\delta}{\delta-\delta_1}K'(\delta_1)$, and
  \item
    $K'(\delta)\leq A(\delta_1)$.
  \end{mylist}
  In particular $\sup_{\delta>0}A(\delta)\keq\sup_{\delta>0}
  K'(\delta)$.
\end{thm}
\begin{proof}
  We begin by proving~(i). Fix $\epsilon\kin (0,1]$. There is a
  normalized, weakly null sequence $(x_i)$ such that
  $A((y_i),\delta)\kge A(\delta)\kminus\epsilon$ for every subsequence
  $(y_i)$ of $(x_i)$. On the other hand, after passing to a
  subsequence if necessary, we may assume that $A((x_i),\delta)\kleq
  A(\delta)\kplus\epsilon$. Set $C\keq
  A(\delta)(\delta\kminus\delta_1)/\delta
  -\epsilon(\delta\kplus\delta_1)/\delta$. For each
  $M\kin\vegtelen{\bn}$ there is a triple $(\ba,x^{\ast},F)$, called a
  \emph{witness for $M$}, where
  \begin{myeqnarray}
    \label{eq:A-iff-K;witness1}
    \ba\keq (a_i)\kin\mc_{00},\quad x^{\ast}\kin B_{X^{\ast}},\quad
    F\ksubset\{i\kin M:\,\abs{x^{\ast}(x_i)}\kgeq\delta_1\},\\
    \label{eq:A-iff-K;witness2}
    \sum_{i\in F}a_ix^{\ast}(x_i)\kge C\norm{x},\quad\text{where}\
    x\keq\sum_{i\in M}a_ix_i.
  \end{myeqnarray}
  Indeed, since $A((x_i)_{i\in M},\delta)\kge
  A(\delta)\kminus\epsilon$, there exist $\ba\keq (a_i)\kin\mc_{00}$
  and $E\kin\veges{M}$ such that $\delta\sum_{i\in E}\abs{a_i}\kleq
  \bnorm{\sum_{i\in E}a_ix_i}$, and $\bnorm{\sum_{i\in E}a_ix_i}\kge
  (A(\delta)\kminus\epsilon)\bnorm{x}$, where $x\keq\sum_{i\in
  M}a_ix_i$. By homogeneity, we may assume that $\sum_{i\in
  E}\abs{a_i}\keq 1$. Let $x^{\ast}\kin B_{X^{\ast}}$ be a support
  functional for $\sum_{i\in E}a_ix_i$, and let $F\keq\{i\kin
  E:\,\abs{x^{\ast}(x_i)}\kgeq\delta_1\}$. An easy computation now shows
  that~\eqref{eq:A-iff-K;witness2} holds.

  We now use Lemma~\ref{lem:selection} in the usual way to obtain a
  continuous selection $M\mapsto (\ba_M,x^{\ast}_M,F_M)$ of witnesses. We
  let $\ba_M\keq (a^M_i)$ and $x_M\keq \sum_{i\in M}a^M_ix_i$.

  Next we find infinite subsets $N_1\ksupset N_2\ksupset N_3$ of $\bn$
  as follows. First, there exists $N_1\kin\vegtelen{\bn}$ such that
  $(x_i)_{i\in N_1}$ is a basic sequence with constant
  $1\kplus\epsilon$. Then we apply
  Theorem~\ref{thm:schreier-unc;n-form} with
  $n\keq 1$ to get $N_2\kin\vegtelen{N_1}$ such that $\abs{a^M_i}\kleq
  (1\kplus\epsilon)\norm{x_M}$ for all $M\kin\vegtelen{N_2}$ and for
  all $i\kin M$. Finally, by Lemma~\ref{lem:support-condition} there
  exists $N_3\kin\vegtelen{N_2}$ such that $\sum_{i\in N_3\setminus
  M}\abs{x^{\ast}_M(x_i)}\kle\epsilon$ for all $M\kin\vegtelen{N_3}$.

  After relabelling we may assume that $N_3\keq\bn$. Let $(e_i)$ be
  the unit vector basis of $\mc_{00}$, and for each
  $M\kin\vegtelen{\bn}$ set
  \[
  t_M=\frac{1}{(1+\epsilon)^2\norm{x_M}} \sum_{i\in M}a^M_ie_i,
  \]
  which is an element of $[-1,1]^{\bn}$ by the choice of $N_2$. We
  endow $[-1,1]^{\bn}$ with the product topology and let $S$ be the
  closure of $\{t_M:\, M\in\vegtelen{\bn}\} \cup\{e_i:\,i\in\bn\}$.
  Note that $S$ is a compact metric space. For each $i\kin\bn$ let
  $f_i$ be the $i^{\text{th}}$ co-ordinate map. By the continuity of
  the choice of witnesses, $S$ contains only sequences of finite
  support. Hence $(f_i)$ is a normalized, weakly null sequence in
  $C(S)$. We claim that $K'((f_i)_{i\in M},\delta_1)\kgeq
  C/(1\kplus\epsilon)^2$ for all $M\kin\vegtelen{\bn}$. Since
  $\epsilon$ was arbitrary, (i) follows from this claim.

  Given $M\kin\vegtelen{\bn}$, set $n_M\keq\max F_M$, and let
  \[
  f_M=\sum_{%
    \begin{subarray}{l}
      i\in M\\[1pt]
      i\leq n_M
    \end{subarray}}x^{\ast}_M(x_i)f_i.
  \]
  For each $L\kin\vegtelen{\bn}$ we have
  \begin{eqnarray*}
    (1+\epsilon)^2\norm{x_L}\abs{f_M(t_L)} &=&
    \Babs{\sum_{%
        \begin{subarray}{l}
          i\in L\cap M\\[1pt]
          i\leq n_M
        \end{subarray}
      }a^L_ix^{\ast}_M(x_i)
    }\\
    &\leq& \Babs{x^{\ast}_M\Big(\sum_{%
        \begin{subarray}{l}
          i\in L\\[1pt]
          i\leq n_M
        \end{subarray}
      }a^L_ix_i \Big)
    }+\epsilon(1+\epsilon)\norm{x_L}
    \leq (1+\epsilon)^2\norm{x_L}
  \end{eqnarray*}
  by the choices of $N_1,N_2$ and $N_3$. It follows that
  $\norm{f_M}\kleq 1$. On the other hand, we have $F_M\ksubset\{i\kin
  M:\,i\kleq n_M,\ \abs{x^{\ast}_M(x_i)}\kgeq\delta_1\}$ and
  \begin{eqnarray*}
  \Bnorm{\sum_{i\in F_M}x^{\ast}_M(x_i)f_i} &\geq&
  \sum_{i\in F_M}x^{\ast}_M(x_i)f_i(t_M)\\
  &=& \frac{1}{(1+\epsilon)^2\norm{x_M}}\sum_{i\in
  F_M}x^{\ast}_M(x_i)a^M_i>\frac{C}{(1+\epsilon)^2},
  \end{eqnarray*}
  as claimed.

  To show~(ii) fix $\epsilon\kin (0,1]$ so that $\delta_1(1\kplus
  3\epsilon)\kle\delta$, and let $(x_i)$ be a
  normalized, weakly null sequence such that $K'((y_i),\delta)\kge
  K'(\delta)\kminus\epsilon$ for every subsequence $(y_i)$ of
  $(x_i)$. So for each $M\kin\vegtelen{\bn}$ there is a triple
  $(\ba,E,x^{\ast})$, called a \emph{witness for $M$}, where
  \begin{myeqnarray}
    \ba\keq (a_i)\kin\mc_{00},\quad E\ksubset\{i\kin
    M:\,\abs{a_i}\kgeq\delta\},\quad x^{\ast}\kin B_{X^{\ast}},\\
    \norm{x}\keq 1,\ \text{where}\ x\keq \sum_{i\in M}a_ix_i,\\
    \sum_{i\in E}a_ix^{\ast}(x_i)>K'(\delta)-\epsilon,\quad \text{and}\
    a_ix^{\ast}(x_i)\kge 0\ \text{for all}\ i\kin E.
  \end{myeqnarray}
  As usual, we then have a continuous choice $M\mapsto
  (\ba_M,E_M,x^{\ast}_M)$ of witnesses, and we let $\ba_M\keq (a^M_i)$
  and $x_M\keq\sum_{i\in M}a^M_ix_i$.

  We now pass to infinite subsets $N_1\ksupset N_2\ksupset N_3$ of
  $\bn$. First, we choose $N_1\kin\vegtelen{\bn}$ so that $(x_i)_{i\in
  N_1}$ is a basic sequence with constant $1\kplus\epsilon$. Then we
  apply Theorem~\ref{thm:schreier-unc;n-form} with $n\keq 1$ to
  find $N_2\kin\vegtelen{N_1}$ such that we have $\abs{a^M_i}\kleq
  (1\kplus\epsilon)\norm{x_M}$ for all $M\kin\vegtelen{N_2}$ and for
  all $i\kin M$. Finally we use Lemma~\ref{lem:support-condition} in
  the usual way to obtain $N_3\kin\vegtelen{N_2}$ so that $\sum_{i\in
  N_3\setminus M}\abs{x^{\ast}_M(x_i)}\kle\epsilon$ for all
  $M\kin\vegtelen{N_3}$.

  We now relabel so that we can take $N_3\keq\bn$. As before, we let
  $(e_i)$ be the unit vector basis of $\mc_{00}$. We define $S$ to be
  the closure in $[-1,1]^{\bn}$ of the set
  \[
  \Big\{\frac{1}{1+3\epsilon} \sum_{i\in M}a^M_ie_i:\,
  M\in\vegtelen{\bn}\Big\} \cup\{e_i:\,i\in\bn\}.
  \]
  As before, it is easy to verify that $S$ is a compact metric
  space containing only sequences of finite support, and that the
  sequence $(f_i)$ of co-ordinate maps is a normalized, weakly null
  sequence in $C(S)$. We now show that $A((f_i)_{i\in M},\delta_1)\kgeq
  (K'(\delta)\kminus\epsilon)(1\kplus3\epsilon)^{-1}$ for all
  $M\kin\vegtelen{\bn}$.

  Given $M\kin\vegtelen{\bn}$, set $n_M\keq\max E_M$, and let
  \[
  f_M=\sum_{%
    \begin{subarray}{l}
      i\in M\\[1pt]
      i\leq n_M
    \end{subarray}}x^{\ast}_M(x_i)f_i.
  \]
  For each $L\kin\vegtelen{\bn}$ we have
  \[
  \Babs{\sum_{%
      \begin{subarray}{l}
        i\in L\cap M\\[1pt]
        i\leq n_M
      \end{subarray}
    }a^L_ix^{\ast}_M(x_i)
  }\leq \Babs{x^{\ast}_M\Big(\sum_{%
      \begin{subarray}{l}
        i\in L\\[1pt]
        i\leq n_M
      \end{subarray}
    }a^L_ix_i \Big)
  }+\epsilon(1+\epsilon)\leq 1+3\epsilon
  \]
  by the choices of $N_1,N_2$ and $N_3$. It follows that
  $\norm{f_M}\kleq 1$. On the other hand, we have
  \[
  \Bnorm{\sum_{i\in E_M}x^{\ast}_M(x_i)f_i}\geq
  \frac{1}{1+3\epsilon}\sum_{i\in
    E_M}x^{\ast}_M(x_i)a^M_i\geq \delta_1\sum_{i\in
    E_M}\abs{x^{\ast}_M(x_i)},
  \]
  as well as
  \[
  \Bnorm{\sum_{i\in E_M}x^{\ast}_M(x_i)f_i}\geq
  \frac{1}{1+3\epsilon}\sum_{i\in
    E_M}x^{\ast}_M(x_i)a^M_i>\frac{K'(\delta)-\epsilon}{1+3\epsilon}\,
  \norm{f_M},
  \]
  which proves the claim.
\end{proof}

\section{Unconditionality in $C(S)$ spaces and duality}
\label{section:C(S)-spaces}

We now turn to questions on finding unconditional basic sequences in
spaces of continuous functions on a compact, Hausdorff space. We will
then relate these to problems considered so far. We start by stating a
result of Rosenthal.
\begin{thm}
  \label{thm:wns-of-indicator-fns}
  For any compact, Hausdorff space~$S$, every weakly null sequence of
  (non-zero) indicator functions in $C(S)$ has an unconditional
  subsequence with constant~$1$.
\end{thm}
In~\cite{AGR} this is presented as a consequence of a combinatorial
lemma. Here we prove a more general version of that, and obtain a more
general version of Theorem~\ref{thm:wns-of-indicator-fns}. Before
stating it we need some notation. Let
$k\kin\bn$ and $M\kin\vegtelen{\bn}$. Given $\ba\keq (a_i)_{i\in M}$
and $\bb\keq (b_i)_{i\in M}$ in $\{0,1,\dots,k\}^M$, we write
$\ba\ksubset\bb$ if for all $i\kin M$ either $a_i\keq 0$ or $a_i\keq
b_i$. Given $j\kin\{1,\dots,k\}$, we write $\ba\ksubset_j \bb$ if for
all $i\kin M$ either $a_i\keq 0$ or $a_i\keq b_i\keq j$. A family
$\mathcal F\ksubset\{0,1,\dots,k\}^M$ is \emph{hereditary} if
$\ba\kin\mathcal F$ whenever $\bb\kin\mathcal F$ and $\ba\ksubset\bb$,
and is \emph{weakly hereditary} if $\ba\kin\mathcal F$ whenever
$\bb\kin\mathcal F$ and there exists $j\kin\{1,\dots,k\}$ such that
$\ba\ksubset_j\bb$. Given $L\kin\vegtelen{M}$, we denote by $\mathcal
F_L$ the set of restrictions to~$L$ of elements of $\mathcal F$. Note
that $\mathcal F_L\ksubset\{0,1,\dots,k\}^L$.
\begin{lemma}
  \label{lem:dual-matching}
  Let $k\kin\bn$ and $\mathcal F\ksubset\{0,1,\dots,k\}^{\bn}$ be a
  compact family of sequences of finite support. Then there exists
  $M\kin\vegtelen{\bn}$ such that $\mathcal F_M$ is weakly
  hereditary.
\end{lemma}
\begin{proof}
  We argue by contradiction. Assuming that the statement is false, for
  each $M\kin\vegtelen{\bn}$ we can find a quadruple $(\ba,\bb,j,K)$
  that we shall call a \emph{witness for $M$}, where
  \begin{myeqnarray}
    \label{eq:dual-matching;witness1}
    \ba\kin\{0,1,\dots,k\}^{\bn},\quad \bb\kin\mathcal F,\quad
    j\kin\{1,\dots,k\},\quad K\kin\bn;\\
    \label{eq:dual-matching;witness2}
    \supp{\ba}\ksubset M,\quad\ba\ksubset_j\bb,\quad
    K\kge\max\supp{\bb};\\
    \text{if $a_i\keq c_i$ for all $i\kin M$ with $i\kleq K$,
    then}\ \bc\knotin\mathcal F.
  \end{myeqnarray}
  Indeed, the assumption that $\mathcal F_M$ is not weakly hereditary
  implies the existence of $\ba,\bb,j$ as
  in~\eqref{eq:dual-matching;witness1} such that $\supp{\ba}\ksubset
  M,\ \ba\ksubset_j\bb$ and there is no $\bc\kin\mathcal F$ such that
  the restrictions to $M$ of $\ba$ and $\bc$ are identical. The
  existence of a suitable~$K$ now follows easily from the compactness
  of~$\mathcal F$.

  Let $\Omega$ denote the set of all witnesses of all infinite subsets
  of~$\bn$. For $r\kin\bn$ let $\Omega _r$ be the set of elements
  $(\ba,\bb,j,K)\kin\Omega$ for which $K\kleq r$. The conditions of
  Lemma~\ref{lem:selection} are now easily verified (which is why we
  needed to introduce the parameter~$K$). So there is a continuous
  selection $\phi\colon\vegtelen{\bn}\to\Omega$ of witnesses. Let
  $\phi(M)\keq(\ba_M,\bb_M,j_M,K_M)$, where $\ba_M\keq(a^M_i)$ and
  $\bb_M\keq(b^M_i)$ for each $M\kin\vegtelen{\bn}$.

  The continuity of~$\phi$ and the compactness of~$\mathcal F$ imply
  that the function $M\mapsto\bb_M\colon\vegtelen{\bn}\to\mc_0$ is
  continuous and its image has compact closure (in the topology of
  pointwise convergence). So applying
  Lemma~\ref{lem:support-condition} with $\epsilon\keq 1/2$, say, we
  find $N_1\kin\vegtelen{\bn}$ such that
  \begin{equation}
    \label{eq:dual-matching;support-condition}
    N_1\cap\supp{\bb_L}\ksubset L\qquad \text{for all}\
    L\kin\vegtelen{N_1}.
  \end{equation}
  An easy application of infinite Ramsey theory then gives
  $j\kin\{1,\dots,k\}$ and $N_2\kin\vegtelen{N_1}$ such that $j_M\keq
  j$ for all $M\kin\vegtelen{N_2}$.
  
  To conclude the proof we apply the Matching Lemma with $n\keq 1$ to
  the function $M\mapsto\supp{\ba_M}$ to find $L,M\kin\vegtelen{N_2}$
  such that
  \[
  L\cap M=\supp{\ba_L}\cap\supp{\ba_M}=\supp{\ba_M}.
  \]
  Now if $i\kin\supp{\ba_M}$, then $a^M_i\keq a^L_i\keq b^L_i\keq j$
  by property~\eqref{eq:dual-matching;witness2} of a witness and by
  the choice of~$N_2$. On the other hand, if $i\kin
  M\ksetminus\supp{\ba_M}$, then $i\knotin L$, and hence
  by~\eqref{eq:dual-matching;support-condition} we have
  $i\knotin\supp{\bb_L}$, so $a^M_i\keq b^L_i\keq 0$. We have shown
  that the restrictions to $M$ of $\ba_M$ and the element $\bb_L$ of
  $\mathcal F$ are identical which gives the required contradiction.
\end{proof}

\begin{thm}
  \label{thm:dual-L}
  For all $\delta\kin(0,1]$ there is a constant~$L^{\ast}$ such that for
  any compact, Hausdorff space~$S$, if $(f_i)$ is a normalized, weakly
  null sequence in $C(S)$ with $\abs{f_i(t)}\kin\{0\}\cup[\delta,1]$
  for all $t\kin S$ and $i\kin\bn$, then $(f_i)$ has an unconditional
  subsequence with constant~$L^{\ast}$. Moreover, $L^{\ast}\kleq
  6\log_2\big(1/\delta\big)$ for $\delta\kle1/4$.
\end{thm}
\begin{proof}
  For $\delta\kin(0,1]$ let $k\keq\lfloor\log_2(1/\delta)\rfloor
  \kplus 1$. Let $I_0\keq\{0\}$ and let $I_1,\dots,I_k$ be closed
  intervals covering $[\delta,1]$ such that $\max I_j\kleq 2\min I_j$
  for each $j\keq 1,\dots ,k$. Furthermore, let $I_{j+k}\keq\kminus
  I_j$ for $j\keq 1,\dots,k$. Let~$S$ be a compact, Hausdorff space
  and $(f_i)$ be a normalized, weakly null sequence in $C(S)$ with
  $\abs{f_i(t)}\kin\{0\}\cup[\delta,1]$ for all $t\kin S$ and
  $i\kin\bn$. Let $\mathcal F$ be the collection of all
  $\bc\kin\{0,1,\dots,2k\}^{\bn}$ for which there exists $t\kin S$
  with $f_i(t)\kin I_{c_i}$ for all $i\kin\bn$. Note that  $\mathcal
  F$ is a compact subset of $\{0,1,\dots,2k\}^{\bn}$ consisting of
  sequences of finite support. By Lemma~\ref{lem:dual-matching} there
  exists $M\kin\vegtelen{\bn}$ such that $\mathcal F_M$ is weakly
  hereditary. We show that the sequence $(f_i)_{i\in M}$ is
  unconditional with constant $L^{\ast}\keq 4k$.

  Fix $\ba\keq (a_i)\kin\mc_{00}$ and
  $E\kin\veges{M}$. Choose $t\kin S$ such that
  \[
  \Bnorm{\sum_{i\in E}a_if_i}\keq\Babs{\sum_{i\in E}a_if_i(t)}.
  \]
  Replacing $\ba$ by $-\ba$ if necessary, we may assume that
  \[
  \Babs{\sum_{i\in E}a_if_i(t)}\leq\sum_{i\in F}a_if_i(t),
  \]
  where $F\keq\{i\kin E:\,a_if_i(t)>0\}$. Now choose $\bc\kin\mathcal
  F$ such that $f_i(t)\kin I_{c_i}$ for all $i\kin\bn$. Note that
  $c_i\kneq 0$ for any $i\kin F$, and so
  \[
  \sum_{i\in F}a_if_i(t)\leq 2k\sum_{i\in F_j}a_if_i(t)
  \]
  for some $j\kin\{1,\dots,2k\}$, where $F_j\keq\{i\kin F:\,c_i\keq
  j\}$. Finally, since $\mathcal F_M$ is weakly hereditary, there
  exists $\bc'\kin\mathcal F$ such that $c'_i\keq c_i\keq j$ for all
  $i\kin F_j$, and $c'_i\keq 0$ for all $i\kin M\ksetminus F_j$. Let
  $t'\kin S$ satisfy $f_i(t')\kin I_{c'_i}$ for all $i\kin\bn$. We
  then have
  \[
  \sum_{i\in F_j}a_if_i(t)\leq 2\sum_{i\in M}a_if_i(t')\leq
  2\Bnorm{\sum_{i\in M}a_if_i}.
  \]
  This completes the proof of our claim.
\end{proof}
\begin{rems}
  1.~If $(f_i)$ is a weakly null sequence of (non-zero) indicator
  functions, then in the proof above we need only to work with two
  intervals $I_0\keq\{0\}$ and $I_1\keq\{1\}$. This way  we do not get
  the factor of~2 at either of the two places where it occurs above,
  and so we obtain a proof of
  Theorem~\ref{thm:wns-of-indicator-fns}.
  We also mention here a quantitative version of Rosenthal's result
  due to Gasparis, Odell and Wahl~\cite{GOW}: if $(f_i)$ is a weakly
  null sequence of (non-zero) indicator functions, then there exists a
  countable ordinal $\alpha$ and a subsequence $(g_i)$ of $(f_i)$
  which is equivalent to a subsequence of the unit vector basis of the
  generalized Schreier space $X^{\alpha}$.

  \noindent
  2.~Lemma~\ref{lem:dual-matching} and Theorem~\ref{thm:dual-L} were
    also proved by Arvanitakis (he uses slightly different language
    and method). In ~\cite[Remark~2.1]{A} he effectively asks if
  weakly hereditary can be replaced by hereditary in
  Lemma~\ref{lem:dual-matching}. It is not hard to see that if that
  was possible, then the proof of Theorem~\ref{thm:dual-L} would give
  a constant $L^{\ast}$ independent of~$\delta$. In turn, by
  Theorem~\ref{thm:duality} below, this would yield a positive
  solution to the $\mc_0$-problem. The following simple example shows
  that Lemma~\ref{lem:dual-matching} cannot be strengthened in this
  way even for $k\keq 2$. For each $M\keq\{m_1\kle
  m_2\kle\dots\}\kin\vegtelen{\bn}$ define $c_M\kin\{0,1,2\}^{\bn}$ by
  letting
  \begin{equation*}
    c_M(m_i)=
    \begin{cases}
      2 & \text{if $i\keq 1$}\\
      1 & \text{if $2\kleq i\kleq m_1\kplus 1$}\\
      2 & \text{if $m_1\kplus 1\kle i\kleq m_2\kplus 1$}
    \end{cases}
  \end{equation*}
  and $c_M$ is zero elsewhere. Now let $\mathcal F$ be the set
  of all $c\kin\{0,1,2\}^{\bn}$ such that there exist
  $M\kin\vegtelen{\bn}$ and $n\kin\bn$ such that $c(i)\keq c_M(i)$ for
  $i\keq 1,\dots, n$ and $c(i)\keq 0$ for all $i\kge n$ --- we denote
  this $c$ by $c_{M,n}$. Then $\mathcal F$ is a compact family of
  sequences of finite support. To see that $\mathcal F_L$ is not
  hereditary for any $L\keq\{l_1\kle l_2\kle\dots\}\kin\vegtelen{\bn}$
  consider $c,c'\kin\{0,1,2\}^L$ defined as follows: $c'(l_1)\keq 0,\
  c(l_1)\keq c_L(l_1)$ and $c'(l_i)\keq c(l_i)\keq c_L(l_i)$ for all
  $i\kgeq 2$. Then $c\kin\mathcal F_L,\ c'\ksubset c$, but given
  $M\keq\{m_1\kle m_2\kle\dots\}\kin\vegtelen{\bn}$, there is no
  $n\kin\bn$ such that $c'$ is the restriction to $L$ of $c_{M,n}$
  (consider the cases $m_2\kle l_2,\ m_2\keq l_2$ and $m_2\kge l_2$).
\end{rems}
We now prove a more general result of which Theorem~\ref{thm:dual-L}
is an immediate consequence.
\begin{thm}
  \label{thm:dual-K}
  For all $\delta\kin(0,1]$ there is a constant $K^{\ast}$ such that
  for any compact, Hausdorff space~$S$, every normalized, weakly null
  sequence $(f_i)$ in $C(S)$ has a subsequence $(g_i)$ such that for
  all $t\kin S$ and $E\ksubset\{i\kin\bn:\,\abs{g_i(t)}\kgeq\delta\}$
  we have
  \[
  \Babs{\sum_{i\in E}a_ig_i(t)}\leq
  K^{\ast}\Bnorm{\sum_{i=1}^{\infty}a_ig_i}\qquad \text{for all}\
  (a_i)\kin\mc_{00}.
  \]
  Moreover, $K^{\ast}\kleq 6\log_2\big(1/\delta\big)$ for $\delta\kle1/4$.
\end{thm}
\begin{proof}
  Fix $\delta\kin(0,1]$ and $K^{\ast}\kin[1,\infty)$. Assume that~$S$
  is a compact, Hausdorff space, and $(f_i)$ is a normalized, weakly
  null sequence in $C(S)$ that has no subsequence satisfying the
  statement of the theorem. We will show that $K^{\ast}\kleq 4k$,
  where $k\keq\lfloor\log_2(1/\delta)\rfloor\kplus 1$.

  Let $I_1,\dots,I_k$ be closed intervals covering $[\delta,1]$ such
  that $\max I_j\kleq 2\min I_j$ for each $j\keq 1,\dots
  ,k$. Furthermore, let $I_{j+k}\keq\kminus I_j$ for $j\keq
  1,\dots,k$. For every $M\kin\vegtelen{\bn}$ there is a
  \emph{witness} $(t,\ba,j,F)$ to the failure of the subsequence
  $(f_i)_{i\in M}$, where
  \begin{myeqnarray}
    t\kin S,\quad \ba\keq(a_i)\kin\mc_{00},\quad
    j\kin\{1,\dots,2k\},\phtm{\sum_{i=1}}\newline
    F\ksubset\{i\kin M:\,f_i(t)\kin I_j,\ a_if_i(t)\kge 0\};\\
    \norm{f}\keq 1,\ \text{where}\ f\keq\sum_{i\in M}a_if_i;\\
    2k\sum_{i\in F}a_if_i(t)\kge K^{\ast}.
  \end{myeqnarray}
  We now use Lemma~\ref{lem:selection} to get a continuous selection
  $M\mapsto (t_M,\ba_M,j_M,F_M)$ of witnesses. Let $\ba_M\keq
  (a^M_i)$ and $f_M\keq\sum_{i\in M}a^M_if_i$ for each
  $M\kin\vegtelen{\bn}$.

  As usual, the next phase of the proof is stabilization. Find
  $N_1\kin\vegtelen{\bn}$ such that $(f_i)_{i\in N_1}$ is a basic
  sequence with constant~$2$, and so $\abs{a^M_i}\kleq 4$ for all
  $i\kin M$ and for all $M\kin\vegtelen{N_1}$. Then pass to
  $N_2\kin\vegtelen{N_1}$ such that $j_L\keq j_M$ for all
  $L,M\kin\vegtelen{N_2}$, which in particular implies that
  $f_i(t_M)$ and $f_i(t_L)$ have the same sign, and differ by a factor
  of at most~$2$ for all $i\kin F_L\cap F_M$. Finally, we fix
  $\epsilon\kge 0$ and use
  Lemma~\ref{lem:support-condition} to obtain $N_3\kin\vegtelen{N_2}$
  such that $\sum_{i\in N_3\setminus P}\abs{f_i(t_P)}\kle \epsilon$
  for all $P\kin\vegtelen{N_3}$.

  The Matching Lemma applied with $n\keq 1$ now yields
  $L,M\kin\vegtelen{N_3}$ such that $L\cap M\keq F_L\cap F_M\keq
  F_M$. Then
  \begin{eqnarray*}
    \abs{f_M(t_L)} &=& \Babs{\sum_{i\in M}a^M_if_i(t_L)}\\
    &\geq& \sum_{i\in F_M}a^M_if_i(t_L) -4\epsilon\\
    &\geq& \frac{1}{2}\sum_{i\in F_M}a^M_if_i(t_M) -4\epsilon\\
    &\geq& \frac{K}{4k}-4\epsilon.
  \end{eqnarray*}
  On the other hand, $\abs{f_M(t_L)}\kleq\norm{f_M}\keq 1$, and hence
  $K\kleq 4k(1+4\epsilon)$.
\end{proof}
We will now establish a relationship between Theorem~\ref{thm:dual-L},
which is a result about finding unconditional subsequences, and the
constant $L'$ (defined on page~\pageref{defn:constants}), which comes
from a certain form of partial unconditionality. We will also show the
close connection between Theorem~\ref{thm:dual-K} and
Problem~\ref{problem:near-unc}. First we need to introduce some
appropriate constants, and then we will express these relationships in
Theorem~\ref{thm:duality} below.

For a basic sequence $(x_i)$ in a Banach space let $C(x_i)$ be the
least real number $C$ such that $(x_i)$ is unconditional with
constant~$C$. Then for each $\delta\kin(0,1]$ we define
\[
L^{\ast}(\delta)=\sup_{S,(f_i)}\underset{(g_i)\subset(f_i)}{\inf
  \phantom{p}} C(g_i),
\]
where the supremum is taken over all compact, Hausdorff spaces~$S$ and
over all normalized, weakly null sequences $(f_i)$ in $C(S)$ with
$\abs{f_i(t)}\kin\{0\}\cup[\delta,1]$ for all $t\kin S$ and
$i\kin\bn$, and the infimum is taken over subsequences $(g_i)$ of
$(f_i)$. Theorem~\ref{thm:dual-L} above claims that $L^{\ast}(\delta)$ is
finite and of order $\log\big(1/\delta\big)$.

Given $\delta\kin (0,1]$, and a normalized, weakly null sequence
$(f_i)$ in $C(S)$ with $S$ a compact, Hausdorff space, we define
$K^{\ast}((f_i),\delta)$ to be the least real number~$K^{\ast}$ such that
whenever $t\kin S$ and $E\ksubset
\{i\kin\bn:\,\abs{f_i(t)}\kgeq\delta\}$, we have
\[
\Babs{\sum_{i\in E}a_if_i(t)}\leq
K^{\ast}\Bnorm{\sum_{i=1}^{\infty}a_if_i}
\qquad \text{for all}\ (a_i)\kin\mc_{00}.
\]
We then set
\[
K^{\ast}(\delta)=\sup_{S,(f_i)}\underset{(g_i)\subset
  (f_i)}{\inf\phantom{p}} K^{\ast}((g_i),\delta),
\]
where the supremum is over all compact, Hausdorff spaces $S$ and all
normalized, weakly null sequences $(f_i)$ in $C(S)$, and the infimum
is over all subsequences $(g_i)$ of $(f_i)$. Note that by
Theorem~\ref{thm:dual-K} above $K^{\ast}(\delta)$ is finite and of
order $\log\big(1/\delta\big)$.
\begin{thm}
  \label{thm:duality}
  For all $0\kle\delta'\kle\delta\kleq 1$ we have
  $K^{\ast}(\delta)\kleq K'(\delta)\kleq K^{\ast}(\delta')$ and
  $L^{\ast}(\delta)\kleq L'(\delta)\kleq L^{\ast}(\delta')$.
\end{thm}
\begin{proof}
  We first show that $K^{\ast}(\delta)\kleq K'(\delta)$. Fix
  $\epsilon\kin (0,1]$. There is a compact Hausdorff space~$S$ and a
  normalized, weakly null sequence $(f_i)$ in $C(S)$ such that
  $K^{\ast}((f_i)_{i\in M},\delta)\kge
  K^{\ast}(\delta)\kminus\epsilon$ for all $M\kin\vegtelen{\bn}$. So
  for each $M\kin\vegtelen{\bn}$ there is a \emph{witness} $(t,E,\ba)$
  for~$M$, where
  \begin{myeqnarray}
    t\kin S,\quad E\ksubset\{i\kin M:\,\abs{f_i(t)}\kgeq\delta\},\quad
    \ba\keq(a_i)\kin\mc_{00};\\
    \norm{f}\keq 1,\quad \text{where}\ f\keq \sum_{i\in M}a_if_i;\\
    \Babs{\sum_{i\in E}a_if_i(t)}> K^{\ast}(\delta)-\epsilon.
  \end{myeqnarray}
  We now proceed as usual. We make a continuous choice
  $M\mapsto (t_M,E_M,\ba_M)$ of witnesses, and let
  $\ba_M\keq(a^M_i)$ and $f_M\keq\sum_{i\in M}a^M_if_i$.
  We then find $N_1\kin\vegtelen{\bn}$ such that $(f_i)_{i\in N_1}$ is
  a basic sequence with constant $1\kplus\epsilon$. By
  Theorem~\ref{thm:schreier-unc;n-form} there exists
  $N_2\kin\vegtelen{N_1}$ such that $\abs{a^M_i}\kleq 1\kplus\epsilon$
  for all $i\kin M$ and for all $M\kin\vegtelen{N_2}$. Finally, we
  pass to a further infinite subset $N_3$ of $N_2$ such that $\sum
  _{i\in N_3\setminus P}\abs{f_i(t_P)}\kle \epsilon$ for all
  $P\kin\vegtelen{N_3}$.

  After relabelling, if necessary, we may assume that $N_3\keq\bn$. We
  define a norm on $\mc_{00}$ by letting
  \[
  \norm{(b_i)}=\sup_{i}\abs{b_i} \vee \frac{1}{(1+\epsilon)^2}\sup
  \Big\{\Babs{\sum_{i\in M}b_ia^M_i}:\,M\in\vegtelen{\bn}\Big\}
  \]
  for each $(b_i)\kin\mc_{00}$. Let $X$ be the completion of the
  resulting normed space. It is easy to check that the unit vector
  basis $(e_i)$ of $\mc_{00}$ is a normalized, weakly null sequence
  in~$X$. Indeed, the continuity of the selection of witnesses implies
  that the closure of $\big\{\sum_{i\in M}a^M_ie_i:\,M\kin
  \vegtelen{\bn}\big\}\cup\{e_i:\,i\kin\bn\}$ in the topology of
  pointwise convergence contains only finitely supported sequences. We
  will now show that $K'((y_i),\delta)\kge
  K^{\ast}(\delta)\kminus\epsilon$ for any subsequence $(y_i)$ of
  $(e_i)$. This then proves the inequality $K^{\ast}(\delta)\kleq
  K'(\delta)$.

  Fix $M\kin\vegtelen{\bn}$, and set $n_M\keq \max E_M$ and
  \[
  x_M=\sum _{%
    \begin{subarray}{l}
      i\in M\\[1pt]
      i\leq n_M
    \end{subarray}} f_i(t_M)e_i.
  \]
  For each $L\kin\vegtelen{\bn}$ we have
  \begin{eqnarray*}
    \Babs{\sum _{%
    \begin{subarray}{l}
      i\in L\cap M\\[1pt]
      i\leq n_M
    \end{subarray}} f_i(t_M)a^L_i} &\leq& \Babs{\sum _{%
    \begin{subarray}{l}
      i\in L\\[1pt]
      i\leq n_M
    \end{subarray}} f_i(t_M)a^L_i} + (1\kplus\epsilon) \sum_{i\in
      L\setminus M}\abs{f_i(t_M)}\\
    &\leq& (1\kplus\epsilon)\norm{f_L}+(1\kplus\epsilon)\epsilon =
    (1\kplus\epsilon)^2.
  \end{eqnarray*}
  It follows that $\norm{x_M}\kleq 1$. On the other hand, on $E_M$ the
  coefficients of $x_M$ are at least~$\delta$ and
  \[
  \Bnorm{\sum_{i\in E_M} f_i(t_M)e_i} \geq \Babs{\sum_{i\in E_M}
  f_i(t_M)a^M_i} >K^{\ast}(\delta)-\epsilon.
  \]
  We now show that $K'(\delta)\kleq K^{\ast}(\delta')$ whenever $0\kle
  \delta'\kle\delta\kleq 1$. Fix $\epsilon\kin (0,1]$ such that
  $(1\kplus\epsilon)\delta'\kle\delta$. Let $(x_i)$ be a normalized,
  weakly null sequence with $K'((y_i),\delta)\kge
  K'(\delta)\kminus\epsilon$ for every subsequence $(y_i)$ of
  $(x_i)$. So for each $M\kin\vegtelen{\bn}$ there is a \emph{witness}
  $(\ba,E,x^{\ast})$ \emph{of} $M$, where
  \begin{myeqnarray}
    \ba\keq (a_i)\kin\mc_{00},\quad E\ksubset\{i\kin
    M:\,\abs{a_i}\kgeq\delta\},\quad x^{\ast}\kin B_{X^{\ast}};\\
    \norm{x}\keq 1,\quad \text{where}\ x\keq\sum_{i\in M}a_ix_i;\\
    \sum_{i\in E}a_ix^{\ast}(x_i) > K'(\delta)-\epsilon.
  \end{myeqnarray}
  Let $M\mapsto (\ba_M,E_M,x^{\ast}_M)$ be a continuous selection of
  witnesses, and let $\ba_M\keq(a^M_i)$ and
  $x_M\keq\sum_{i\in M}a^M_ix_i$. Choose $N_1\kin\vegtelen{\bn}$ such
  that $(x_i)_{i\in N_1}$ is a basic sequence with constant
  $1\kplus\epsilon$. Use Theorem~\ref{thm:schreier-unc;n-form} to find
  $N_2\kin\vegtelen{N_1}$ so that $\abs{a^M_i}\kleq 1\kplus\epsilon$
  for all $i\kin M$ and $M\kin\vegtelen{N_2}$. Finally, by
  Lemma~\ref{lem:support-condition} there exists
  $N_3\kin\vegtelen{N_2}$ such that $\sum_{i\in N_3\setminus
  P}\abs{x^{\ast}_P(x_i)}\kle\epsilon$ for all $P\kin\vegtelen{N_3}$.

  Relabel so that we can take $N_3\keq \bn$, and set
  $t_M\keq\frac{1}{1\kplus\epsilon}\sum_{i\in M}a^M_ie_i$ for each
  $M\kin\vegtelen{\bn}$, where $(e_i)$ is the unit vector basis of
  $\mc_{00}$. Let $S$ be the closure of the set
  $\{t_M:\,M\kin\vegtelen{\bn}\}\cup\{e_i:\,i\kin\bn\}$ in the product
  space $[-1,1]^{\bn}$. As before, it is easy to verify that $S$
  consists only of finitely supported sequences, and hence the
  sequence $(f_i)$ of co-ordinate maps is a normalized, weakly null
  sequence in $C(S)$. We will show that
  \[
  K^{\ast}((g_i),\delta')>
  \frac{K'(\delta)-\epsilon}{(1+\epsilon)^2}
  \]
  for every subsequence $(g_i)$ of $(f_i)$, which then implies that
  $K^{\ast}(\delta')\kgeq K'(\delta)$.

  Fix $M\kin\vegtelen{\bn}$ and let $n_M\keq\max\supp{\ba_M}$ and
  \[
  f_M=\sum_{%
    \begin{subarray}{l}
      i\in M\\[1pt]
      i\leq n_M
    \end{subarray}} x^{\ast}_M(x_i)f_i.
  \]
  For each $L\kin\vegtelen{\bn}$ we have
  \begin{eqnarray*}
    (1+\epsilon)\abs{f_M(t_L)} &=& \Babs{\sum_{%
    \begin{subarray}{l}
      i\in L\cap M\\[1pt]
      i\leq n_M
    \end{subarray}} x^{\ast}_M(x_i)a^L_i}\\
    &\leq& \Babs{\sum_{%
    \begin{subarray}{l}
      i\in L\\[1pt]
      i\leq n_M
    \end{subarray}} x^{\ast}_M(x_i)a^L_i}+(1\kplus\epsilon)\sum_{i\in
      L\setminus M}\abs{x^{\ast}_M(x_i)}\\
    &\leq& (1\kplus\epsilon)\norm{x_L}+(1\kplus\epsilon)\epsilon\leq
    (1\kplus\epsilon)^2.
  \end{eqnarray*}
  It follows that $\norm{f_M}\kleq 1\kplus\epsilon$. On the other
  hand, we have
  \[
  \abs{f_i(t_M)}=\abs{a^M_i}/(1\kplus\epsilon)>\delta'\qquad\text{for
  all}\ i\kin E_M,
  \]
  and moreover
  \[
  \sum_{%
    \begin{subarray}{l}
      i\in E_M\\[1pt]
      i\leq n_M
  \end{subarray}}x^{\ast}_M(x_i)f_i(t_M)=\sum_{i\in
  E_M}x^{\ast}_M(x_i)a^M_i/(1\kplus\epsilon) >
  \frac{K'(\delta)-\epsilon}{(1+\epsilon)^2} \norm{f_M}.
  \]
  This completes the proof of the inequalities involving $K'$ and
  $K^{\ast}$. The argument for the functions $L'$ and $L^{\ast}$ is
  similar and is omitted.
\end{proof}
Recall that if $(x_i)$ is a normalized, weakly null sequence with
spreading model not equivalent to the unit vector basis of $\mc_0$,
then for any $\epsilon\kge0$ and for any $\delta\kin (0,1]$ there is a
$\delta$-near-unconditional subsequence of $(x_i)$ with constant
$1\kplus\epsilon$. There are dual versions of this corresponding to
Theorems~\ref{thm:dual-L} and~\ref{thm:dual-K} above. For example, for
any compact, Hausdorff space~$S$ and for any $\delta\kin(0,1]$, if
$(f_i)$ is a normalized, weakly null sequence in $C(S)$ with
$\abs{f_i(t)}\kin\{0\}\cup[\delta,1]$ for all $t\kin S$ and 
$i\kin\bn$, and $(f_i)$ has spreading model not equivalent to the unit
vector basis of $\ell_1$, then for any $\epsilon\kge0$ there is a
subsequence of $(f_i)$ that is unconditional with
constant~$1\kplus\epsilon$. The proof (which we omit here) uses a
similar argument to that of \cite[Theorem~5.4]{DKK}.

\section{The combinatorics of patterns and resolutions}
\label{section:resolutions}

In this section we consider combinatorial structures that arise in our
approach to Problem~\ref{problem:near-unc}. We begin by setting up
witnesses for the constant $K'(\delta)$ (\cf
Definition~\ref{defn:constants}). The notation will be used throughout
this section. We fix $\delta\kin(0,1]$, set $k\keq\lfloor\log_2(1/
\delta)\rfloor\kplus 1$, and choose $\epsilon\kin(0,1)$ so that
$2^k\delta\kge 1\kplus\epsilon$. We then select closed intervals
$I_1,\dots,I_k$ covering $[\delta,1\kplus\epsilon]$ so that $\max
I_j\kleq 2\min I_j$ for each $j\keq 1,\dots,k$. By the definition of
$K'(\delta)$ there is a normalized, weakly null sequence $(x_i)$ in
some Banach space $X$ such that $K'((y_i),\delta)\kge
\frac{1}{2}\,K'(\delta)$ for every subsequence $(y_i)$ of
$(x_i)$. After passing to a subsequence if necessary we can assume, as
usual, that
\begin{myeqnarray}
  \label{eq:patterns;basic}
  (x_i)\ \text{is a basic sequence with constant}\ 1\kplus\epsilon,\\
  \label{eq:patterns;schreier}
  \sup_i\abs{a_i}\kleq(1\kplus\epsilon)\Bnorm{\sum_{i=1}^{\infty}
  a_ix_i}\ \text{for all}\ (a_i)\kin\mc_{00}.
\end{myeqnarray}
Recall that the latter property is achieved using
Theorem~\ref{thm:schreier-unc;n-form}. We now make a continuous
selection $M\mapsto (\ba_M,x_M^{\ast},F_M)$ of witnesses in the usual
manner using Lemma~\ref{lem:selection}, where
\begin{myeqnarray}
  \label{eq:patterns;witness1}
  \ba_M\keq(a^M_i)\kin\mc_{00},\quad x^{\ast}_M\kin
  B_{X^{\ast}},\quad F_M\ksubset\{i\kin M:\,a^M_i\kgeq\delta,\
  x^{\ast}_M(x_i)\kge 0\},\\
  \norm{x_M}\keq 1,\ \text{where}\ x_M\keq\sum_{i\in M}a^M_ix_i,\\
  \label{eq:patterns;witness3}
  \sum_{i\in F_M}a^M_ix^{\ast}_M(x_i)> K'(\delta)/4.
\end{myeqnarray}
Note that $\abs{a^M_i}\kin [\delta,1\kplus\epsilon]$ for all
$M\kin\vegtelen{\bn}$ and for all $i\kin F_M$. For each
$M\kin\vegtelen{\bn}$ let us define
$\bc_M\keq(c^M_i)\kin\{0,1,\dots,k\}^{\bn}$ by letting $c^M_i$ be the
least $j\kin\{1,\dots,k\}$ such that $a^M_i\kin I_j$ if $i\kin F_M$,
and letting $c^M_i\keq 0$ otherwise. Set
$F^M_j\keq\{i\kin\bn:\,c^M_i\keq j\}$ for each $j\keq1,\dots,k$. Note
that
\begin{myeqnarray}
  \label{eq:patterns;colour-class}
  F^M_1,\dots,F^M_k\ \text{are pairwise disjoint, finite subsets of}\
  M\ \text{with}\ F_M\keq \bigcup_{j=1}^k F^M_j,\\
  \label{eq:patterns;continuity}
  \text{for each}\ j\keq 1,\ldots,k\ \text{the function}\ M\mapsto
  F^M_j\colon\vegtelen{\bn}\to \veges{\bn}\ \text{is continuous}.
\end{myeqnarray}
Note that we have $\osc{\ba_M}{F^M_j}\kleq 2$ for all
$M\kin\vegtelen{\bn}$ and for each $j\keq 1,\ldots,k$. Moreover, for
any two infinite subsets $L, M$ of $\bn$ we have
\begin{equation}
  \frac{a^M_i}{a^L_i}\geq\frac{1}{2}\qquad \text{for all}\ i\in
  F^L_j\cap F^M_j,\ j=1,\ldots,k.
\end{equation}
Using the usual Ramsey type arguments
(Lemma~\ref{lem:support-condition} and the infinite Ramsey theorem)
and relabeling, if necessary, we may assume the following
stabilizations.
\begin{myeqnarray}
  \label{eq:patterns;support-condition}
  \sum_{i\notin M}\abs{x^{\ast}_M(x_i)}<\epsilon\ \text{for all}\
  M\kin\vegtelen{\bn};\\
  \label{eq:patterns;weights}
  \text{for each}\ j\keq 1,\ldots,k\ \text{there exists}\ w_j\
  \text{such that for all}\ M\kin\vegtelen{\bn}\ \text{we
  have}\newline
  w_j\kleq\sum_{i\in F^M_j}a^M_ix^{\ast}_M(x_i)\kleq
  w_j\kplus\epsilon/k.
\end{myeqnarray}
Observe that~\eqref{eq:patterns;witness3}
and~\eqref{eq:patterns;weights} give
\begin{myeqnarray}
  \label{eq:patterns;total-weight}
  \sum_{j=1}^k w_j>\frac{K'(\delta)}{4}-\epsilon.
\end{myeqnarray}
We now give a simple necessary and sufficient condition for a positive
answer to Problem~\ref{problem:near-unc}.
\begin{prop}
  \label{prop:matching}
  We have $\sup_{\delta>0}K'(\delta)\kle\infty$ if and only if there
  is a constant $c$ such that for all $\delta\kin(0,1]$ whenever
  $(x_i)$ is a normalized, weakly null sequence in a Banach space and
  $M\mapsto (\ba_M,x_M^{\ast},F_M)$ is a continuous selection of witnesses
  so that~\eqref{eq:patterns;basic}--\eqref{eq:patterns;witness3}
  and~\eqref{eq:patterns;support-condition} hold, then there exist
  infinite subsets $L,M$ of $\bn$ such that $\abs{x^{\ast}_L(x_M)}\kgeq
  cK'(\delta)$.
\end{prop}
\begin{proof}
  Sufficiency is clear: for any $\delta\kin(0,1]$ there is a
  normalized, weakly null sequence $(x_i)$ and a continuous selection
  $M\mapsto (\ba_M,x_M^{\ast},F_M)$ of witnesses so
  that~\eqref{eq:patterns;basic}--\eqref{eq:patterns;witness3}
  and~\eqref{eq:patterns;support-condition} hold. The assumption then
  gives $K'(\delta)\kleq 1/c$.

  Now assume that $K'\keq\sup_{\delta >0}K'(\delta)$ is finite. We
  show that the condition is necessary with $c\keq\frac{1}{16K'}$. Let
  $\delta\kin(0,1]$ and assume that we are given a normalized, weakly
  null sequence $(x_i)$ and a continuous selection $M\mapsto
  (\ba_M,x_M^{\ast},F_M)$ of witnesses so
  that~\eqref{eq:patterns;basic}--\eqref{eq:patterns;witness3}
  and~\eqref{eq:patterns;support-condition} hold. Set
  $t\keq\frac{K'(\delta)}{16K'}$. For $\bb\keq(b_i)\kin\mc_{00}$
  define
  \[
  \tnorm{\bb}=t\norm{\bb}_{\ell_{\infty}}\vee
  \sup\bigg\{\Babs{\sum_{i=1}^{\infty}b_ix^{\ast}_L(x_i)}:\,
  L\kin\vegtelen{\bn}\bigg\}.
  \]
  Let $Z$ be the completion of $\mc_{00}$ in the norm
  $\tnorm{\cdot}$. The unit vector basis $(e_i)$ of $\mc_{00}$ is a
  semi-normalized, weakly null sequence of $Z$. So by the definition
  of $K'(t\delta)$ there is an infinite subset $M$ of $\bn$ such that
  $K'((\hat{e}_i)_{i\in M},t\delta)<2K'(t\delta)$, say, where
  $\hat{e}_i\keq\frac{e_i}{\tnorm{e_i}}$ for all
  $i\kin\bn$. From~\eqref{eq:patterns;witness3} we get
  \[
  \Btnorm{\sum_{i\in F_M}a^M_ie_i}\geq \Babs{\sum_{i\in
  F_M}a^M_ix^{\ast}_M(x_i)} > K'(\delta)/4.
  \]
  Now let $\bb_M\keq\sum_{i\in
  M}a^M_ie_i$. By~\eqref{eq:patterns;schreier} we have
  $\norm{\bb_M}_{\ell_{\infty}}\kleq 1\kplus\epsilon$, which in turn
  gives $\tnorm{\bb_M}\kleq 1$ since $x^{\ast}_L$ has norm at most one
  for all $L\kin\vegtelen{\bn}$. Now since
  $\babs{a^M_i\tnorm{e_i}}\kgeq t\delta$ for each $i\kin F_M$, we have
  \[
  \Btnorm{\sum_{i\in F_M}a^M_ie_i} =\Btnorm{\sum_{i\in
  F_M}a^M_i\tnorm{e_i}\hat{e}_i} \leq 2K'(t\delta)\tnorm{\bb_M}.
  \]
  We can now conclude that
  \[
  \tnorm{\bb_M}>\frac{K'(\delta)}{8K'(t\delta)}>
  t\norm{\bb_M}_{\ell_{\infty}}.
  \]
  Hence there exists $L\kin\vegtelen{\bn}$ such that
  \[
  \abs{x^{\ast}_L(x_M)}=\Babs{\sum_{i\in M}a^M_ix^{\ast}_L(x_i)}>
  \frac{1}{2}\tnorm{\bb_M}\geq cK'(\delta).
  \]
\end{proof}
A selection of witnesses as defined
in~\eqref{eq:patterns;witness1}--\eqref{eq:patterns;witness3}
associates to each $M\kin\vegtelen{\bn}$ a certain combinatorial data
that is made up of two parts. One part is the sequence $(c^M_i)_{i\in
  F_M}$ in the set $\{1,\ldots,k\}$, which is a discretized version of
the coefficients $(a^M_i)_{i\in F_M}$ of the vector
$x_M$. Equivalently, this part can also be viewed as the partition
$(F^M_j)_{j=1}^k$ of $F_M$. The other part is the sequence
$\big(x^{\ast}_M(x_i)\big)_{i\in F_M}$ of dual coefficients. To solve
Problem~\ref{problem:near-unc} in the affirmative we would like to
show the existence of $L,M\kin\vegtelen{\bn}$ whose combinatorial data
``match'' in a suitable way to give the necessary and sufficient
condition of Proposition~\ref{prop:matching}. For example, if we could
assume that the sets $F^M_1,\dots,F^M_k$ are successive for all
$M\kin\vegtelen{\bn}$, then the Matching Lemma would provide suitable
sets $L$ and $M$. Indeed, we can generalize this as follows.
\begin{prop}
  \label{prop:pure-matching}
  The following is a sufficient condition for
  $\sup_{\delta>0}K'(\delta)\kle\infty$. There exists a constant $c$
  such that for all $k\kin\bn$ and for all positive real numbers
  $p_1,\ldots,p_k$ with $\sum_{j=1}^kp_j\keq 1$ if for all
  $M\kin\vegtelen{\bn}$ we are given finite subsets
  $F^M_1,\ldots,F^M_k$ of $M$ such
  that~\eqref{eq:patterns;colour-class}
  and~\eqref{eq:patterns;continuity} hold, then there exist $L,
  M\kin\vegtelen{\bn}$ and $J\ksubset\{1,\ldots,k\}$ such that
  $\sum_{j\in J}p_j\kgeq c,\ F^L_j\ksubset F^M_j$ for all $j\kin J$,
  and $L\cap M\ksubset F_L\cap F_M$.
\end{prop}
\begin{rem}
  The Matching Lemma implies that the above sufficient condition is
  satisfied with $c\keq\frac{1}{2}$ provided that we also
  require $F^M_1\kle\ldots\kle F^M_k$ for all $M\kin\vegtelen{\bn}$.
\end{rem}
\begin{proof}
  We will verify that the stated condition implies the sufficient and
  necessary condition of Proposition~\ref{prop:matching}. Given
  $\delta\kin(0,1]$, assume that we are given a normalized, weakly
  null sequence $(x_i)$, a continuous selection $M\mapsto
  (\ba_M,x_M^{\ast},F_M)$ of witnesses so
  that~\eqref{eq:patterns;basic}--\eqref{eq:patterns;witness3}
  and~\eqref{eq:patterns;support-condition} hold. After passing to a
  subsequence, if necessary, we may assume that $\epsilon\kle c/48$
  and all of the
  conditions~\eqref{eq:patterns;basic}--\eqref{eq:patterns;weights}
  hold. Let $w\keq\sum_{j=1}^kw_j$, and set $p_j\keq w_j/w$ for each
  $j\keq1,\ldots,k$. By our assumption we can find $L,
  M\kin\vegtelen{\bn}$ and $J\ksubset \{1,\ldots,k\}$ such that
  $\sum_{j\in J}p_j\kgeq c,\ F^L_j\ksubset F^M_j$ for all $j\kin J$,
  and $L\cap M\ksubset F_L\cap F_M$. Note that $\sum_{j\in J}w_j\kgeq
  cw\kgeq cK'(\delta)/4\kminus\epsilon$. We now obtain a sequence of
  inequalities in a way very similar to that at the end of the proof
  of Theorem~\ref{thm:bdd-osc-unc}.
\begin{eqnarray*}
  x^{\ast}_L(x_M) & \geq & \sum _{i\in L\cap M} a^M_i x^{\ast}_L(x_i) -
  2\epsilon\\
  & \geq & \frac{1}{2} \sum
  _{\phtm{i\in F^L_j\cap F^M_j j=1}j=1}^{k} \sum _{\phtm{i\in
      F^L_j\cap F^M_j j=1}i\in F^L_j\cap F^M_j} a^L_i x^{\ast}_L(x_i) -
  2\epsilon\\
  & \geq & \frac{1}{2}\sum _{j\in J} w_j - 2\epsilon\\
  & \geq & \frac{cK'(\delta)}{8} -
  3\epsilon\geq\frac{c}{16}K'(\delta).
\end{eqnarray*}
\end{proof}
The discrete nature of the sufficient condition of
Proposition~\ref{prop:pure-matching} makes it very attractive: it
reduces Problem~\ref{problem:near-unc} to a combinatorial, Ramsey type
problem. The conclusion in this condition is about ``matching'' the
part of the combinatorial data of $L$ and $M$ that comes from the
discretization of the coefficients of $x_L$ and $x_M$, and it
``ignores'' the dual coefficients. We will now study the entire
combinatorial data as an abstract object (\ie we forget about the
underlying Banach space). This leads to the introduction of
resolutions. We will use them to discuss the possibility of a negative
answer to Problem~\ref{problem:near-unc}. To conclude this section we
shall produce an example to show that $\sup_{\delta>0}K'(\delta)$ is
strictly greater than~$1$ (recall that if $(x_i)$ is a normalized,
weakly null sequence with spreading model not equivalent to the unit
vector basis of $\mc_0$, then for any $\epsilon\kge 0$ there is a
subsequence $(y_i)$ of $(x_i)$ such that $K'((y_i),\delta)\kle
1\kplus\epsilon$).

Let $k\kin\bn$. A \emph{$k$-pattern} is a finite sequence in the set
$\{1,\ldots,k\}$ (the numbers $1,\ldots,k$ will be called
\emph{colours}). A \emph{$k$-resolution} is a pair
$r\keq((c_i)_{i=1}^n,(\alpha_i)_{i=1}^n)$, where $(c_i)_{i=1}^n$ is a
$k$-pattern, and $(\alpha_i)_{i=1}^n$ are positive, real numbers. When
we work with a fixed $k$ we shall simply say pattern and resolution,
respectively.

Let $r$ be a $k$-resolution. The \emph{weight of colour $j$ in $r$} is
\[
w_j(r)=\sum_{i:\,c_i=j} \alpha _i,\qquad j=1,\ldots,k,
\]
and the \emph{weight of $r$} is $w(r)\keq\sum_{j=1}^kw_j(r)$. A pair
$(x,x^{\ast})$ of elements of $\mc_{00}$ \emph{has resolution}~$r$ (or
$(x,x^{\ast})$ is a \emph{representation of} $r$) if the non-zero
co-ordinates of $x$ are $(2^{-c_i})_{i=1}^n$ in this order, and the
non-zero co-ordinates of $x^{\ast}$ are $(2^{c_i}\alpha_i)_{i=1}^n$ in
this order, and moreover $x$ and $x^{\ast}$ have the same support. In
other words, we have $x\keq\sum_{i=1}^n 2^{-c_i}e_{l_i}$ and
$x^{\ast}\keq\sum_{i=1}^n 2^{c_i}\alpha_ie_{l_i}$  for some $1\kleq
l_1\kle \ldots\kle l_n$. Note that
$x^{\ast}(x)\keq\sum_{i=1}^n\alpha_i\keq w(r)$.

Given $k\kin\bn$ and non-negative, real numbers $w_1,\ldots,w_k$
(called \emph{weights}) with $\sum_{j=1}^k w_j\keq 1$, we let
$\mathcal R\keq\mathcal R(w_1,\ldots,w_k)$ be the class of all
$k$-resolutions $r$ with $w_j(r)\keq w_j$ for each $j\keq 1,\ldots,
k$. The necessary and sufficient condition of
Proposition~\ref{prop:matching} motivates the following
definition. Given $r,s\kin\mathcal R$ we let
\[
[r,s]=\max x^{\ast}(y),
\]
where the maximum is over all pairs $(x,x^{\ast})$ and $(y,y^{\ast})$ of
elements of $\mc_{00}$ that have resolutions $r$ and $s$,
respectively. We also let $\langle r,s\rangle\keq\max \big\{
[r,s],[s,r]\big\}$. Note that $[r,s]\kleq\sum_{j=1}^k2^{j-1}w_j$
for all $r,s\kin\mathcal R$.

Given $k$-patterns $c\keq(c_i)_{i=1}^m$ and $d\keq(d_i)_{i=1}^n$, we
write $c\ksubset d$ if there exist $1\kleq l_1\kle\ldots\kle l_m\kleq
n$ such that $c_i\keq d_{l_i}$ for $i\keq 1,\ldots,m$. Observe that if
$r\keq(c,\alpha)$ and $s\keq (d,\beta)$ are elements of $\mathcal R$
and $c\ksubset d$, then $[r,s]\kgeq 1$. More generally, if we can find
representations $(x,x^{\ast})$ and $(y,y^{\ast})$ of $r$ and $s$,
respectively,
and a set $J\ksubset\{1,\ldots,k\}$ so that $\{i\kin\bn:\,x_i\keq
2^{-j}\}\ksubset\{i\kin\bn:\,y_i\keq 2^{-j}\}$ for each $j\kin J$,
then we have $[r,s]\kgeq x^{\ast}(y)\kgeq \sum_{j\in J}w_j$ (this
observation is motivated by
Proposition~\ref{prop:pure-matching}). Since for any
$j\kin\{1,\ldots,k\}$ we can find representations $(x,x^{\ast})$ and
$(y,y^{\ast})$ such that the sets $\{i\kin\bn:\,x_i\keq
2^{-j}\}$ and $\{i\kin\bn:\,y_i\keq 2^{-j}\}$ are comparable, we have
$\langle r,s\rangle\kgeq\max w_j\kgeq 1/k$ for all $r,s\kin\mathcal
R$.

Given $r,s\kin\mathcal R$ and $\eta\kin (0,1)$, we say that $r$ and
$s$ are $\eta$-\emph{orthogonal}, in symbols $r\perp_{\eta}s$, if
$\langle r,s\rangle\kle\eta$. Note that this can only happen for
$\eta\kge 1/k$. Roughly speaking, if one could find for each
$k\kin\bn$ an infinite set of pairwise $\eta(k)$-orthogonal
resolutions with $\eta(k)\to 0$ as $k\to\infty$, then one
could `code' an example in a way reminiscent of the Maurey-Rosenthal
construction~\cite{MR} to show that $\sup_{\delta >0}
L(\delta)\keq\infty$, where $L$ is the function given in
Definition~\ref{defn:constants}. We sketch this next.
\begin{ex}
  \label{ex:maurey-rosenthal}
  Let $k\kin\bn,\ \eta\keq\eta(k)\kin(0,1)$ and $C\keq
  C(k)\kgeq1$. Assume that we can find weights $w_1,\ldots,w_k$ and a
  sequence $(r_i)$ in $\mathcal R\keq\mathcal R(w_1,\ldots,w_k)$ so
  that $\langle r_i,r_j\rangle \kle \eta$ whenever $i\kneq j$, and
  $\langle r_i,r_i\rangle\kleq C$ for all $i\kin\bn$. Assume also
  that if $r_i\keq (c^{(i)},\alpha^{(i)})$, then
  $\max_j2^{c^{(i)}_j}\alpha^{(i)}_j\kleq 1$ for all $i\kin\bn$, and
  $\max_j2^{c^{(i)}_j}\alpha^{(i)}_j\to 0$ as $i\to\infty$. (note that
  this is not a serious assumption: the resolutions in a large family
  of pairwise orthogonal elements of $\mathcal R$ are necessarily
  ``flat'' --- \cf proof of
  Proposition~\ref{prop:low-resolution-index}). We will now show that
  $L(2^{-k})\kgeq 1/(2C\kplus 6)\eta$. In particular, if
  $(C(k))_{k=1}^{\infty}$ is bounded and $\eta(k)\to 0$ as
  $k\to\infty$, then this solves Problem~\ref{problem:near-unc} in the
  negative.

  Let $Q$ be the set of all representations of the resolutions $r_i,\
  i\kin\bn$. Let us fix an injective function $\phi$ (the coding
  function) that maps finite sequences of elements of $Q$ to positive
  integers. A sequence $(x_j,x_j^{\ast})_{j=1}^k$ of pairs of elements of
  $\mc_{00}$ is called a \emph{special sequence} if there exist
  positive integers $l_j$ for $j\keq 1,\ldots,k$ such that the
  following hold.
  \begin{myeqnarray}
    x^{\ast}_1<\ldots<x^{\ast}_k,\\
    (x_j,x^{\ast}_j)\ \text{has resolution}\ r_{l_j}\ \text{for}\ j\keq
    1,\ldots,k,\\
    l_j\keq\phi\big((x_1,x^{\ast}_1),\ldots,(x_{j-1},x^{\ast}_{j-1})\big)\
    \text{for}\ j\keq 1,\ldots,k.
  \end{myeqnarray}
  We then call the sum $\sum_{j=1}^kx^{\ast}_j$ a \emph{special
  functional}. Let $\mathcal F$ be the set of all special functionals,
  and let us define a norm on $\mc_{00}$ by letting
  \[
  \norm{x}=\norm{x}_{\ell_{\infty}}\vee \sup\big\{
  \abs{x^{\ast}(Ex)}:\,x^{\ast}\kin\mathcal F,\ E\kin\mathcal I\big\}.
  \]
  Here $\mathcal I$ denotes the set of intervals of positive integers,
  and $Ex$ is the projection of $x$ onto $E$. Let $X$ be the
  completion of $\mc_{00}$ in this norm. Then $(e_i)$ is a normalized,
  bimonotone, weakly null basis of $X$. Let $M\kin\vegtelen{\bn}$. One
  can clearly choose a special sequence $(x_j,x^{\ast}_j)_{j=1}^k$ such
  that $\supp{x_j}\ksubset M$ for each $j\keq 1,\ldots ,k$. Using the
  injectivity of $\phi$ and the orthogonality of the resolutions
  $r_i$, it is not difficult to show that $\Bnorm{\sum_{j=1}^k
  (-1)^jx_j}\kleq 1\kplus C\kplus 2k\eta\kleq (C\kplus 3)k\eta$,
  whereas
  \[
  \Bnorm{\sum_{j\,\text{odd}}x_j} +
  \Bnorm{\sum_{j\,\text{even}}x_j}\geq
  \sum_{j=1}^kx^{\ast}_j\Big(\sum_{j=1}^k x_j\Big)=k.
  \]
  This shows that $L((e_i)_{i\in M},2^{-k})\kgeq 1/(2C\kplus 6)\eta$.
\end{ex}
Our next result together with an earlier observation shows that a
Maurey-Rosenthal-type example as described above is far from
possible. Indeed, it shows that for all $k\kin\bn$ and for all weights
$w_1,\ldots,w_k$, any infinite subset $\mathcal S$ of $\mathcal
R(w_1,\ldots,w_k)$ contains a further infinite subset $\mathcal S'$
such that $\langle r,s\rangle\kgeq 1$ for all $r,s\kin\mathcal S'$.
\begin{prop}
Let $k\kin\bn$. Given $k$-patterns $c^{(i)},\ i\kin\bn$, there exist
$1\kleq l_1\kle l_2\kle\ldots$ such that $c^{(l_i)}\ksubset
c^{(l_{i+1})}$ for all $i\kin\bn$.
\end{prop}
\begin{proof}
  We apply induction on $k$. When $k\keq 1$ the result is trivial. Now
  assume that $k\kge 1$. For each $i\kin\bn$ we can write
  \[
  c^{(i)}=
  (c^{(i,1)},c^{(i)}_1,c^{(i,2)},c^{(i)}_2\ldots,
  c^{(i,m_i)},c^{(i)}_{m_i},c^{(i,m_i+1)}),
  \]
  where $m_i$ is a non-negative integer, $c^{(i)}_j$ is a single
  colour (\ie an element of $\{1,\ldots,k\}$) and $c^{(i,j)}$ is a
  $k$-pattern using exactly the $k\kminus 1$ colours
  $\{1,\ldots,k\}\ksetminus \{c^{(i)}_j\}$ for $1\kleq j\kleq m_i$, and
  finally $c^{(i,m_i+1)}$ is a pattern (possibly of length zero) using
  strictly less than $k$ colours. To see this simply trace the pattern
  $c^{(i)}$ from left to right and stop every time you have seen all
  $k$ colours.

  We consider two cases. In the first case $\sup_i
  m_i\keq\infty$. Let $\lambda_i$ be the length of $c^{(i)}$ for each
  $i\kin\bn$. We can find $1\kleq l_1\kle l_2\kle\ldots$ such
  that $m_{l_{i+1}}\kge\lambda_{l_i}$ for all $i\kin\bn$. Then for
  each $i\kin\bn$ the pattern $c^{(l_{i+1})}$ is the concatenation of
  more than $\lambda_{l_i}$ patterns each using all $k$ colours, from
  which $c^{(l_i)}\ksubset c^{(l_{i+1})}$ is clear.

  In the second case the sequence $(m_i)_{i=1}^{\infty}$ is
  bounded. Then after passing to a subsequence we may assume that for
  all $i\kin\bn$ we have $m_i\keq m,\ c^{(i)}_j\keq c_j$ for $1\kleq
  j\kleq m$, and $c^{(i,m+1)}$ uses exactly the colours from a proper
  subset $S$ of $\{1,\ldots,k\}$. Then by the induction hypothesis we
  find $1\kleq l_1\kle l_2\kle\ldots$ such that $c^{(l_i,j)}\ksubset
  c^{(l_{i+1},j)}$ for all $i\kin\bn$ and for each $j\keq
  1,\ldots,m\kplus 1$. It follows that $c^{(l_i)}\ksubset
  c^{(l_{i+1})}$ for all $i\kin\bn$.
\end{proof}
Having seen that there are no pairwise $\eta$-orthogonal, infinite
sets of resolutions for any $\eta\kin (0,1]$, we now introduce
so-called \emph{Rademacher resolutions} that form arbitrarily large,
finite sets of pairwise $\eta$-orthogonal resolutions for $\eta$ of
the order $1/\sqrt{k}$. This kills any hope of obtaining a positive
answer to Problem~\ref{problem:near-unc} by proving a version of
the Matching Lemma that allows us to find for any $N\kin\bn$, infinite
sets $L_1,\ldots,L_N$ whose combinatorial data match in a suitable
way.

Let $k\kgeq 4$ and $k_0\keq\lfloor\sqrt{k}\rfloor$. Set
$w_{jk_0}\keq 1/k_0$ for each $j\keq 1,\ldots,k_0$, and let $w_j\keq
0$ when $j$ is not a multiple of $k_0$. We will now consider certain
special elements of $\mathcal R\keq\mathcal R(w_1,\ldots,w_k)$. Fix
positive integers $n_1\kle\ldots\kle n_{k_0}$ satisfying
\begin{equation}
  \label{eq:rademacher-ris}
  \sum_{1\leq j<j'\leq k_0} \frac{n_j}{n_{j'}} <2^{-k}.
\end{equation}
For $n\in\bn$ we denote by $R_n$ the resolution $(c,\alpha)$, where
\[
c=(k_0,\ldots,k_0,2k_0,\ldots,2k_0,\ldots,k_0^2,\ldots,k_0^2),
\]
where colour $jk_0$ appears $n n_j$ times, and $\alpha_i\keq
1/nn_jk_0$ whenever $c_i\keq jk_0$ (\ie we distribute each
weight uniformly over the corresponding colour). We will use the
following notation: given $m\kin\bn$ and a resolution $r\keq
(c,\alpha)$ we write $(r,\ldots,r)_m$ for the resolution
$s\keq(d,\beta)$, where $d\keq (c,\ldots,c)$ with $c$ repeated $m$
times, and $\beta\keq (\alpha/m,\ldots,\alpha/m)$ with $\alpha/m$ also
repeated $m$ times. Note that if $r$ belongs to $\mathcal R$, then so
does $(r,\ldots,r)_m$ (indeed, this is true for any choice of weights
$w_1,\ldots,w_k$). Now given $l,n\kin\bn$, we define the
\emph{Rademacher} $R_{n,l}$ to be the resolution
$(R_n,\ldots,R_n)_{k_0^{l-1}}$. Note that $R_{n,l}\kin\mathcal R$ for
all $l,n\kin\bn$.
\begin{prop}
  \label{prop:rademacher}
  For all $m,n\kin\bn$, the Rademachers $R_{nk_0^{m-l},l}$, $l\keq
  1,\ldots,m$, are pairwise $5/k_0$-orthogonal. Moreover, $\langle
  R_{nk_0^{m-l},l},R_{nk_0^{m-l},l}\rangle\kleq 1\kplus \frac{2}{k_0}$
  for each $l\keq 1,\ldots,m$.
\end{prop}
\begin{proof}
  Fix $l,l'\kin\{1,\ldots,m\}$, let $r\keq R_{nk_0^{m-l},l}\keq
  (c,\alpha)$ and $s\keq R_{nk_0^{m-l'},l'}\keq (d,\beta)$. Choose
  representatives $(x,x^{\ast})$ and $(y,y^{\ast})$ of $r$ and $s$,
  respectively, so that
  \[
  [r,s]=x^{\ast}(y)=\sum_{i\in \supp{x}\cap\supp{y}} x^{\ast}_iy_i.
  \]
  Note that each term $x^{\ast}_iy_i$ is equal to $2^{c_u}\alpha
  _u2^{-d_v}$ for some $u$ and $v$. Let $S_1$ (respectively, $S_2$ and
  $S_3$) be the set of all $i\kin\bn$ for which $c_u\kle d_v$
  (respectively, $c_u\kge d_v$ and $c_u\keq d_v$). It is clear that
  \[
  \sum_{i\in S_1} x^{\ast}_iy_i\leq
  2^{-k_0}\sum_u\alpha_u=2^{-k_0}.
  \]
  For each $i\kin S_2$ there exist $1\kleq j\kle j'\kleq k_0$ such
  that $x^{\ast}_i\keq 2^{c_u}\alpha_u$ and $y_i\keq 2^{-d_v}$, where
  $c_u\keq j'k_0,\ \alpha_u\keq 1/nn_{j'}k_0^m$ and $d_v\keq
  jk_0$. Moreover, colour $jk_0$ occurs $nn_jk_0^{m-1}$ times in
  $s$. It follows that
  \[
  \sum_{i\in S_2} x^{\ast}_iy_i\leq \sum_{1\leq j<j'\leq
  k_0}2^{(j'-j)k_0}
  \frac{1}{nn_{j'}k_0^m} nn_jk_0^{m-1} <\frac{1}{k_0},
  \]
  by the choice of $n_1,\ldots,n_{k_0}$ \eqref{eq:rademacher-ris}. So
  the only significant contribution to $[r,s]$ comes from the set
  $S_3$ of co-ordinates, \ie where the colours match. Here we always
  have the trivial estimate
  \[
  \sum_{i\in S_3} x^{\ast}_iy_i\leq \sum_u\alpha_u=w(r)=1.
  \]
  In particular, when $l\keq l'$ this gives $\langle r,s\rangle\kleq
  1\kplus \frac{2}{k_0}$, as required. Note that for $i\kin S_3$ we
  have $x_i\keq y_i\keq 2^{-jk_0}$ and $x^{\ast}_i\keq y^{\ast}_i\keq
  2^{jk_0}/nn_jk_0^m$ for some $j\kin\{1,\ldots,k_0\}$. In particular
  $x^{\ast}_iy_i\keq y^{\ast}_ix_i$, so when $l\kneq l'$ we may
  without loss of
  generality assume that $l\kle l'$. Recall that for each
  $j\kin\{1,\ldots,k_0\}$ colour $jk_0$ in $r$ comes in $k_0^{l-1}$
  blocks, each block having length $nn_jk_0^{m-l}$. Consider such a
  block $B$, and suppose that a $\delta$-proportion of the block
  corresponds to co-ordinates $i$ of $x$ that belong to $S_3$. The
  corresponding co-ordinates of $y$ in turn correspond to
  colour-$jk_0$ bits of $s$. Since $s$ is made up of $k_0^{l'-1}$
  copies of $R_{nk_0^{m-l'}}$ and since colour $jk_0$ appears
  $nn_jk_0^{m-l'}$ times in each copy, the number of copies used up in
  this matching is at least
  \begin{equation}
    \label{eq:rademacher;block}
    \frac{\delta nn_jk_0^{m-l}}{nn_jk_0^{m-l'}}=\delta k_0^{l'-l},
  \end{equation}
  which is strictly greater than~$1$ if $\delta\kge 1/k_0$. Also the
  contribution of a $\delta$-proportion of block $B$ to $\sum_{i\in
  S_3}x^{\ast}_iy_i$ is
  \begin{equation}
    \label{eq:rademacher;weight}
    \delta nn_jk_0^{m-l}\frac{1}{nn_jk_0^m}=\delta k_0^{-l}.
  \end{equation}
  Let $\Delta$ be the sum of the $\delta$'s that are greater than
  $1/k_0$ over all colour-$jk_0$ blocks $B$ of $r$ and over all
  $j\kin\{1,\ldots,k_0\}$. It follows from~\eqref{eq:rademacher;block}
  that $\Delta k_0^{l'-l}$ is at most $2k_0^{l'-1}$, since the number
  of copies of $R_{nk_0^{m-l'}}$ that make up $s$ is $k_0^{l'-1}$ and
  each copy is counted at most twice. Hence
  from~\eqref{eq:rademacher;weight} we obtain the estimate
  \[
  \sum_{i\in S_3}x^{\ast}_iy_i\leq \Delta k_0^{-l}+\frac{1}{k_0}w(r)\leq
  \frac{3}{k_0}.
  \]
  This finally shows that $[r,s]\kleq 5/k_0$, as required.
\end{proof}
\begin{rems}
  1.~We observed earlier that for any $r,s\kin\mathcal R$ we have
  $[r,s]\kgeq \max w_j$, which in the above situation is
  $1/k_0$. Moreover, we always have $\langle r,r\rangle\kgeq 1$ for
  all $r\kin\mathcal R$. So the measure of orthogonality we achieve is
  essentially best possible.\\
  2.~In Example~\ref{ex:maurey-rosenthal} we required the resolutions
  in the pairwise orthogonal family to be `flat'. Note that this holds
  for the Rademachers. Given $m,n\kin\bn$ and $l\kin\{1,\ldots,m\}$,
  if $R_{nk_0^{m-l},l}\keq (c,\alpha)$, then $\max _i
  2^{c_i}\alpha_i\kleq 2^k/nn_1k_0^m\to 0$ as $m\to\infty$.
\end{rems}

It is possible to measure, for each $\eta\kin (0,1)$, the complexity
of the family of finite sets of pairwise $\eta$-orthogonal resolutions
by introducing a suitable ordinal index. We shall not do that, but
simply comment that the above result would then say that for $\eta\kge
5/\sqrt{k}$ and under the assumption that we only use colours that are
multiples of $\sqrt{k}$ and carry equal weights, this complexity is at
least $\omega$. Whereas our next result shows that the complexity never
exceeds $\omega$ (and this holds for general weights). So in some
sense the set of resolutions has just enough complexity to allow the
possibility of a negative asnwer to Problem~\ref{problem:near-unc}.
\begin{prop}
\label{prop:low-resolution-index}
Assume that $k\kin\bn$ and $w_1,\ldots,w_k$ are arbitrary weights. Let
$\mathcal R\keq\mathcal R(w_1,\ldots,w_k)$ and $\eta\kin (0,1/4)$. For
all $r\kin\mathcal R$ there exists $n\kin\bn$ so that whenever
$s_1,\ldots,s_n\kin\mathcal R$ are pairwise $\eta$-orthogonal, we have
$[r,s_i]\kgeq 1/2$ for each $i\keq 1,\ldots,n$.
\end{prop}
\begin{proof}
  Choose $j_0$ and $j_1$ minimal so that
  \[
  \sum_{j=1}^{j_0}w_j\geq 1/4\qquad and\qquad
  \sum_{j=1}^{j_1}w_j\geq 1/2.
  \]
  We then have
  \[
  \sum_{j=j_0}^{j_1}w_j\geq 1/4\qquad and\qquad
  \sum_{j=j_1}^kw_j\geq 1/2.
  \]
  Now assume the result is false. Then there exists $r\kin\mathcal R$
  such that for all $n\kin\bn$ we have $\mathcal R_n\ksubset\mathcal
  R$ and $t_n\kin\mathcal R_n$ such that $\babs{\mathcal R_n}\kgeq n$,
  $\langle t,t'\rangle\kle\eta$ for all $t,t'\kin\mathcal R_n$ with
  $t\kneq t'$ and $[r,t_n]\kle 1/2$. We now verify two claims.

  First observe that for each $n\kin\bn$ the number of co-ordinates of
  $t_n$ of colours $1,\ldots,j_1$ is at most the length of
  $r$. Indeed, otherwise we can choose representatives $(x,x^{\ast})$ of
  $r$ and $(y,y^{\ast})$ of $t_n$ so that whenever $x_i\keq 2^{-j}$ for
  some $j\kgeq j_1$, then $y_i\keq 2^{-j'}$ for some $j'\kleq j_1$,
  and this would give $[r,t_n]\kgeq w_{j_1}\kplus\dots\kplus w_k\kgeq
  1/2$.

  Secondly, we claim that for all $n\kin\bn$ and for all
  $t\kin\mathcal R_n$ the number of co-ordinates of $t$ of colours
  $1,\ldots,j_0$ is at most the length of $r$. Otherwise by the first
  claim we can find representatives $(x,x^{\ast})$ of $t_n$ and
  $(y,y^{\ast})$
  of $t$ so that whenever $x_i\keq 2^{-j}$ for some $j_0\kleq j\kleq
  j_1$, then $y_i\keq 2^{-j'}$ for some $j'\kleq j_0$, and this would
  give $[t_n,t]\kgeq w_{j_0}\kplus\dots\kplus w_{j_1}\kgeq 1/4$.

  Now by simple pigeonhole principle, if $n$ is greater than the
  number of patterns of length at most the length of $r$ in
  colours $1,\ldots,j_0$, then there exist distinct $t,t'\kin\mathcal
  R_n$ so that the patterns in $t$ and $t'$ formed by the colours
  $1,\ldots,j_0$ are identical. It follows that there exist
  representatives $(x,x^{\ast})$ of $t$ and $(y,y^{\ast})$ of $t'$ so that
  $\{i\kin\bn:\,x_i\keq 2^{-j}\}\keq\{i\kin\bn:\,y_i\keq 2^{-j}\}$ for
  each $j\keq 1,\ldots,j_0$, and hence we obtain the contradiction
  $\langle t,t'\rangle\kgeq w_1\kplus\dots\kplus w_{j_0}\kgeq 1/4$.
\end{proof}
We conclude by constructing a relatively simple example using
Rademachers to show that $\sup_{\delta>0} K'(\delta)\kgeq 5/4$.
\begin{ex}
  \label{ex:elton>1}
  Let $\epsilon\kin(0,1)$. Fix positive integers $n_1\kle n_2$ and $K$
  such that
  \[
  \frac{n_1}{2n_2}+2^{-K}<\epsilon\qquad\text{and}\qquad
  \frac{2n_1+n_2}{n_12^K}<1.
  \]
  For an infinite subset $M\keq\{m_1\kle m_2\kle \ldots\}$ of $\bn$
  set $n_M\keq (n_1\kplus n_2)2^{Km_2-1}$ and let
  \[
  E_M=\{m_3, m_4, \ldots,m_{n_M+2}\}.
  \]
  Now write $E_M$ as a union
  \[
  E_M=\bigcup_{j=1}^{2^{Km_1-1}}I^M_j\ \cup\
  \bigcup_{j=1}^{2^{Km_1-1}}J^M_j,
  \]
  where $I^M_1\kle J^M_1\kle I^M_2\kle J^M_2\kle\ldots\kle
  I^M_{2^{Km_1-1}}\kle J^M_{2^{Km_1-1}}$ and $\abs{I^M_j}\keq n_1
  2^{Km_2-Km_1}$ and $\abs{J^M_j}\keq n_2 2^{Km_2-Km_1}$ for each
  $j\keq 1,\ldots,2^{Km_1-1}$. Finally, set
  \[
  E^M_1=\bigcup_{j=1}^{2^{Km_1-1}}I^M_j\qquad \text{and}\qquad
  E^M_2=\bigcup_{j=1}^{2^{Km_1-1}}J^M_j,
  \]
  so we have $\abs{E^M_1}\keq n_1 2^{Km_2-1}$ and $\abs{E^M_2}\keq n_2
  2^{Km_2-1}$. Note that if we let $c_i\keq 2$ whenever $m_{i+2}\kin
  E^M_1$ and $c_i\keq 4$ whenever $m_{i+2}\kin E^M_2$, then $(c_i)$ is
  the pattern of the Rademacher resolution $R_{2^{Km_2-Km_1},Km_1}$ as
  defined preceding Proposition~\ref{prop:rademacher} when $k\keq 4$.

  We shall denote by $\bi_F$ the indicator function of a set
  $F\ksubset\bn$, which is also the element $\sum_{i\in F}e_i$ of
  $\mc_{00}$. Given $M\keq\{m_1\kle m_2\kle
  \ldots\}\kin\vegtelen{\bn}$, let
  \[
  \begin{array}{lcrcrcrcr}
    x_M &=& -\frac{1}{2}e_{m_1} &+& \frac{1}{2}e_{m_2} &+&
    \frac{1}{2}\bi_{E^M_1} &+& \frac{1}{4}\bi_{E^M_2},\\[6pt]
    x^+_M &=& &&\frac{1}{2}e_{m_2} &+& \frac{1}{2}\bi_{E^M_1} &+&
    \frac{1}{4}\bi_{E^M_2},\\[6pt]
    x^{\ast}_M &=& \frac{1}{2}e_{m_1} &+& e_{m_2} &+&
    \frac{1}{\abs{E^M_1}}\bi_{E^M_1} &+&
    \frac{1}{\abs{E^M_2}}\bi_{E^M_2}.
  \end{array}
  \]
  Define a norm on $\mc_{00}$ by setting
  \[
  \norm{x}=\norm{x}_{\ell_{\infty}}\vee\sup\big\{\abs{x^{\ast}_M(Ex)}:
  \,M\kin \vegtelen{\bn},\ E\kin\mathcal I\big\}
  \]
  for each $x\kin\mc_{00}$. Here $\mathcal I$ denotes the set of
  initial segments of $\bn$. Let $X$ be the completion of
  $(\mc_{00},\norm{\cdot})$. It is easy to verify that $(e_i)$ is a
  normalized, weakly null, monotone basis of $X$. We are going to show
  that for any subsequence $(f_i)$ of $(e_i)$ we have
  $K((f_i),1/4)\kgeq 5/4(1\kplus\epsilon)$. Since $\epsilon$ was
  arbitrary, this shows that $K(1/4)\kgeq 5/4$.

  Fix $M\keq\{m_1\kle m_2\kle\ldots\}\kin\vegtelen{\bn}$. On the one
  hand we have
  \[
  \norm{x^+_M}\geq x^{\ast}_M(x^+_M)=\frac{5}{4}.
  \]
  On the other hand, we are going to show that $\norm{x_M}\kleq
  1\kplus\epsilon$. So let us fix $L\keq\{l_1\kle l_2\kle
  \ldots\}\kin\vegtelen{\bn}$. We need to estimate $x^{\ast}_L(Ex_M)$ for
  any $E\kin\mathcal I$. This is always at least $-\frac{1}{2}$. To
  get an upper bound, we may clearly assume that $\supp{x_M}\ksubset
  E$. We now split into four cases. The first three of these use only
  the trivial estimate
  \[
  x^{\ast}(y)=\sum_ix^{\ast}_iy_i\kleq \min
  \big\{\norm{x^{\ast}}_{\ell_{\infty}}\cdot\norm{y}_{\ell_1},\
  \norm{x^{\ast}}_{\ell_1}\cdot\norm{y}_{\ell_{\infty}}\big\}
  \]
  for any $x^{\ast},y\kin\mc_{00}$.

  \noindent
  \textit{Case~1.} If $l_1\keq m_1$ and $l_2\keq m_2$, then we have
  \begin{eqnarray*}
    x^{\ast}_L(x_M) &=& -\frac{1}{4} +\frac{1}{2} +\frac{1}{\abs{E^L_1}}
    \bi_{E^L_1} \Big(\frac{1}{2}\bi_{E^M_1} +
    \frac{1}{4}\bi_{E^M_2}\Big)\\[6pt]
    & & +\frac{1}{\abs{E^L_2}} \bi_{E^L_2} \Big(\frac{1}{2}\bi_{E^M_1}
    + \frac{1}{4}\bi_{E^M_2}\Big)\\[6pt]
    &\leq& -\frac{1}{4} +\frac{1}{2} + 1\cdot\frac{1}{2} +
    \frac{1}{n_22^{Km_2-1}} \Big(\frac{1}{2}n_1 2^{Km_2-1} +
    \frac{1}{4}n_2 2^{Km_2-1}\Big)\\[6pt]
    &=& 1+\frac{n_1}{2n_2} < 1+\epsilon.
  \end{eqnarray*}
  \textit{Case~2.} If $l_2\kle m_2$, then we have
  \begin{eqnarray*}
    x^{\ast}_L(x_M) &=& x^{\ast}_L\Big(-\frac{1}{2}e_{m_1}\Big)\\[6pt]
    & & + \bigg(\frac{1}{\abs{E^L_1}}\bi_{E^L_1} +
    \frac{1}{\abs{E^L_2}}\bi_{E^L_2}\bigg) \bigg(
    \frac{1}{2}e_{m_2}+\frac{1}{2}\bi_{E^M_1} +
    \frac{1}{4}\bi_{E^M_2}\bigg)\\[6pt]
    &\leq& 0+2\cdot\frac{1}{2}=1.
  \end{eqnarray*}
  \textit{Case~3.} If $l_2\kge m_2$, then we have
  \begin{eqnarray*}
    x^{\ast}_L(x_M) &=& \Big( \frac{1}{2}e_{l_1}+e_{l_2}\Big)(x_M)\\[6pt]
    & & + \bigg(\frac{1}{\abs{E^L_1}}\bi_{E^L_1} +
    \frac{1}{\abs{E^L_2}}\bi_{E^L_2}\bigg) \bigg(
    \frac{1}{2}\bi_{E^M_1} +
    \frac{1}{4}\bi_{E^M_2}\bigg)\\[6pt]
    &\leq& \frac{3}{2}\cdot\frac{1}{2} + \frac{1}{n_12^{Kl_2-1}} \Big(
    \frac{1}{2}n_12^{Km_2-1}+\frac{1}{4}n_22^{Km_2-1}\Big)\\[6pt]
    &=& \frac{3}{4} +\frac{2n_1+n_2}{4n_12^K}\leq 1.
  \end{eqnarray*}
  \textit{Case~4.} If $l_2\keq m_2$ and $l_1\kneq m_1$, then we have
  to use the structure of the Rademacher patterns to get an upper
  bound. The argument is along similar lines to the proof of
  Proposition~\ref{prop:rademacher}. First we have
  \begin{equation}
    \label{eq:large-constant;to-begin}
    x^{\ast}_L(x_M)=\frac{1}{2} +\bigg(\frac{1}{\abs{E^L_1}}\bi_{E^L_1} +
    \frac{1}{\abs{E^L_2}}\bi_{E^L_2}\bigg) \bigg(
    \frac{1}{2}\bi_{E^M_1} + \frac{1}{4}\bi_{E^M_2}\bigg),
  \end{equation}
  and
  \begin{equation}
    \label{eq:large-constant;big-on-small}
    \frac{1}{\abs{E^L_2}}\bi_{E^L_2}\bigg(
    \frac{1}{2}\bi_{E^M_1}\bigg) \leq
    \frac{1}{n_22^{Kl_2-1}}\cdot\frac{1}{2}\cdot n_12^{Km_2-1}=
    \frac{n_1}{2n_2}.
  \end{equation}
  Also, since $\abs{E^L_1\cap E^M_2}\kleq \abs{E^L_1}\kminus
  \abs{E^L_1\cap E^M_1}$, we have
  \begin{equation}
    \label{eq:large-constant;main}
    \begin{array}{rcl}
      \ds
      \bi_{E^L_1}\Big(\frac{1}{2}\bi_{E^M_1}+\frac{1}{4}\bi_{E^M_2}\Big)
      &=& \ds \frac{1}{2}\abs{E^L_1\cap E^M_1}
      +\frac{1}{4}\abs{E^L_1\cap E^M_2}\\[12pt]
      &\leq& \ds\frac{1}{4}\abs{E^L_1\cap E^M_1} +
      \frac{1}{4}\abs{E^L_1}.
    \end{array}
  \end{equation}
  Let us now assume that $l_1\kle m_1$. For each $j\keq
  1,\ldots,2^{Kl_1-1}$ set
  \[
  A_j=\big\{ i\kin\{1,\ldots,2^{Km_1-1}\}:\,I^M_i\cap I^L_j\kneq
  \emptyset\big\}.
  \]
  We now have
  \[
  \abs{E^L_1\cap E^M_1}=\sum_{j=1}^{2^{Kl_1-1}} \sum_{i\in A_j}
  \abs{I^L_j\cap I^M_i}\leq \sum_{j=1}^{2^{Kl_1-1}} \abs{A_j}
  n_12^{Km_2-Km_1}.
  \]
  Hence from~\eqref{eq:large-constant;main} we obtain
  \begin{equation}
    \label{eq:large-constant;colour-two}
    \frac{1}{\abs{E^L_1}}\bi_{E^L_1}\Big(\frac{1}{2}\bi_{E^M_1}+
    \frac{1}{4}\bi_{E^M_2}\Big) \leq \frac{1}{2}\cdot
    2^{-Km_1}\sum_{j=1}^{2^{Kl_1-1}} \abs{A_j}+\frac{1}{4}.
  \end{equation}
  Since $E^L_2\cap J^M_i\keq \emptyset$ whenever $\min A_j\kleq i\kle
  \max A_j$ for some $j\kin\{1,\ldots,2^{Kl_1-1}\}$, we have
  \[
  \abs{E^L_2\cap E^M_2}=\sum_{i=1}^{2^{Km_1-1}} \abs{E^L_2\cap
  J^M_i}\leq \bigg( 2^{Km_1-1}-\sum_{j=1}^{2^{Kl_1-1}} (\abs{A_j}-1)
  \bigg) n_22^{Km_2-Km_1}.
  \]
  It follows that
  \begin{equation}
    \label{eq:large-constant;colour-four}
    \frac{1}{\abs{E^L_2}}\bi_{E^L_2} \Big(\frac{1}{4}\bi_{E^M_2}\Big)
    \leq \frac{1}{4}-\frac{1}{2}\cdot 2^{-Km_1}\sum_{j=1}^{2^{Kl_1-1}}
    \abs{A_j} + 2^{Kl_1-Km_1}.
  \end{equation}
  Note that $2^{Kl_1-Km_1}\kleq 2^{-K}$ since we are assuming that
  $l_1\kle m_1$.
  Putting together \eqref{eq:large-constant;to-begin},
  \eqref{eq:large-constant;big-on-small},
  \eqref{eq:large-constant;colour-two} and
  \eqref{eq:large-constant;colour-four} we finally obtain
  \begin{equation}
    \label{eq:large-constant;final}
    x^{\ast}_L(x_M)\leq 1+\frac{n_1}{2n_2}+2^{-K}<1+\epsilon,
  \end{equation}
  as required. The case when $l_1\kge m_1$ is very similar. For each
  $j\keq 1,\ldots,2^{Km_1-1}$ set
  \[
  A_j=\big\{ i\kin\{1,\ldots,2^{Kl_1-1}\}:\,I^L_i\cap I^M_j\kneq
  \emptyset\big\}.
  \]
  We then proceed as before making the obvious changes in the various
  summations.
\end{ex}
\begin{rems}
  1.~Since $\norm{x^{\ast}_M}_{\ell_1}\kleq \frac{7}{2}$ for all
  $M\kin\vegtelen{\bn}$, the basis $(e_i)$ of $X$ is $7/2$-equivalent
  to the unit vector basis of $\mc_0$, yet no subsequence is
  $C$-unconditional for $C\kle 5/4(1\kplus\epsilon)$. So the above
  example also shows that $C(\delta)\kgeq 5/4$ whenever $\delta\kleq
  2/7$, where $C(\delta)$ is the constant introduced in
  Section~\ref{section:c_0-problem} in relation to the
  $\mc_0$-problem.\\
  2.~The basis $(e_i)$ of the space $X$ constructed above is also an
  example of a normalized, weakly null sequence that has no
  quasi-greedy basic subsequence with constant strictly less than
  $8/7$. To see this let $\alpha=2/3$ and let
  \[
  \begin{array}{lcrcrcrcl}
    y_M &=& -\alpha e_{m_1} &+& e_{m_2} &+& \bi_{E^M_1} &+&
    \frac{2}{3}\bi_{E^M_2},\\[6pt]
    y^+_M &=& && e_{m_2} &+& \bi_{E^M_1} &+&
    \frac{2}{3}\bi_{E^M_2}
  \end{array}
  \]
  for each $M\kin\vegtelen{\bn}$ (following the notation in the proof
  above). Given $\epsilon\kge 0$ we may choose the parameters $n_1,
  n_2$ and $K$ so that
  \begin{equation}
    \label{eq:quasy-greedy;large-constant}
    \frac{\norm{y^+_M}}{\norm{y_M}}> \frac{8}{7}-\epsilon
  \end{equation}
  for all $M\kin\vegtelen{\bn}$. This is proved by exactly the same
  calculation as in the proof above.
  
  Now if $\alpha\keq 2/3\kminus\eta$ for some $\eta\kge 0$,
  then~\eqref{eq:quasy-greedy;large-constant} still holds provided
  $\eta$ is sufficiently small. Then $y^+_M$ is the projection of
  $y_M$ onto the set of co-ordinates where the size of the coefficient
  is at least $2/3$. It follows that $(e_i)_{i\in M}$ is not
  quasi-greedy with constant $8/7\kminus\epsilon$ for any
  $M\kin\vegtelen{\bn}$.
\end{rems}

\vspace{1ex}

\noindent
Department of Mathematics, University of South Carolina, Columbia, SC
29208, USA.\\
\textit{E-mail adress:} \texttt{dilworth@math.sc.edu}\\

\noindent
Department of Mathematics, The University of Texas at Austin,
1~University Station C1200, Austin, TX 78712-0257, USA.\\
\textit{E-mail address:} \texttt{odell@math.utexas.edu}\\

\noindent
Department of Mathematics, Texas A\&M University, College Station,\\
TX 77843-3368, USA.\\
\textit{E-mail address:} \texttt{schlump@math.tamu.edu}\\

\noindent
Fitzwilliam College, Cambridge, CB3 0DG, England.\\
\textit{E-mail address:} \texttt{a.zsak@dpmms.cam.ac.uk}

\end{document}